\documentclass[12pt,a4paper]{article}
\usepackage[utf8]{inputenc}
\usepackage{amsmath}
\usepackage{amsfonts}
\usepackage{amssymb}
\usepackage{amsthm}
\usepackage{graphicx}

\usepackage{color, fullpage, float, epsf}
\usepackage{caption, subcaption}
\usepackage{enumitem} 
\usepackage{hyperref} 

\usepackage{titlesec}
\titlelabel{\thetitle. }

\usepackage{tikz}	
\usetikzlibrary{shapes}
\usetikzlibrary{calc}
\usetikzlibrary{decorations.pathreplacing}
\usetikzlibrary{hobby}
\usetikzlibrary{patterns}
\usetikzlibrary{arrows.meta}
\usepackage{tikz-3dplot}

\usepackage[left=2cm,right=2cm,top=2.5cm,bottom=2.5cm]{geometry}

\newcommand{\R}{\mathbb{R}}
\newcommand{\Sp}{\mathbb{S}}

\newcommand{\supp}{\mathrm{supp}}

\theoremstyle{definition}
\newtheorem{remark}{Remark}[section]

\newtheorem{problem}{Problem}

\providecommand{\keywords}[1]
{
	\small
	\noindent\textbf{{Keywords:}} #1\\
}

\providecommand{\classification}[1]
{
	\small
	\textbf{{AMS Mathematics Subjects Classification (2010):}} #1
}

\title{A geometric based preprocessing for weighted ray transforms with applications in SPECT}

\author{Fedor Goncharov\thanks{CMAP, CNRS, Ecole Polytechnique, Institut Polytechnique de Paris, 91128 Palaiseau, France  \newline \indent \indent email: fedor.goncharov.ol@gmail.com}}

\begin{document}
	
	\maketitle
	\abstract{In this work we investigate numerically the reconstruction approach proposed in Goncharov, Novikov, 2016, for weighted ray transforms (weighted Radon transforms along oriented straight lines) in 3D. In particular, the approach is based on a geometric reduction of the data modeled by weighted ray transforms  to new data modeled by weighted Radon transforms along two-dimensional planes in 3D.
		Such reduction could be seen as a preprocessing procedure which could be further completed by any preferred reconstruction algorithm. 	
		In a series of numerical tests on modelized and real SPECT (single photon emission computed tomography) data we demonstrate that such procedure can significantly reduce the impact of noise on reconstructions.\\}
	
	\\
	\classification{44A12, 35R30, 49N45}
	
	\section{Introduction}
	We consider weighted Radon transforms $P_W$ along oriented straight lines and $R_W$ along oriented planes in $\R^3$ which are defined by the formulas:
	
	\begin{align}
	\label{eq:wraytransform.def}
	&P_Wf(x,\theta) = \int\limits_{\gamma(x,\theta)}W(y,\theta)f(y)\,dy, \, 
	(x,\theta)\in T\Sp^2,\\
	\label{eq:wradontransform.def}
	&R_Wf(s,\theta) = \int\limits_{y\theta = s} W(y,\theta)f(y) \, dy, \, 
	(s,\theta) \in \R\times \Sp^2, 
	\end{align}
	where $f=f(x)$ is a test-function on $\R^3$, $W=W(x,\theta)$ is the weight. The manifold of oriented straight lines in $\R^3$ is identified with $T\Sp^2$ (tangent bundle of $\Sp^2$) via formulas:
	\begin{align}
	&T\Sp^2 = \{(x,\theta)\in \R^3\times \Sp^2 : x\theta = 0\}, \\
	\begin{split} 
	&(x,\theta)\in T\Sp^2 \text{ corresponds to ray } 
	\gamma(x,\theta) = \{y \in \R^3 : y = x + t\theta, \, t\in \R\}, \\
	&\text{where } \theta \text{ gives the orientation on } \gamma(x,\theta).
	\end{split}
	\end{align}
	Oriented planes in $\R^3$ also constitute a manifold which can be identified with $\R\times \Sp^2$ as follows:
	\begin{align}
	\begin{split}
	&(s,\theta) \in \R\times \Sp^2 \text{ corresponds to plane }
	\Sigma(s,\theta) = \{y \in \R^3 : y\theta = s\}, \\
	&\text{where } \theta \text{ gives the orientation on } \Sigma(s,\theta).
	\end{split}
	\end{align}
	Weighted Radon transforms defined in \eqref{eq:wraytransform.def}, \eqref{eq:wradontransform.def} and some of their generalizations are well-known in many domains of pure mathematics: in theory of groups (\cite{ggg1959geomhom}, \cite{ggv1962generalized}, \cite{harishchandra1958a}, \cite{helgason1965}, \cite{helgason1999}, \cite{ilmavirta2016}), 
	harmonic analysis (\cite{strichartz1982variations}, \cite{strichartz1991harmonic}), PDEs (\cite{beylkin1984inversion}, \cite{john1955book}), integral geometry (\cite{sharafutdinov2012}), microlocal analysis (\cite{beylkin1984inversion}, \cite{quinto2014}, \cite{quinto2018}) and can be also of self-interest (\cite{quinto1980measures}, \cite{frigyik2008xray}, \cite{boman2011localnoninj}, \cite{ilmavirta2019}, \cite{paternain2016geodesic}). At the same time, transformations $P_W$ (and a lot less $R_W$) are used as an important tool in many computerized tomographies (\cite{quinto1983rotinv}, \cite{natterer1986book}, \cite{miller1987newslant}, \cite{kuchment2014book}, \cite{novikov2002formula}, \cite{quinto2006}, \cite{deans2007book}, \cite{nguyen2009compton}, \cite{bal2011spect}, \cite{miqueles2011}).\\

	In this article we consider the two following problems:
	
	\begin{problem}\label{pr:problem-pw}
		Reconstruct $f$ from $P_Wf$ on rays parallel to some fixed 2-plane $\Sigma$.
	\end{problem}
	
	\begin{problem}\label{pr:problem-rw}
		Reconstruct $f$ from $R_W f$ on $\R\times \Sp^2$.
	\end{problem}
	In Problems~\ref{pr:problem-pw}, \ref{pr:problem-rw} we assume that weight $W$ is known apriori.
	
	Problem 1 naturally arises in different tomographies such as X-ray transmission tomography, PET, SPECT, fluorescence tomography and even in ultrasound-modulated optical tomography (UMOT); see, for example, \cite{natterer1986book}, \cite{kuchment2014book}, \cite{miqueles2011}, \cite{bocoum2019structured}. On the other hand, Problem 2 is less-known in tomographic applications, probably because the measured data often correspond to signals propagating along straight lines rather than planes. The only known to the author an example of using Problem~\ref{pr:problem-rw} directly in applications are migration-type algorithms in geophysics \cite{miller1987newslant}. 
	
	On the other hand, recently, in \cite{goncharov2016analog}, \cite{goncharov2017iterative} the authors
	had addressed the fact that Problem~\ref{pr:problem-pw} can be directly reduced to Problem~\ref{pr:problem-rw} using a simple geometrical formula (see formulas of \eqref{eq:reduction.main-formula} in Section~\ref{sect:mathmodel}). Such reduction makes Problem~\ref{pr:problem-rw} relevant for all types of tomographies, where the measured data are described by transformations $P_W$ from formula \eqref{eq:wraytransform.def}.
	In particular, in \cite{goncharov2016analog} it was noted that the aforementioned reduction can be used in SPECT to decrease the impact of strong Poisson noise on reconstructions. In this work we bring in for the first time a numerical evidence for this suggestion. 
	
	More precisely, by our experiments we show that reconstructions from reduced SPECT data $R_Wf$ are more stable in presence of noise than slice-by-slice reconstructions from $P_Wf$ as in Problem~\ref{pr:problem-pw}. For reconstructions from $P_Wf, R_Wf$ we consider two types of methods: (1) Chang-type methods \cite{chang1978correciton}, \cite{novikov2011chang}, \cite{goncharov2016analog} and (2) iterative Kunyansky-type methods  \cite{kunyansky1992generalized}, \cite{guillement2014finite}, \cite{goncharov2017iterative}. For each type we compare reconstructions from $P_Wf$ and $R_Wf$ ($R_Wf$ being the reduction of $P_Wf$) for different levels of noise, attenuation strengths and activity phantoms. It appears, that in all of the above cases with non-zero noise level the reconstructions from reduced data $R_Wf$ have significantly better quality than their respective reconstructions from $P_Wf$. For the proof of concept we also apply the new reconstruction approach on real SPECT data provided by \textit{Service Hospitalier Frédéric Joliot, CEA (Orsay)}. 
	Though in the latter case any comparison is impossible because the real isotope distribution is unknown, by this experiment we show that the new reconstruction approach can be applied in a realistic SPECT framework.

	The organization of this article is as follows. In Section~\ref{sect:mathmodel} we recall formulas from \cite{goncharov2016analog} which give the reduction of   Problem~\ref{pr:problem-pw} to Problem~\ref{pr:problem-rw}. In Section~\ref{sect:numericalexp-modelized} we explain the organization of  numerical experiment on modelized SPECT data and present corresponding reconstruction examples. In Section~\ref{sect:numericalexp-real} we test our new method on real SPECT data. Discussion of our results is given in Section~\ref{sect:discussion}.
	
	The method we consider can serve as a preprocessing procedure for tomographical data modeled by weighted ray transforms along straight lines. 
	The point is that more complicated than Chang-type or Kunyansky-type methods  could be further applied to perform reconstructions from such preprocessed data. For example, algorithms capable of correcting strong non-uniform attenuation and noise effects (for example \cite{hudson1994accelerated}, \cite{slambrouck2014reconstruction}, \cite{vslambrock2015pet}, \cite{filipovic2019pet}) could be tested on new data, which makes it an interesting topic for future research.

	\section{Preprocessing formulas}\label{sect:mathmodel}
\begin{figure}[H]
		\centering
		\scalebox{0.7}{\begin{tikzpicture}[use Hobby shortcut, closed=true]
		\draw [fill=gray!20] (-1,0) .. (-0.7,2.5) .. (1.5,2).. (3,3.0) .. (5,2.5).. (5,0.5) ..(3,-1);
		\draw[dashed] (3, -1) arc (-90:90:0.5 and 2);
		\draw[semithick] (3, 3) arc (90:270:0.2 and 2);
		
		\draw[dashed] (1.1, -1.16) arc (-90:90:0.5 and 1.59);
		\draw[semithick] (1.1, 2.02) arc (90:270:0.2 and 1.59);	
		
		\draw[dashed] (-0.5, -0.46) arc (-90:90:0.5 and 1.51);
		\draw[semithick] (-0.5, 2.53) arc (90:270:0.2 and 1.5);		
		
		\draw (-4,-2)--(-4, 3) -- (-2.5, 5) -- (-2.5, 0) -- (-4, -2);
		\fill (-3.0, 1.5) circle (1.2pt) node [right,scale=3] {};
		\draw[->, semithick] (-3.0, 1.5) -- (-2, 1.5);
		
		\node at (-3, 1.0) {\large $O$};
		\node at (-2, 1.9) {\large $\eta$};
		\node at (-1.0, 4.5) {\large $\Sigma_{\eta} = \Sigma(0,\eta)$};
		\node at (5, -1) {\large $\mathrm{supp} \, f$};
		
		\draw[->, semithick, dashed, red] (1.3, 3) -- (1, -2);
		\node at (1.6, -2.3) {\large $P_Wf(\gamma), \gamma \in T\Sp^{1}(Y)$};
		\node at (1.8, 2.8) {\large $f|_{Y}$};	
		\node at (0.6, 1.0) {\large $Y$};
		
\end{tikzpicture}}
		\caption{Reconstructions slice-by-slice}
		\label{fig:slice-by-slice-reconstructions}
	\end{figure}
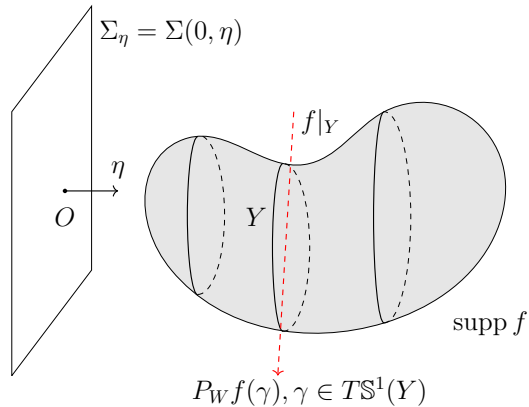

	In many tomographies data modeled by $P_Wf$ is often restricted to the subset of rays which are parallel to some fixed plane (see Figure~\ref{fig:slice-by-slice-reconstructions}). This situation occurs, for example, in X-ray tomography and in emission tomographies (e.g., in PET, SPECT), where a gantry of detectors moves along some fixed direction (along $\eta$ on Figure~\ref{fig:slice-by-slice-reconstructions}), and in each  two-dimensional slice $Y$ spectral data $P_Wf(\gamma), \gamma \in T\Sp^1(Y)$ are collected. 
	The standard way to reconstruct $f$ is to invert $P_Wf(\gamma), \, \gamma\in T\Sp^1(Y)$ on each slice $Y$ and then stack all slices together. In literature this is also  called slice-by-slice reconstruction approach which in this article we associate with solving Problem~\ref{pr:problem-pw}.
	
	In \cite{goncharov2016analog} authors proposed another approach which is based on a reduction of Problem~\ref{pr:problem-pw} to another inverse problem for weighted Radon transforms along oriented two-dimensional planes in $\R^3$. The reduction is performed via the following formulas:
	\begin{align}\label{eq:reduction.main-formula}
	\begin{split}
	&R_wf(s,\theta) = \int\limits_{-\infty}^{+\infty}P_{W}f(s\theta + \tau [\theta, \alpha_\eta(\theta)], \, \alpha_\eta(\theta)) \, d\tau, \, (s,\theta)\in \R\times \Sp^2 \\
	& w(x,\theta) = W(x,\alpha_\eta(\theta)), \, \alpha_\eta(\theta) = \dfrac{[\eta,\theta]}{|[\eta,\theta]|}, \, [\eta, \theta] \neq 0, \, x\in \R^3,
	\end{split}
	\end{align}
	where $R_w$ is the weighted Radon transform defined in \eqref{eq:wradontransform.def}, $[\cdot,\cdot]$ denotes the standard vector product in $\R^3$.
	
	Using formulas from \eqref{eq:reduction.main-formula} one can reduce initial data $P_Wf(\gamma), \, \gamma \, || \, \Sigma_\eta$, to new one given by $R_wf(s,\theta), (s,\theta)\in \R\times \Sp^2$, which gives the desired reduction of Problem~\ref{pr:problem-pw} to Problem~\ref{pr:problem-rw}. One can note that $w$ and $R_wf$ are not defined for direction $\theta$ which is parallel to $\eta$ (see formulas of \eqref{eq:reduction.main-formula}).
	For directions $\theta = \pm \eta$, weight $w$ could be set identically to zero, which automatically sets $R_wf(\cdot, \pm \eta)=0$. From a numerical point of view such setting is not crucial because it is performed on a set of  measure zero, so in practice one can always choose sampling grids which do not contain directions $\theta = \pm \eta$.
	
	\subsection{Discretization of preprocessing formulas}
	\label{subsect:preprocessing.discretizations}
	The preprocessing procedure consists of implementing formulas of \eqref{eq:reduction.main-formula}, therefore, a set of coordinate systems and discretizations must be introduced. 
	
	We assume that $\eta = (0, 0, 1)^T\in \R^3$, hence, plane $\Sigma_\eta$ (see Figure~\ref{fig:slice-by-slice-reconstructions}) corresponds to $XY$-plane in $xyz$-coordinate system in $\R^3$. Data $P_Wf$ are sampled uniformly in $z$-variable (i.e., for rays which belong to slices $z=const$) and also uniformly in $(s, \varphi)$-variables ($s\in [-1,1], \, \varphi \in [0, 2\pi)$) which stand for ray coordinates in slices $z=const$. To sum up, spectral data $P_Wf$ was given on grid $\Gamma$ defined as follows:
	\begin{align}
	\nonumber
	&\Gamma = \{\gamma = \gamma(z_i, s_j, \varphi_k) : z_i = -R + (i-1)\Delta z, \, 
	s_j = -R + (j-1) \Delta s, \, \varphi_k = (k-1) \Delta \varphi \}, \\
	\label{eq:numerical.simulation-noiseless.gamma1}
	&i=1,\dots, \, n_z, \, j = 1,\dots, \, n_s, \, k = 1,\dots, \, n_{\varphi}, \\
	\nonumber
	&\Delta z = 2R / (n_z-1), \, \Delta s = 2R / (n_s-1), \, \Delta \varphi = 2\pi / n_{\varphi},
	\end{align}
	where $R$ is the radius of image support,  ray parametrization $\gamma(z, s, \varphi)$ is given by the formula:
	\begin{align}
	\label{eq:numerical.simulation-noiseless.raygamma}
	\begin{split}
	\gamma(\mathfrak{z},s,\varphi) = \{(x,y,z) : x = s\cos(\varphi) - t\sin(\varphi), \, 
	y = s\sin(\varphi) + t\cos(\varphi), \, z = \mathfrak{z}, \,t\in \R \}.
	\end{split}
	\end{align}
	
	\noindent The sampling grid on the set of oriented 2-planes in $\R^3$ was defined as follows:
	\begin{align}\label{eq:numerical.grid.pi}
	&\Pi = \{(s_i, \theta(\varphi_j, \psi_k)) : s_i - \text{ grid on }[-R,R],\, \theta(\varphi_j, \psi_k) - \text{grid on }\Sp^2\},
	\end{align} 
	where
	\begin{align}
	\begin{split}\label{eq:numerical.grid.pi.props}
	&s_i = -R + (i-1) \Delta s, \, i = 1, \dots, \, n_s, \, \Delta s = 2R / (n_s-1),\\
	&\varphi_j = (j-1) \Delta \varphi, \, j = 1, \dots, \, n_\varphi, \, \Delta \varphi = 2\pi / n_\varphi, \\
	&\psi_k = \arccos(t_k), \, k = 1, \dots, n_\psi,\\
	&\{t_k\}_{k=1}^{n_\psi} \text{-- points for Gauss-Legendre quadrature on $[-R,R]$}
	\end{split}
	\end{align} 
	and $\theta(\varphi, \psi) = (\sin(\psi)\cos(\varphi), \sin(\psi)\sin(\varphi), \cos(\psi)), \, \varphi \in [0,2\pi], \, \psi \in [0,\pi]$ is the parametri-zation on the sphere~$\Sp^2$. 
	
	Note that in formulas \eqref{eq:numerical.simulation-noiseless.gamma1}, \eqref{eq:numerical.grid.pi.props} we use the same notations for variables $s_i, \varphi_j, n_\varphi, n_s, R$. This is because for our numerical tests we set these values the same for both grids $\Gamma$ and~$\Pi$. In principle, one can choose parameters for $\Gamma$ and $\Pi$ independently, so the only restriction to satisfy is the Shannon-Nyquist type condition for weighted Radon transforms; see, e.~g., \cite{natterer1986book} (Chapter 3), \cite{kunyansky2001formula}.
	
	For grids $\Gamma, \Pi$ formulas of \eqref{eq:reduction.main-formula} can be rewritten as follows:
	
	\begin{align}
	\label{eq:numerical.reduction.formulas-discr}
	\begin{split}
	&R_wf(s_i,\theta(\varphi_j, \psi_k)) = 
	\int\limits_{-\sqrt{R^2-s^2}}^{\sqrt{R^2-s^2}}P_{W}f(\gamma(z(s_i,\varphi_j,\psi_k,\tau), \, \sigma(s_i,\varphi_j,\psi_k,\tau), \, \varphi_j)) \, d\tau,\\
	&\gamma = \gamma(z, \sigma, \varphi) \text{ is defined in \eqref{eq:numerical.simulation-noiseless.raygamma}}, \, z\in [-R,R], \, \sigma\in [-R,R], \, \varphi\in [0,2\pi),\\
	&z(s,\varphi, \psi,\tau) = s\cos(\psi) + \tau \sin(\psi),\\
	&\sigma(s,\varphi,\psi,\tau) = s\sin(\psi) - \tau \cos(\psi),
	\end{split}\\
	\label{eq:numerical.weight.formulas-discr}
	&w(x,\theta(\varphi,\psi)) = W\left(x,\theta\left(\varphi + \frac{\pi}{2}, \frac{\pi}{2}\right)\right).
	\end{align}
	Parameter $\tau$ in \eqref{eq:numerical.reduction.formulas-discr} plays the role of  parametrization of rays $\gamma(z,\sigma,\varphi)$ which are parallel to $XY$ plane and ``fiber'' oriented plane $(s,\theta(\varphi, \psi))$. The limits of integration in formula \eqref{eq:numerical.reduction.formulas-discr} (i. e., $\pm\sqrt{R^2-s^2}$) correspond to the assumption that $f$ is supported in a ball of radius~$R$.
	
	To evaluate integrals $R_w(s_i, \theta(\phi_j, \psi_k))$ from  \eqref{eq:numerical.reduction.formulas-discr} we apply trapezoidal rule in $\tau$ with step $\Delta \tau = \Delta s = 2R / (n_s-1)$. Note that for given $\tau$ term $P_Wf(\gamma(z(\dots, \tau), \sigma(\dots, \tau), \varphi_j)$ from  \eqref{eq:numerical.reduction.formulas-discr}, generally, is not contained in  grid $\Gamma$, so some interpolations of $P_Wf(\gamma), \, \gamma \in \Gamma$ must be used. For example, in our experiments on modeled and real SPECT data we used, first, spline interpolation in variable~$\sigma$ for three different nodes of $z$ (which surround $z(\dots, \tau)$ from formula \eqref{eq:numerical.reduction.formulas-discr}) and then we used them for  quadratic interpolation in $z$-variable. Finally, interpolation in variable $\varphi$  was not needed since we used same nodes $\{\varphi_j\}_{j=1}^{n_\varphi}$ in both grids $\Gamma$, $\Pi$.
	
	We would like to note here that the aforementioned interpolations, strictly speaking, are not fully correct, because they do not respect the structure of image of $P_{W}$. For our tests on SPECT data, the impact of using the values slightly of outside of $\mathrm{Im}(P_W)$ is not so crucial compared to the impact of strong Poisson noise in emission data. 
	
	In addition, such high orders of interpolations (quadratic, splines) were used because of artifacts in reconstructions when using piecewise linear interpolations. This could be explained by the property that $P_{W}$, in general, is a smoothing operator (see, for example, \cite{natterer1986book} for the case of $W\equiv 1$). Therefore, the image of operator $P_{W}$ contains smoother functions and higher orders of interpolations than linear must be used.\\

	\section{Experiment on modeled SPECT data}
	\label{sect:numericalexp-modelized}
	
	The goal of the test on modeled SPECT data is to see if reconstructions from $R_wf$ given in  \eqref{eq:numerical.reduction.formulas-discr}, \eqref{eq:numerical.weight.formulas-discr}, become more stable in presence of noise than slice-by-slice reconstructions from $P_Wf$. For this, we begin by recalling the model for data in SPECT.

	\subsection{Model for emission data in SPECT}
	\label{subsect:numericalexp-modelized.wrt-spect}
	
	In SPECT the measured data are described by photon counts $N(\gamma)$, where $\gamma$ are the rays from the grid $\Gamma$ defined in \eqref{eq:numerical.simulation-noiseless.gamma1}. On a mathematical level,  photon counts $N(\gamma)$ are related to transforms $P_W$ by the following formula:
	\begin{equation}\label{eq:spect.statistics}
	N(\gamma) \sim \mathrm{Po}(Ct P_{W_a}f(\gamma)), \, \gamma \in \Gamma,
	\end{equation}
	where 
	\begin{align}
	\begin{split}\label{eq:spect.statistics.params}
	&\mathrm{Po}(\cdot) \text{ denotes the Poisson distribution}, \\
	&C \text{ is a positive constant which depends on parameters of the scanner  setup},\\
	&t \text{ is the acquisition time per ray},\\
	&P_{W_a} \text{ is the transform from formula~\eqref{eq:wraytransform.def} for weight }W_a,\\
	&f=f(x), \, x\in \R^3 \text{ is the density distribution of the isotope}
	\end{split}
	\end{align}
	and weight $W_a$ is given by the formulas:
	\begin{align}
	\label{eq:spect.w.def}
	&W_a(x,\theta) = \exp(-Da(x,\theta)),\\
	\label{eq:spect.d.def}
	&Da(x,\theta) = \int\limits_{0}^{+\infty}a(x+t\theta) \, dt, \, 
	x\in \R^3, \, \theta\in \Sp^2, \\
	\label{eq:spect.a.def}
	&a=a(x), \, x\in\R^3 \text{ is the attenuation map}.
	\end{align}
	
	\begin{remark}
		In case when one has access to $P_{W_a}f$ directly, isotope distribution $f$ can be reconstructed via efficient backprojection-type methods based on Novikov's analytic inversion formula \cite{novikov2002formula}, \cite{kunyansky2001formula}, \cite{natterer2001inversion}. However, this  corresponds to the case of weak noise in data which in practice is never fulfilled. Low signal-to-noise ratio for data described by $N(\gamma), \, \gamma \in \Gamma$, is one of the main problems in applying  backprojection-type algorithms in SPECT \cite{guillement2008wiener}, \cite{kunyansky2001formula}, \cite{guillement2002spect}. In particular, the preprocessing procedure proposed here can be seen as a method to increase signal-to-noise ratio in the measured data.
	\end{remark}
	
	We use formulas \eqref{eq:spect.statistics}-\eqref{eq:spect.a.def} to modelize ray data $P_{W_a}f(\gamma), \, \gamma \in \Gamma$, for different isotope distributions, attenuation maps and noise-levels. Modelization of isotope distributions and attenuation maps is straightforward via the choice of different functions $f, \, a$, noise levels in data can be controlled using parameter $t$ in \eqref{eq:spect.statistics}, \eqref{eq:spect.statistics.params}. The details of our modelizations are given in Subsections~\ref{subsect.numerical.modelisation_spect.attphantom}-\ref{subsect.numerical.modelisation_spect.noiseless}.

	\subsection{Organization of the experiment on modeled data}
	\label{subsect:numericalexp-modelized.organization}
	
	Numerical experiment on modeled SPECT data consists of the following steps:
	
	\begin{enumerate}
		\item Modelization of SPECT data $W_a, \, P_{W_a}f(\gamma), \gamma\in \Gamma$ via formulas from Subsection~\ref{subsect:numericalexp-modelized.wrt-spect} for two different distributions $f_1, f_2$, two attenuation maps $a_1, a_2$ (weak and strong attenuation levels, respectively) and for three different noise-levels -- no noise, weak and strong levels, respectively.
		
		\item Reduction of $P_{W_a}f(\gamma), \, \gamma \in \Gamma$ to $R_wf(s,\theta), \, (s,\theta) \in \Pi$ via formulas from Subsection~\ref{subsect:preprocessing.discretizations}. 
		\item Reconstructions of $f$ from $P_{W_a}f(\gamma), \, \gamma \in \Gamma$ and from $R_wf(s,\theta), \, (s,\theta)\in \Pi$, using 2D and 3D versions of Chang-type methods \cite{chang1978correciton}, \cite{novikov2011chang}, \cite{goncharov2016analog} and Kunyansky-type methods \cite{kunyansky1992generalized}, \cite{guillement2014finite},
		\cite{goncharov2017iterative}.
		
		\item Computation of reconstruction errors due to noise impact for 2D and 3D methods, respectively.
	\end{enumerate}
	The details of our modelizations are presented in  Subsections~\ref{subsect.numerical.modelisation_spect.attphantom}-\ref{subsect.numerical.modelisation_spect.noise}. The reconstruction results for Chang-type and Kunyansky-type methods with analysis of errors are given in Subsections~\ref{subsect:reconstruct.chang}, \ref{subsect.numerical.inversion.iterative}. 
	
	For our experiments we used parameters $n_z = n_s = 129, \, R = 1.0, \, n_\psi = 128$ (see Subsection~\ref{subsect:preprocessing.discretizations}). The number of directions $n_\varphi$ was chosen 
	equals to $128$. Note, that this number is less than in \cite{kunyansky2001formula} ($n_\varphi = 400$), where 
	the value for $n_{\varphi}$ was chosen according to Shannon-Nyquist type sampling condition for ray transforms. Our choice for $n_\varphi = 128$ was motivated by the experiment on real data for which we had only this number of directions in each slicing plane.

	\subsection{Attenuation phantoms}
	\label{subsect.numerical.modelisation_spect.attphantom}
	
\begin{figure}[H]
		\centering
		\begin{subfigure}{.3\textwidth}
			\centering
			\includegraphics[scale=0.25]{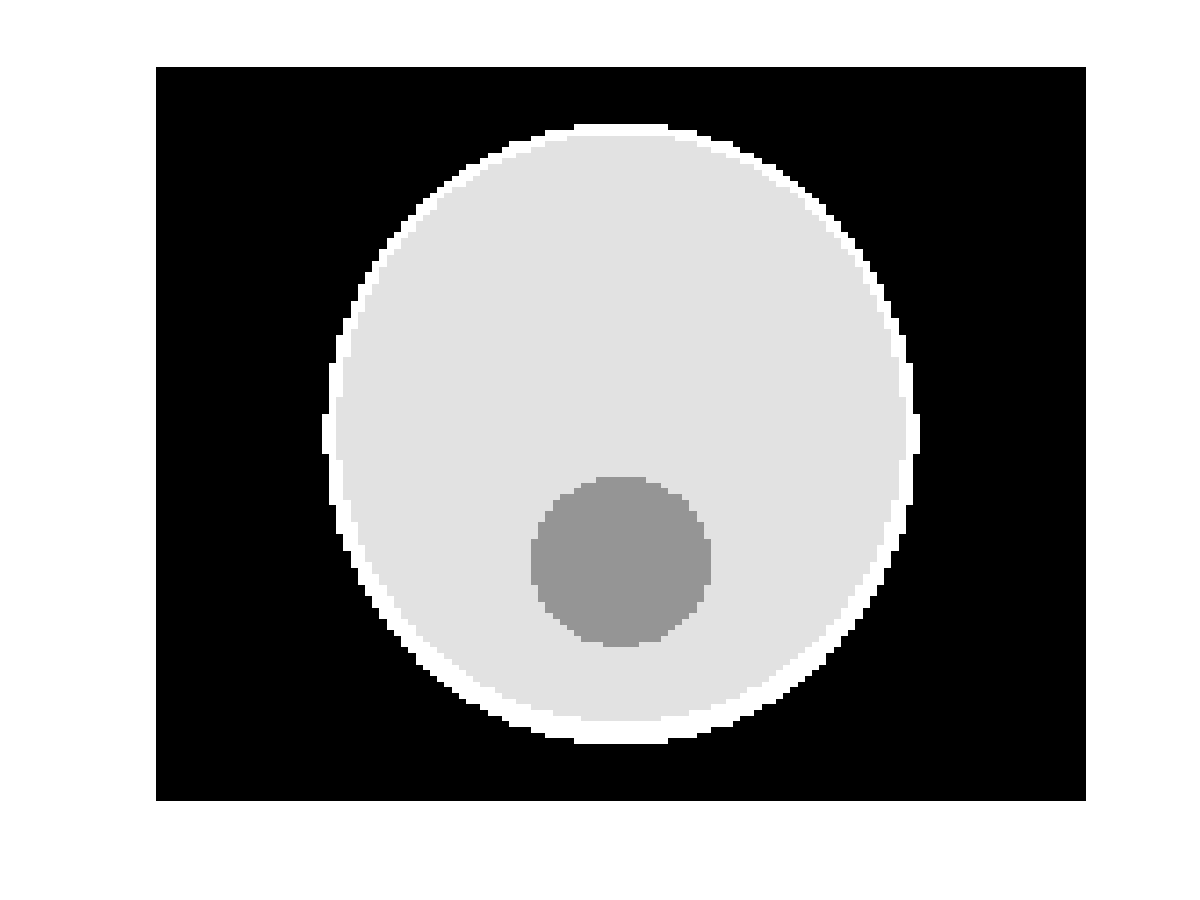}
			\caption{ $i_z = 40$ }
		\end{subfigure}
		\begin{subfigure}{.3\textwidth}
			\centering
			\includegraphics[scale=0.25]{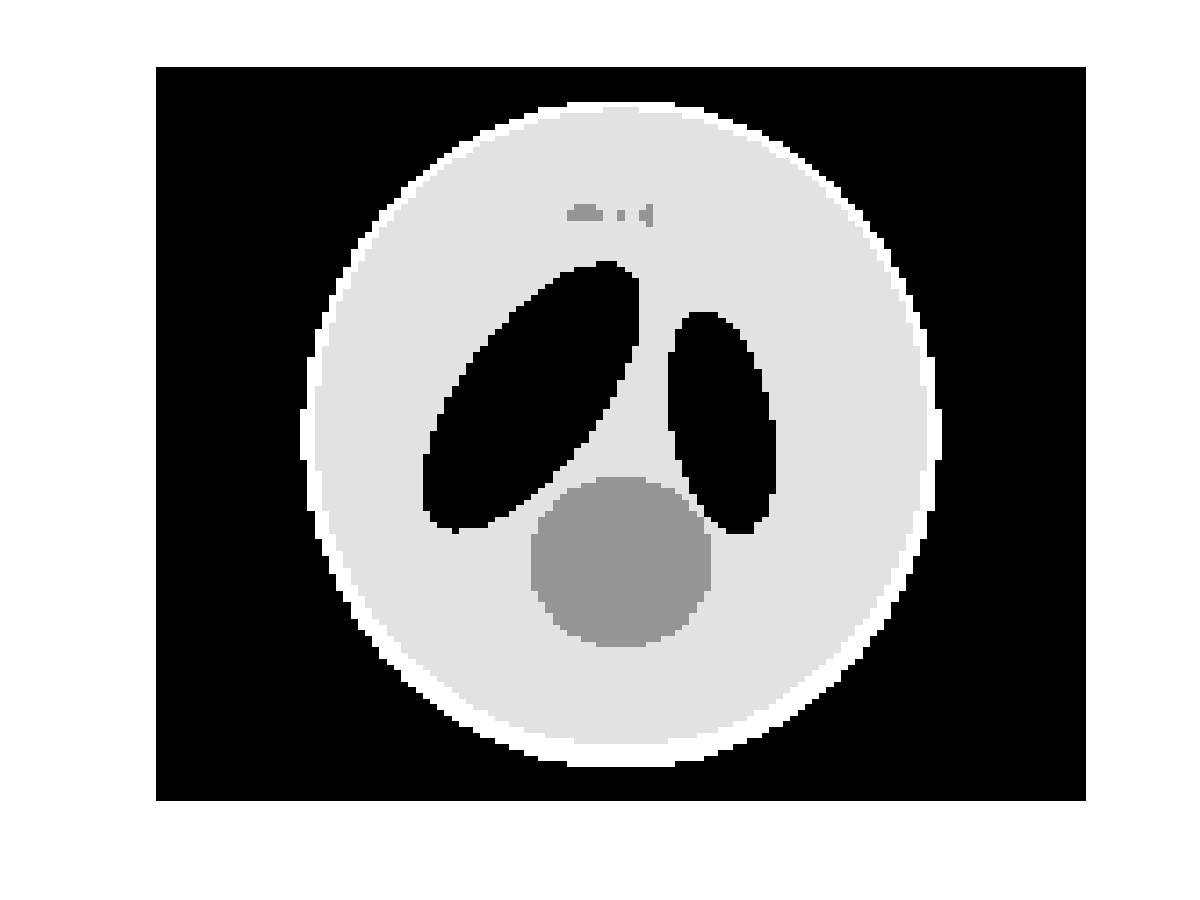}
			\caption{$i_z = 64$}
		\end{subfigure}
		\begin{subfigure}{.3\textwidth}
			\centering
			\includegraphics[scale=0.25]{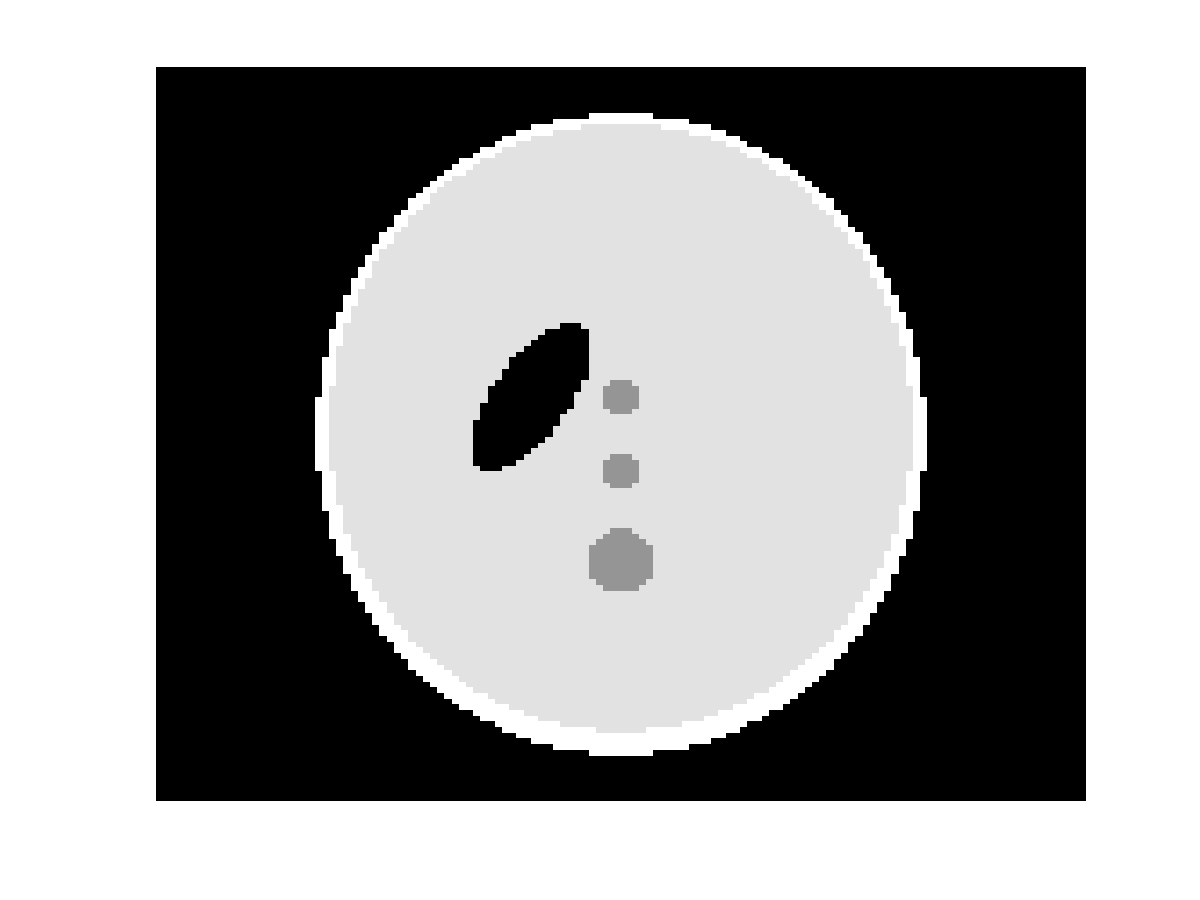}
			\caption{$i_z = 80$}
		\end{subfigure}
		\caption{cross-sections of the Shepp-Logan phantom\\ by planes $z=const$}
		\label{fig:numerical.shepp_logan}
	\end{figure}

	For the attenuation model we used the extension of the famous of Shepp-Logan phantom to 3D (\cite{gach2008shepplogan}, \cite{shepp1974fourier}); see Figure~\ref{fig:numerical.shepp_logan}. 
	This phantom imitates the human head and is standard for testing reconstruction algorithms for tomographies. 
	
	The model consists of a large ellipsoid which represents the brain and of several smaller ellipsoids which represent its features. The outer ellipsoidal layer with axes of length 13.8cm, 18cm, and 18.4cm corresponds to bone of the cranium; the standard attenuation of the bone material is equal to $0.17\text{cm}^{-1}$. The attenuation within smaller inner ellipsoid representing the brain material and is equal to $0.15\text{cm}^{-1}$. Smaller spheres (dark gray spheres on Figure~\ref{fig:numerical.shepp_logan}) represent inclusions in the brain and have attenuation of $0.10\text{cm}^{-1}$. The model also includes cavities which are given by two ellipsoidal regions near the center of the phantom (see the ellipsoids in black color, Figure~\ref{fig:numerical.shepp_logan} (B), (C)). In these regions the attenuation was set to zero.
	\par Quantitative feature which describes the general strength of the attenuation is the optical length of attenuation along central axes $X$, $Y$, $Z$ (the optical length along a path is an integral of the attenuation coefficient along 
	this path). In the given model the optical lengths along $X$, $Y$, $Z$ axes are equal to 2.44, 3.89 and 4.81, respectively. This corresponds to the case of  strong attenuation. 
	\par To investigate the efficiency of our reconstructions methods for different attenuation regimes we used two models. The first one, which is referred by $a_1$ or \textit{``strong attenuation''}, is the Shepp-Logan phantom with parameters described above. The second model, which is referred by $a_2$ or \textit{``weak attenuation''}, corresponds to the same phantom as the first one but the attenuation values are multiplied by a factor of $1/10$ (i.e., $a_2 =a_1\cdot 10^{-1}$). Hence, the optical lengths along $X,Y,Z$ axis for the second attenuation model $a_2$ are equal to 0.244, 0.389 and 0.481, respectively. From the practical point of view, this corresponds to the case of very weak attenuation.
	\par Before to proceed, we note that all our images in this subsection will be presented in a linear grayscale, where darker colors correspond to smaller values (black color corresponds to zero).

	\subsection{Phantoms for isotope distributions}
	\label{subsect.numerical.nucphantom}
\begin{figure}[H]
		\centering
		\begin{subfigure}{.4\textwidth}
			\centering
			\includegraphics[height=50mm, width=55mm]{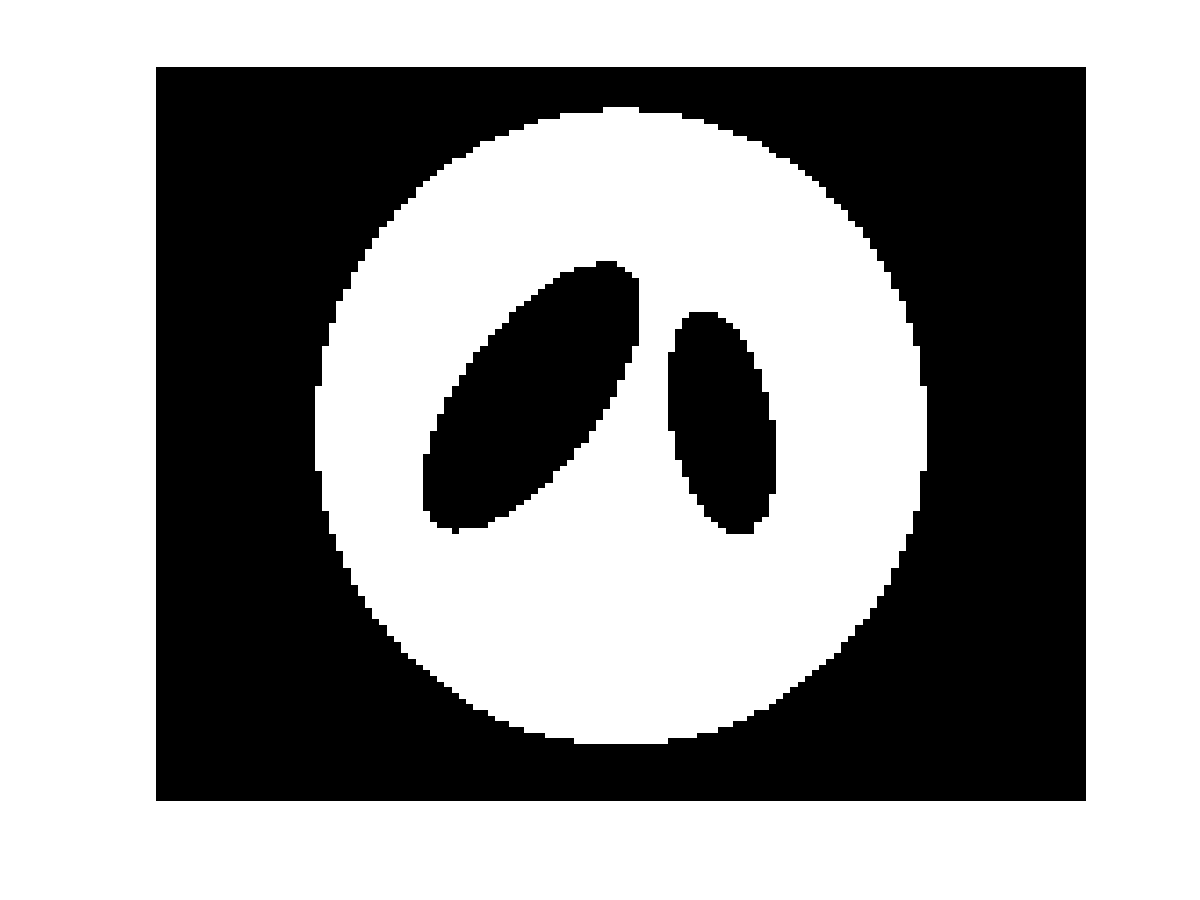}
			\caption{ Phantom 1, $f_1, \, i_z = 64$ }
		\end{subfigure}
		\begin{subfigure}{.4\textwidth}
			\centering
			\includegraphics[height=50mm, width=55mm]{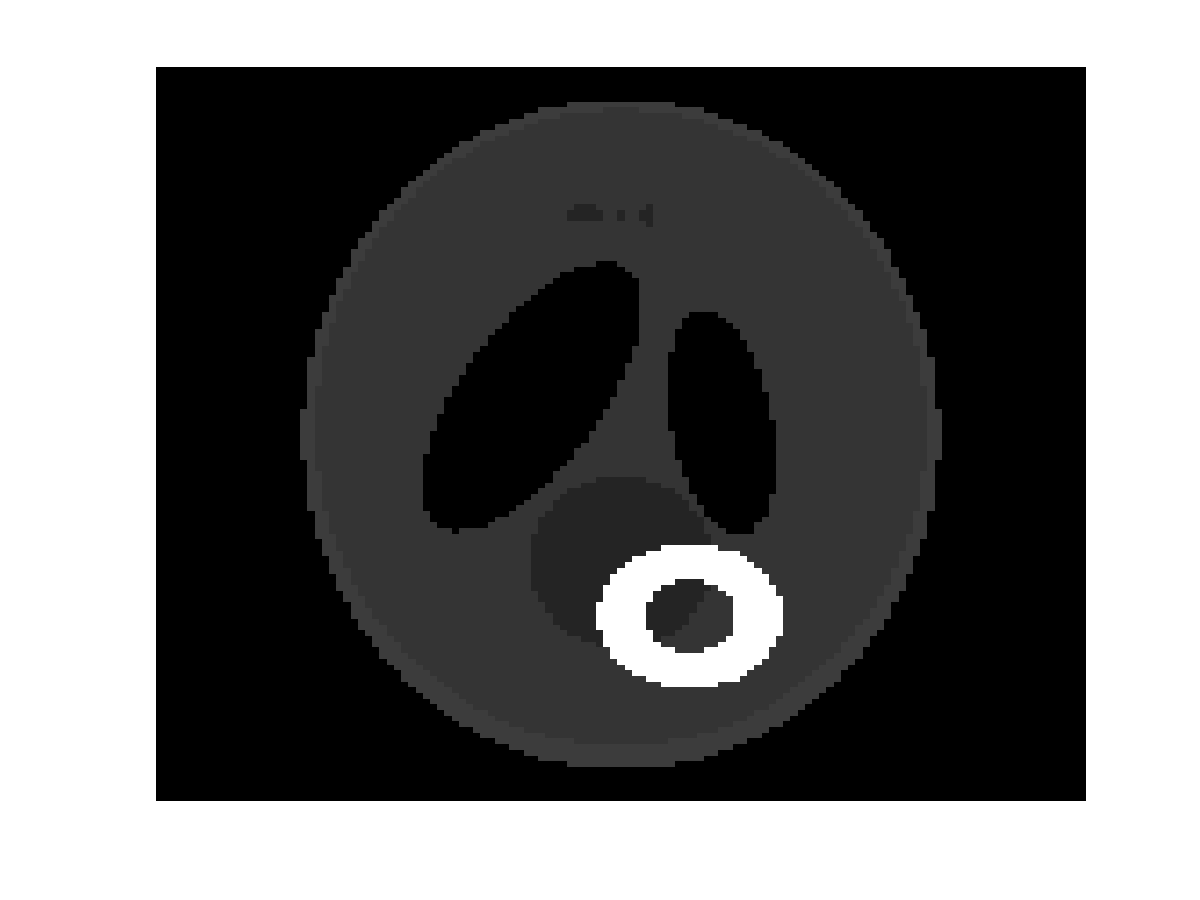}
			\caption{Phantom 2, $f_2, \, i_z = 64$}
		\end{subfigure}
		\caption{isotope distributions for modeled data}
		\label{fig:numerical.att}
	\end{figure}

	We used two models for isotope distributions. The first one, further denoted by $f_1$ or \textit{Phantom 1}, is given by the characteristic function of the inner ellipsoid of the Shepp-Logan phantom from the previous paragraph; see Figure~\ref{fig:numerical.att} (A). 
	
	\noindent The purpose of Phantom 1 is to test the ability to reconstruct spatially uniform and slowly varying distributions of the nucleotide. 
	\par At the same time it is important to test the algorithms when activity distribution has highly non-uniform distribution in space. For these reasons we use phantom $f_2$ or\textit{ Phantom 2}, which is described by a spherical layer placed almost in the center of the attenuation phantom; see Figure~\ref{fig:numerical.att} (B) (the background gray model is displayed only for visualization of position of $f_2$ with respect to the attenuation model). The outer radius of the spherical layer is $4\text{cm}$ 
	and the radius of the inner one is equal to $2\text{cm}$.
	In particular, Phantom 2 is used 
	to simulate the classical SPECT setting: reconstruct the distribution of an  isotope which tends to concentrate in the brain of a patient.
	\par Finally, note that on Figure~\ref{fig:numerical.att} only slices of Phantoms 1, 2 are presented, whereas the latter are complete three-dimensional objects.

	\subsection{Simulation of noiseless data $P_{W_a}f$ }
	\label{subsect.numerical.modelisation_spect.noiseless}
	
	\par Activity phantoms $f_1, \, f_2$ and attenuation maps $a_1, \, a_2$ are  given by their values on a Cartesian grid of the unit cube $[-1,1]^3$ and it is assumed that these functions are supported in the centered ball of radius $R = 1.0$.
	
	In particular, Cartesian grid $\Omega_N$ on $[-1,1]^3$ was defined as follows:
	\begin{align}\nonumber
	&\Omega_N = \{(x_i, y_j, z_k) : x_i = -R + (i-1)\Delta x, \, y_j = -R + (j-1)\Delta y, \, z_k = -R + (k-1)\Delta z\}, \\
	& \Delta x = \Delta y = \Delta z = 2R / (N-1), \, i,\, j, \, k\in \{1, \dots N\},
	\end{align}
	where $N$ is the number points in the grid in one direction. In all our simulations we used $N=129$. For evaluations of $W_a, P_{W_a}f$ from Subsection~\ref{subsect:numericalexp-modelized.wrt-spect} we assumed that functions $f_1, \, f_2$ and $a_1, \, a_2$ are continuous piecewise linear  between the grid points. For computation of integrals along rays we used the trapezoidal integration rule.
	
\begin{figure}[H]
		\begin{subfigure}[b]{0.5\textwidth}
			\centering
			\includegraphics[scale=.35]{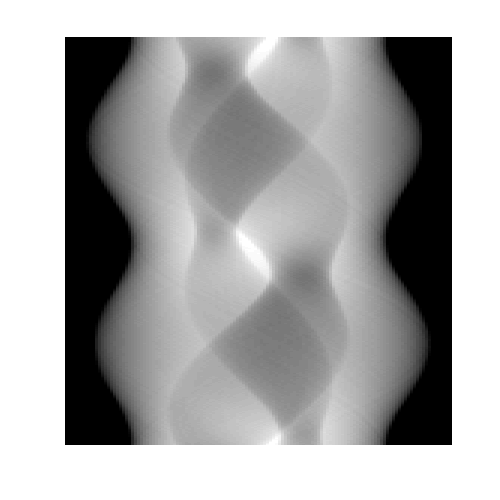}
			\caption{ $P_{W_a}f_1$, weak attenuation }
		\end{subfigure}%
		\begin{subfigure}[b]{0.5\textwidth}
			\centering
			\includegraphics[scale=0.35]{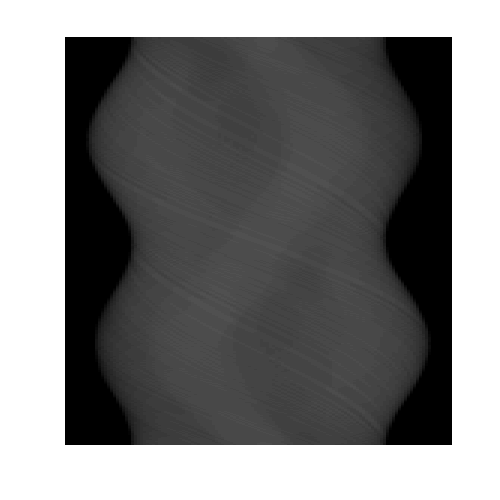}
			\caption{ $P_{W_a}f_1$, weak attenuation }
		\end{subfigure}
		\centering
		\begin{subfigure}[b]{0.5\textwidth}
			\centering
			\includegraphics[scale=.35]{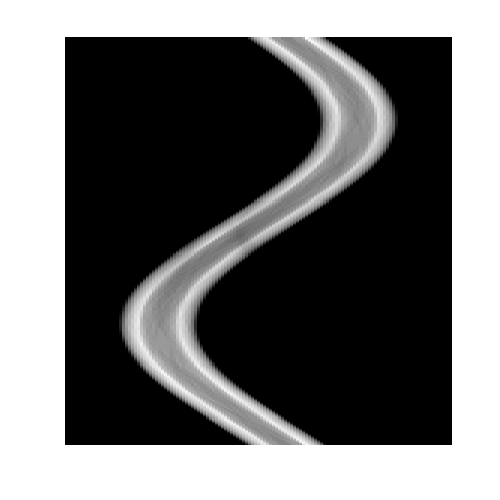}
			\caption{ $P_{W_a}f_2$, weak attenuation}
		\end{subfigure}%
		\begin{subfigure}[b]{0.5\textwidth}
			\centering
			\includegraphics[scale=0.35]{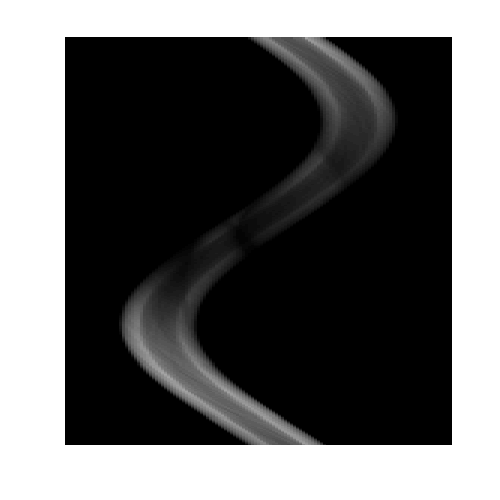}
			\caption{ $P_{W_a}f_2$, strong attenuation}
		\end{subfigure}
		\caption{ data $P_{W_a}f$ for section $z=0$; horizontal and vertical axes correspond to variables $(s,\varphi)\in [-1,1]\times [0, 2\pi)$, 
			respectively}
		\label{fig:numerical.phantomsatt}
	\end{figure}

	In Figure~\ref{fig:numerical.phantomsatt} we give examples for values of weighted ray transforms $P_{W_a}f(z, s, \varphi)$ for section $z = 0$ (such images are also called \textit{sinograms}).
	\par In the next paragraph, we explain how we modelized different noise levels for ray data $P_{W_a}f(\gamma), \, \gamma\in \Gamma$.
	
	\subsection{Modelling of the noise in emission data}
	\label{subsect.numerical.modelisation_spect.noise}
	
	In view of considerations from  Subsection~\ref{subsect:numericalexp-modelized.wrt-spect}, we modelize noise in emission data as follows.
	\par First, for given attenuation $f\in \{f_1, f_2\}$ and activity phantom $a \in \{a_1, a_2\}$ we compute noiseless data $P_{W_a}f(\gamma), \, \gamma \in \Gamma$. Then we choose a normalization constant $C=C_{n}$
	such that 
	\begin{equation}\label{eq:numerical.modelisation-noise.lmax}
	C_{n} \max_{\gamma\in \Gamma} P_{W_a}f(\gamma) = n,
	\end{equation}
	where $n$ corresponds to the maximal number of photons registered along rays in $\Gamma$ if the time of exposures was $t=1$ for all rays. In our experiments we used two different values for $n$: $n_{strong} = 50$ and $n_{weak} = 500$ which correspond to strong and weak noise levels, respectively (see also Remarks~\ref{rem:numerical.modelisation-noise.approximation}, \ref{rem:numerical.modelisation-noise.noise-control} below).
	
	\par Next, according to formula \eqref{eq:spect.statistics}, we define intensities $\lambda(\gamma), \, \gamma \in \Gamma$, for Poisson process by the formula 
	\begin{equation}\label{eq:numerical.modelisation-noise.intensitydef}
	\lambda_n(\gamma) = C_n P_{W_a}f(\gamma), 
	\end{equation}
	where $C_n$ is a constant of \eqref{eq:numerical.modelisation-noise.lmax}. Finally, we generate independent Poisson random variables $N(\gamma)$ with intensities $\lambda_n(\gamma)$ for each $\gamma \in \Gamma$.
	
	\begin{remark}\label{rem:numerical.modelisation-noise.approximation}
		In realistic SPECT procedures setup-dependent constant $C$ and exposure time $t$ from \eqref{eq:spect.statistics.params}, 
		are known. Therefore, for the experiment on modeled data it is fine to assume that $C_n$ is also known
		and it can be used to approximate unknown emission data. More precisely, using $C_n$ and $N(\gamma), \, \gamma\in \Gamma$, one can approximate $P_{W_a}f(\gamma)$ as follows:
		\begin{equation}\label{eq:numerical.modelisation-noise.approximation}
		P_{W_a}f(\gamma) \approx \dfrac{N(\gamma)}{C_n}, \, \gamma \in \Gamma.
		\end{equation}
		In particular, approximation \eqref{eq:numerical.modelisation-noise.approximation} is based on the following property:
		\begin{equation}\label{eq:numerical.modelisation-noise.validexpectation}
		N(\gamma) \sim \mathrm{Po}\, (C_n P_{W_a}f(\gamma)), 
		\text{ then } \mathbb{E}\left(\dfrac{N(\gamma)}{C_n}\right) = P_{W_a}f(\gamma), \, \gamma\in \Gamma,
		\end{equation}
		where $\mathrm{Po}(\cdot)$ denotes the Poisson distribution, $\mathbb{E}$ - denotes the mathematical expectation. To obtain formula \eqref{eq:numerical.modelisation-noise.validexpectation} we used that for $\xi\sim \mathrm{Po}(\lambda)$ its expectation $\mathbb{E}\xi$ is $\lambda$.
		We recall that, after having modelized noise for the emission data, one  cannot longer assume that values for $P_{W_a}f$ are known.
	\end{remark}
	
	\begin{remark}\label{rem:numerical.modelisation-noise.noise-control}
		The choice of constant $n$ in \eqref{eq:numerical.modelisation-noise.lmax} allows to control the noise level in generated emission data. More precisely, by choosing intensities $\lambda(\gamma)$ as in \eqref{eq:numerical.modelisation-noise.intensitydef} we have that
		\begin{equation}\label{eq:numerical.modelisation-noise.variance-control}
		\mathrm{Var}\left(
		\dfrac{N(\gamma)}{C_n}
		\right) = \dfrac{P_{W_a}f(\gamma)}{C_n} = \dfrac{P_{W_a}f(\gamma) \max\limits_{\gamma\in \Gamma}P_{W_a}f(\gamma)}{n}, \, \gamma \in \Gamma,
		\end{equation}
		where $\mathrm{Var}(\cdot)$ denotes the variance. Identity \eqref{eq:numerical.modelisation-noise.variance-control} is based
		on the following properties of Poisson distribution: $\xi \sim \mathrm{Po}(\lambda), \, \mathrm{Var}(\xi) = \lambda$. From \eqref{eq:numerical.modelisation-noise.intensitydef}-\eqref{eq:numerical.modelisation-noise.variance-control} is follows that
		that higher values for $n$ correspond to better approximations in \eqref{eq:numerical.modelisation-noise.approximation}. In particular, our choice for $n_{strong} = 50$ corresponds to the case of strong noise and it is actually close to the setting of our experiment on real data where the maximal number of photons per measurement was 60; see also  \cite{guillement2008wiener}. Finally, the choice for large $n$ in gives better approximations in \eqref{eq:numerical.modelisation-noise.approximation} because this corresponds to \textit{gaussian approximation of Poisson's law} (which is valid, for example, in X-ray transmission tomography).
	\end{remark}

\begin{figure}[H]
		\vspace{-2.0cm}
		\begin{subfigure}[b]{0.5\textwidth}
			\centering
			\includegraphics[scale=.3]{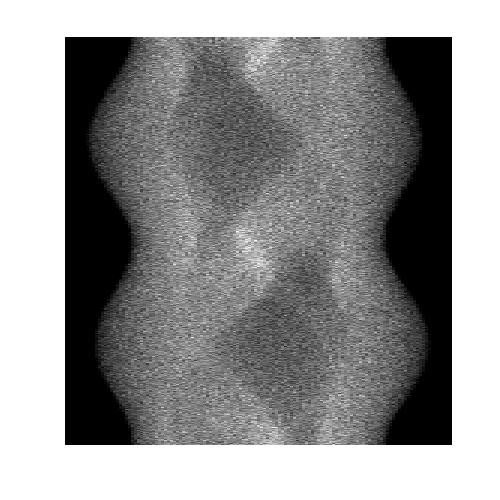}
			\caption{ $P_{W_a}f_1$, weak attenuation }
		\end{subfigure}%
		\begin{subfigure}[b]{0.5\textwidth}
			\centering
			\includegraphics[scale=0.3]{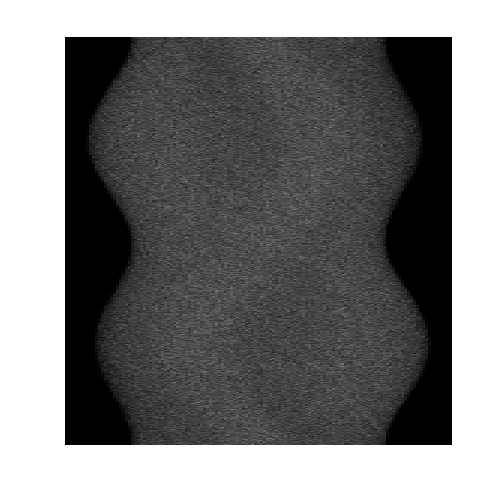}
			\caption{ $P_{W_a}f_1$, strong attenuation}
		\end{subfigure}
		\centering
		\begin{subfigure}[b]{0.5\textwidth}
			\centering
			\includegraphics[scale=.3]{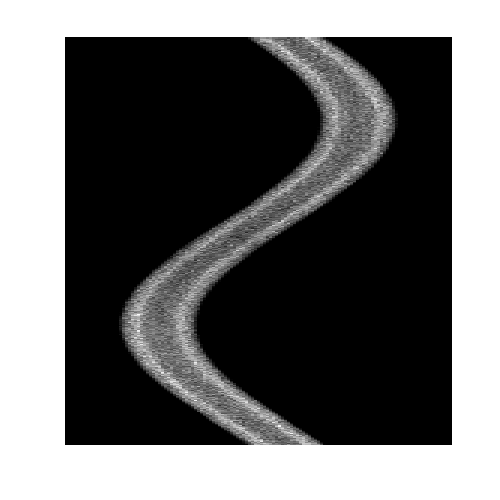}
			\caption{ $P_{W_a}f_2$, weak attenuation}
		\end{subfigure}%
		\begin{subfigure}[b]{0.5\textwidth}
			\centering
			\includegraphics[scale=0.3]{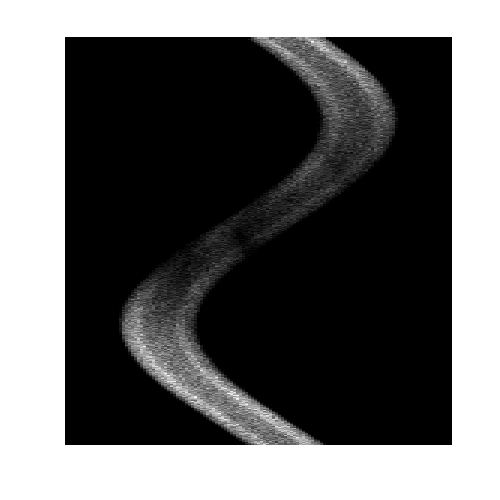}
			\caption{ $P_{W_a}f_2$, strong attenuation}
		\end{subfigure}
		\caption{ Modelling of strong noise $n = n_{strong} = 50$}
		\label{fig:numerical.phantomsatt.noise}
	\end{figure}

\begin{figure}[H]
		\begin{subfigure}[b]{0.5\textwidth}
			\centering
			\includegraphics[scale=.30]{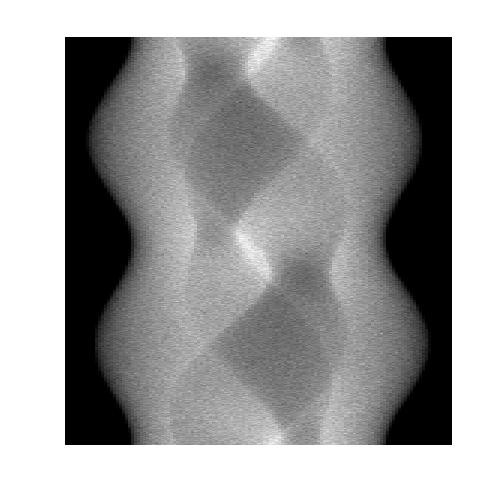}
			\caption{ $P_{W_a}f_1$, weak attenuation }
		\end{subfigure}%
		\begin{subfigure}[b]{0.5\textwidth}
			\centering
			\includegraphics[scale=0.30]{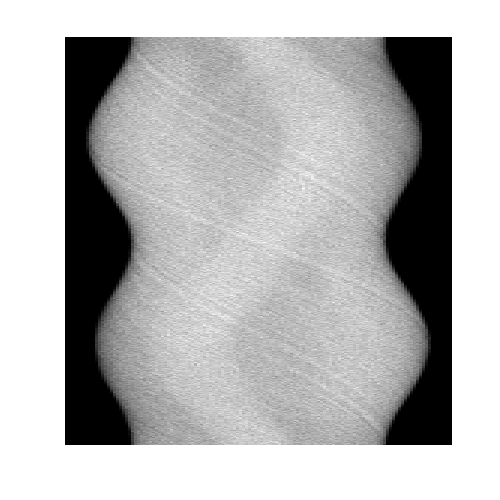}
			\caption{ $P_{W_a}f_1$, strong attenuation}
		\end{subfigure}
		\centering
		\begin{subfigure}[b]{0.5\textwidth}
			\centering
			\includegraphics[scale=.30]{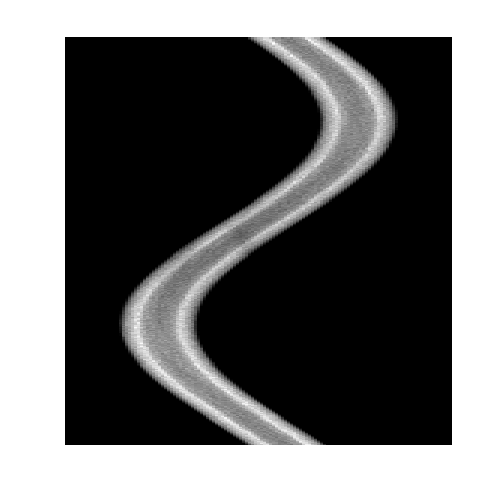}
			\caption{ $P_{W_a}f_2$, weak attenuation}
		\end{subfigure}%
		\begin{subfigure}[b]{0.5\textwidth}
			\centering
			\includegraphics[scale=0.30]{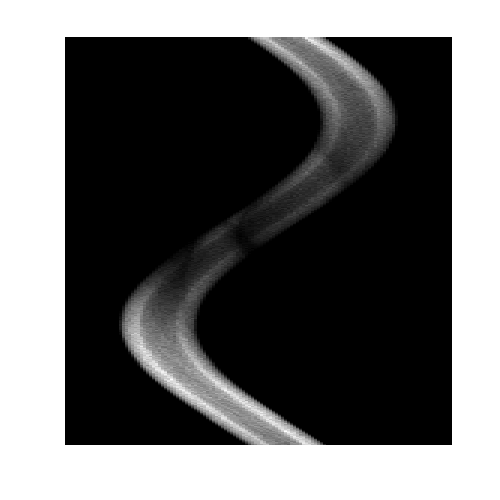}
			\caption{ $P_{W_a}f_2$, strong attenuation}
		\end{subfigure}
		\caption{ Modelisation of weak noise $n = n_{weak} = 500$}
		\label{fig:numerical.phantomsatt.noise.weak}
	\end{figure}

	In Figures~\ref{fig:numerical.phantomsatt.noise},~\ref{fig:numerical.phantomsatt.noise.weak} the modeling of noisy data $P_{W_a}f$ is given for  
	rays in section $z = 0$.
	
	Note also that for the case of weak noise the images almost do not differ from the denoised originals in Figure~\ref{fig:numerical.phantomsatt}. This was explained in Remark~\ref{rem:numerical.modelisation-noise.noise-control} that noise level decreases fast with increase of number of registered photons.

	\subsection{Reconstruction methods}\label{subsect:reconstruct.methods}
	The details for Chang-type and Kunyansky-type methods we used can be found in \cite{goncharov2016analog} and in \cite{kunyansky1992generalized}, \cite{guillement2014finite}, \cite{goncharov2017iterative}, respectively. Here we briefly discuss these methods in order to show that both admit similar variants in 2D and 3D for inversions of weighted Radon transforms.
	
	\subsubsection*{Chang-type methods}
	The original method of L. Chang \cite{chang1978correciton} was proposed in order to correct non-uniform attenuation for reconstructions in SPECT. At the same time, the method admits direct extensions for the case of weighted Radon transforms $R_W$ for arbitrary weights and for all dimensions.
	
	Given weight $W$ and $R_Wf$ in $\R^d$, where $R_W$ is defined analogously to formula \eqref{eq:wradontransform.def}, the method returns an approximation $f_{appr}\approx f$ defined by the formulas:
	\begin{align}\label{eq:reconstruction-methods.chang.fappr}
	&f_{appr}(x) = \dfrac{R^{-1}R_Wf(x)}{w_0(x)}, \, 
	x\in \R^d,\\ \label{eq:reconstruction-methods.w0-def}
	&w_0(x) = \dfrac{1}{\mathrm{Vol}(\Sp^{d-1})}\int\limits_{\Sp^{d-1}}W(x,\theta)\, d\theta, \, w_0(x) \neq 0,
	\end{align}
	where $R^{-1}$ is inverse of the classical Radon transform in $\R^d$. 
	
	For $d=2$ the method is usually used in the framework of slice-by-slice reconstruction approach. For example, in our experiments we apply it on data given by $P_{W_a}f$ for rays which belong to plane slices $z=const$ (see  formula \eqref{eq:numerical.simulation-noiseless.gamma1}). For $d=3$ 
	we apply it on preprocessed data given by $R_wf(s,\theta), \, (s,\theta)\in \Pi$ (see also formulas \eqref{eq:numerical.reduction.formulas-discr}, \eqref{eq:numerical.weight.formulas-discr}).

	\subsubsection*{Kunyansky-type methods}
	Kunyansky-type methods can be seen as an extension of Chang-type methods when approximation $W\approx w_0$ (see formula \eqref{eq:reconstruction-methods.w0-def}) is not sufficient. In this case the reconstructions from $W, \, R_Wf$ are given by functions $f_m, \, m\in \mathbb{N}\cup \{0\}$, which are defined as solutions of the following integral equation:
	\begin{align}\label{eq:reconstruction-methods.kunyansky.integr-eq}
	(I + Q_{W, \, D, \,m})(w_0f_m) = R^{-1}R_Wf,
	\end{align}
	where
	\begin{align}\nonumber
	&w_0 \text{ is defined in \eqref{eq:reconstruction-methods.w0-def}}, \, I \text{ is the identity operator in }L^2(\R^d), \\
	\begin{split}\label{eq:reconstruction-methods.qwd.def}
	&Q_{W, \, D, \, m} \text{ is a linear bounded integral operator in }L^2(\R^d)
	\text{ which depends }\\
	&\text{on weight } W, \text{ domain $D$ surrounding the support of $f$ (i.e., $\supp f\subset D$)}\\
	&\text{and on parameter } m\in \mathbb{N}\cup \{0\}.
	\end{split}
	\end{align} 
	Formulas with explicit expressions for $Q_{W, \, D, \, m}$ and theoretical analysis of Kunyansky's type methods can be found in \cite{kunyansky1992generalized}, \cite{guillement2014finite} for dimension $d=2$ and in \cite{goncharov2017iterative} for $d\geq 3$. The main idea of the method is to solve equation \eqref{eq:reconstruction-methods.kunyansky.integr-eq} by the method of successive approximations, which converge, for example, when 
	\begin{equation}\label{eq:reconstruction-methods.qwd.estimate}
	\|Q_{W,\,D, \,m}\|_{L^2(\R^d)\rightarrow L^2(\R^d)} < 1.
	\end{equation}
	Parameter $m\in \mathbb{N}\cup \{0\}$ from \eqref{eq:reconstruction-methods.qwd.def} controls the error between $f$ and approximation $f_m$ and should be chosen as large as possible provided condition \eqref{eq:reconstruction-methods.qwd.estimate} being satisfied. Larger values for $m$ correspond to better approximations $f_m \approx f$ and, in particular, for $m=+\infty$ function $f_m$ coincides with $f$. At the same time, the operator norm of $Q_{W, \, D, \, m}$ grows with $m$, so for $m$ too large condition \eqref{eq:reconstruction-methods.qwd.estimate} may not be  satisfied. To choose parameter $m$ efficiently one uses the following properties of $Q_{W, \, D, \, m}$:
	\begin{align}\label{eq:reconstruction-methods.kunyansky.qwd-zero}
	&Q_{W, \, D, \, m} = 0 \text{ for } m=0, \\
	\label{eq:reconstruction-methods.kunyansky.qwd-estimate}
	&\|Q_{W, \, D, \, m}\|_{L^2(\R^d) \rightarrow L^2(\R^d)} \leq \sigma_{W, \, D, \, m},
	\end{align}
	where $\sigma_{W, \, D, \, m}$ is an explicit and efficiently computable function of $W$, $D$, $m$ (see \cite{kunyansky1992generalized}, \cite{guillement2014finite}, \cite{goncharov2017iterative}). Therefore, the main  point is to  choose $m$ as large as possible so that 
	\begin{equation}\label{eq:reconstruction-methods.kunyansky.conv-condition}
	\sigma_{W, \, D, \, m} < 1.
	\end{equation}

	Having condition \eqref{eq:reconstruction-methods.kunyansky.conv-condition} satisfied, estimate \eqref{eq:reconstruction-methods.kunyansky.qwd-estimate} implies that 
	equation \eqref{eq:reconstruction-methods.kunyansky.integr-eq} can be solved by the method of successive approximations. Moreover, formulas \eqref{eq:reconstruction-methods.kunyansky.qwd-estimate},  \eqref{eq:reconstruction-methods.kunyansky.conv-condition} guarantee that the algorithm will converge with geometric speed, so only a few iterations of successive approximations are needed in practice. Finally, from \eqref{eq:reconstruction-methods.chang.fappr}-\eqref{eq:reconstruction-methods.kunyansky.qwd-zero} one can see that Kunyansky's method is a direct extension of Chang's method, where the latter one corresponds to the case $m=0$.
	
	Similar to Chang-type methods, for $d=2$ we apply Kunyansky-type method in the framework of slice-by-slice reconstructions on $P_{W_a}f$ reduced to rays in slices $z=const$. For $d=3$ the method is applied on preprocessed data $R_{w}f(s,\theta), \, (s,\theta)\in \Pi$ (see also formulas \eqref{eq:numerical.reduction.formulas-discr}, \eqref{eq:numerical.weight.formulas-discr}).\\
	
	Before going to numerical results we want to highlight again that in our experiments we do not use any regularization or noise-filtering schemes. We apply the aforementioned reconstruction methods directly to noisy data. Only such design of the experiment allows to see the effect of the preprocessing on reconstructions in presence of noise.
	
	\subsection{Reconstructions via Chang-type methods}
	\label{subsect:reconstruct.chang}
	
	\subsubsection*{Noiseless case}
	\label{subsubsect:reconstruct.chang.noiseless}
	\begin{figure}[H]
		\begin{subfigure}{0.24\textwidth}
			\includegraphics[width=45mm, height=40mm]{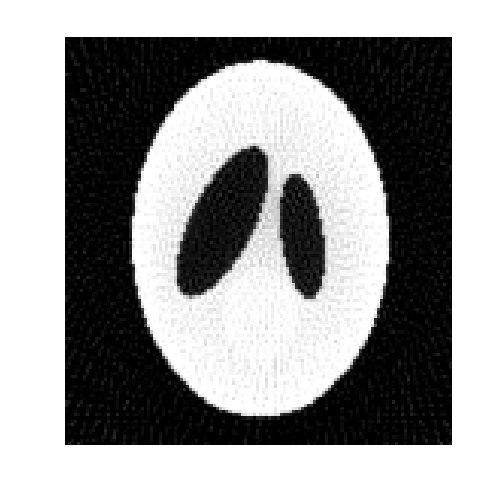}
			\caption{}
			
		\end{subfigure}
		\begin{subfigure}{0.24\textwidth}
			\includegraphics[width=44mm, height=40mm]
			{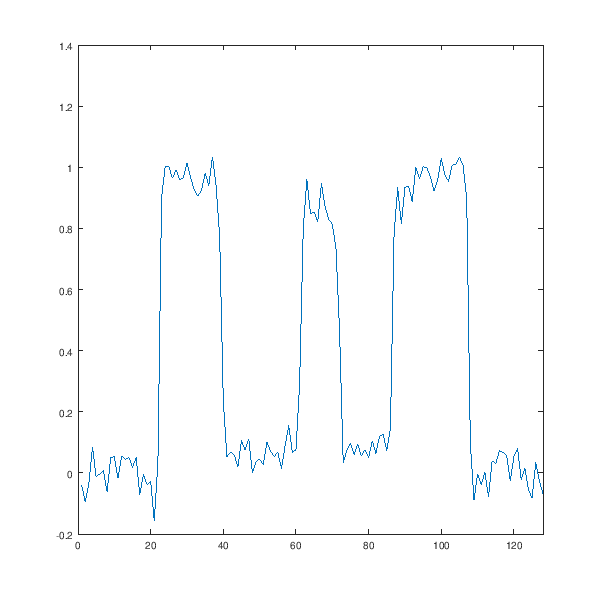}
			\caption{}
			
		\end{subfigure}
		\begin{subfigure}{0.24\textwidth}
			\includegraphics[width=45mm, height=40mm]{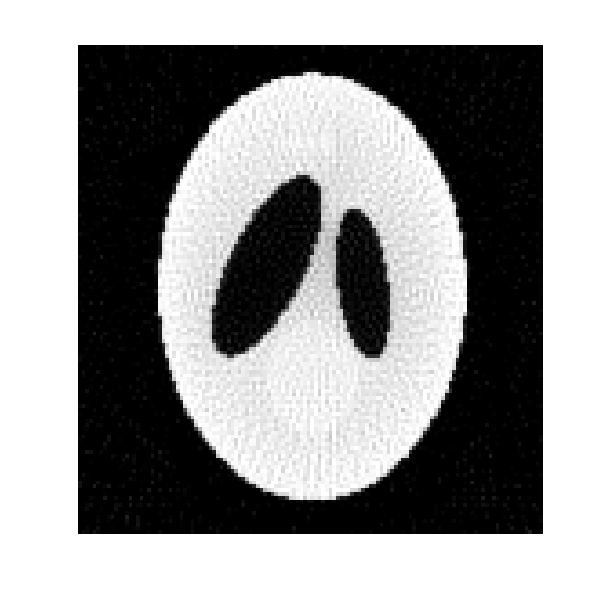}
			\caption{}
			
		\end{subfigure}
		\begin{subfigure}{0.24\textwidth}
			\includegraphics[width=45mm, height=40mm]
			{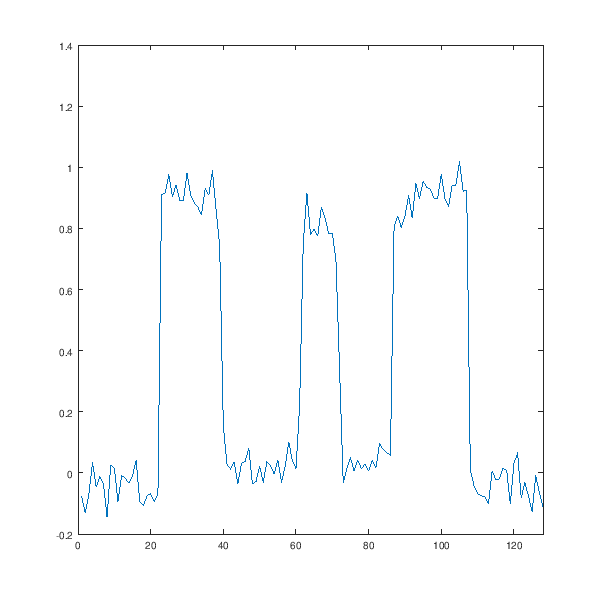}
			\caption{}
			
		\end{subfigure}
		
		\begin{subfigure}{0.24\textwidth}
			\includegraphics[width=45mm, height=40mm]
			{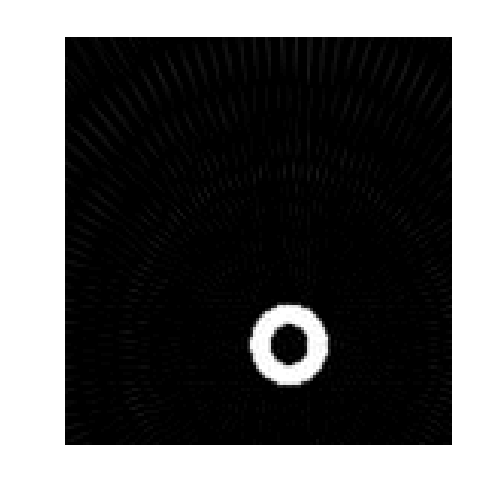}
			\caption{}
			
		\end{subfigure}
		\begin{subfigure}{0.24\textwidth}
			\includegraphics[width=45mm, height=40mm]
			{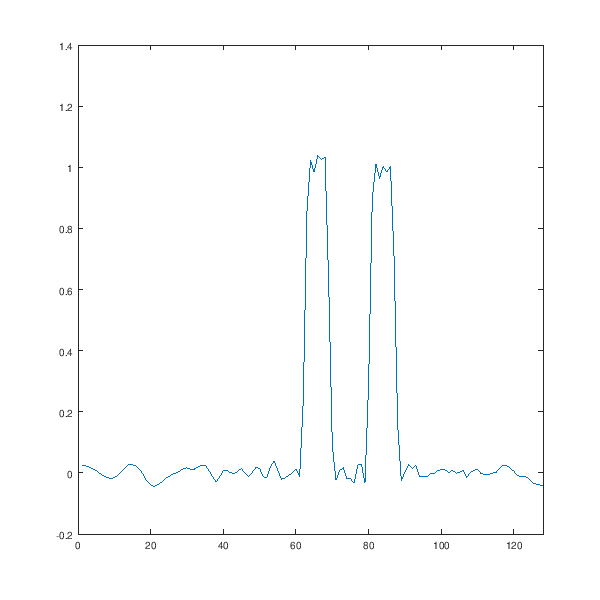}
			\caption{}
			
		\end{subfigure}
		\begin{subfigure}{0.24\textwidth}
			\includegraphics[width=45mm, height=40mm]{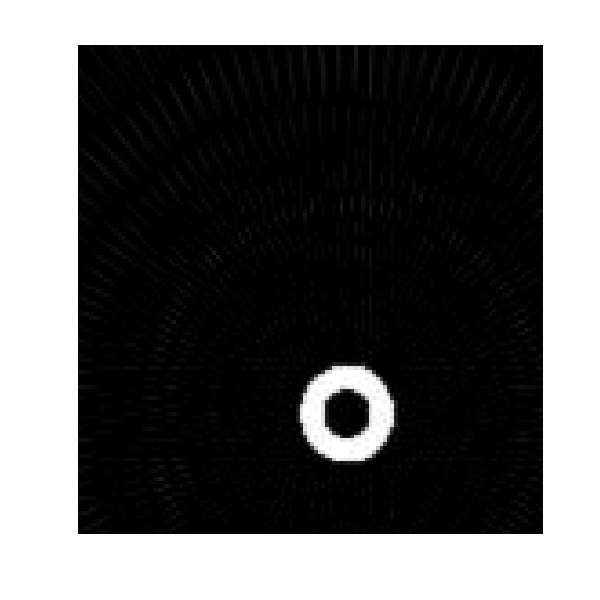}
			\caption{}
			
		\end{subfigure}
		\begin{subfigure}{0.24\textwidth}
			\includegraphics[width=45mm, height=40mm]
			{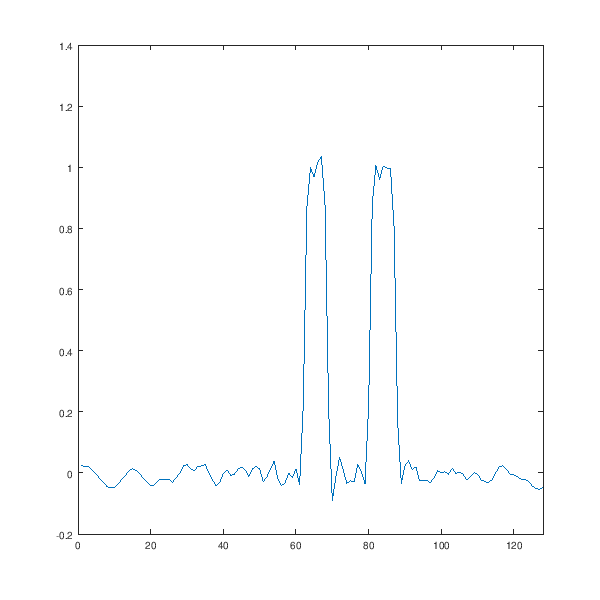}
			\caption{}
			
		\end{subfigure}
		\caption{ weak attenuation, no noise; reconstructions of $f_1, f_2$ using Chang-type methods in 3D (A), (E) and in 2D (C), (G); (B), (D), (F), (B) -- sections along $X$-axis}
		\label{fig:chang.weak.att.nonoise}
	\end{figure}
	\begin{figure}[H]
		\begin{subfigure}{0.24\textwidth}
			\includegraphics[width=45mm, height=40mm]{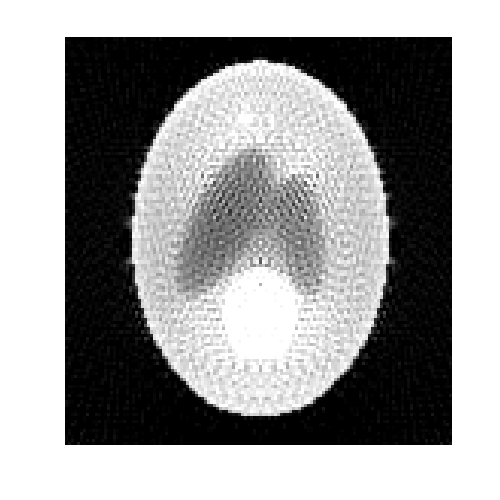}
			\caption{}
			
		\end{subfigure}
		\begin{subfigure}{0.24\textwidth}
			\includegraphics[width=45mm, height=40mm]
			{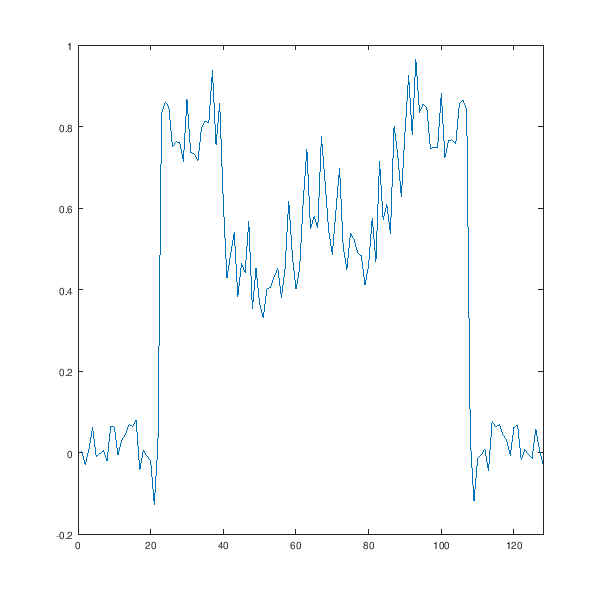}
			\caption{}
			
		\end{subfigure}
		\begin{subfigure}{0.24\textwidth}
			\includegraphics[width=45mm, height=40mm]{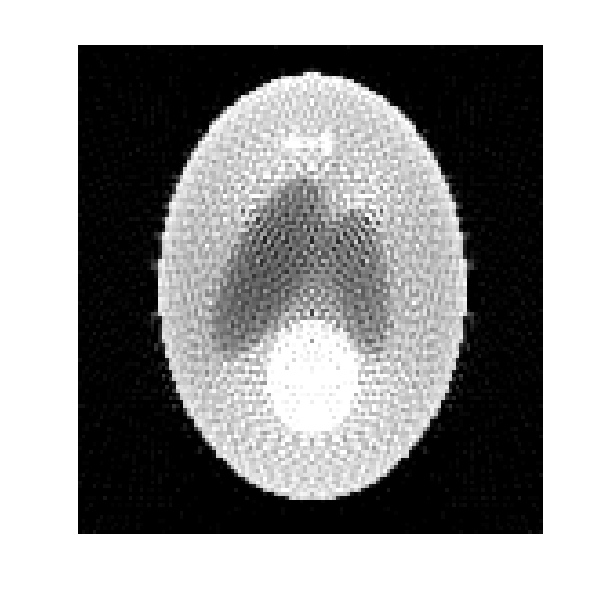}
			\caption{}
			
		\end{subfigure}
		\begin{subfigure}{0.24\textwidth}
			\includegraphics[width=45mm, height=40mm]
			{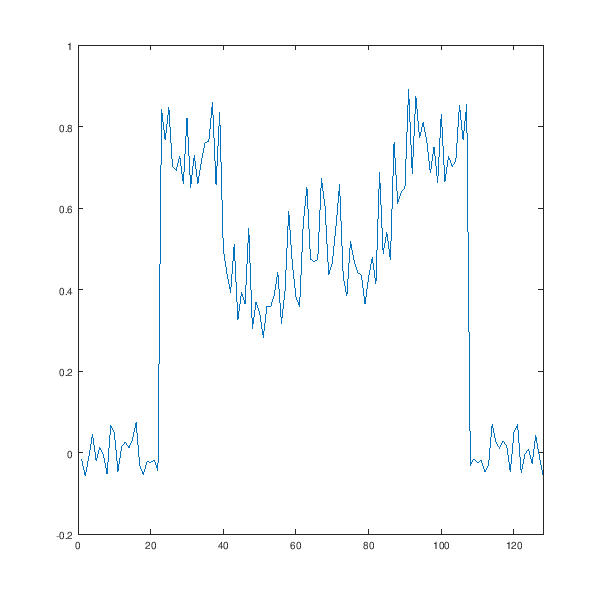}
			\caption{}
			
		\end{subfigure}
		
		\begin{subfigure}{0.24\textwidth}
			\includegraphics[width=45mm, height=40mm]
			{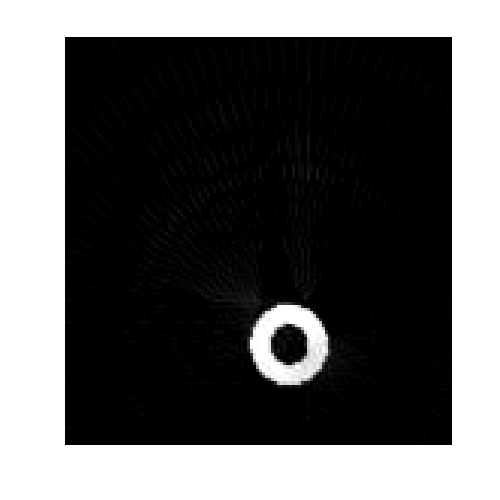}
			\caption{}
			
		\end{subfigure}
		\begin{subfigure}{0.24\textwidth}
			\includegraphics[width=45mm, height=40mm]
			{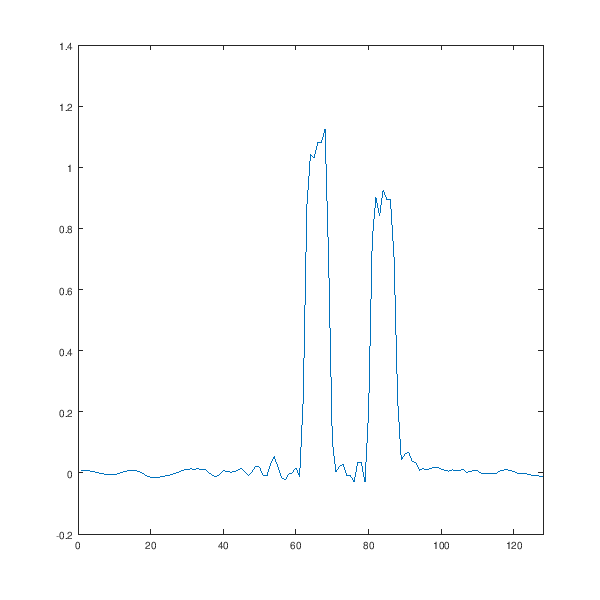}
			\caption{}
			
		\end{subfigure}
		\begin{subfigure}{0.24\textwidth}
			\includegraphics[width=45mm, height=40mm]{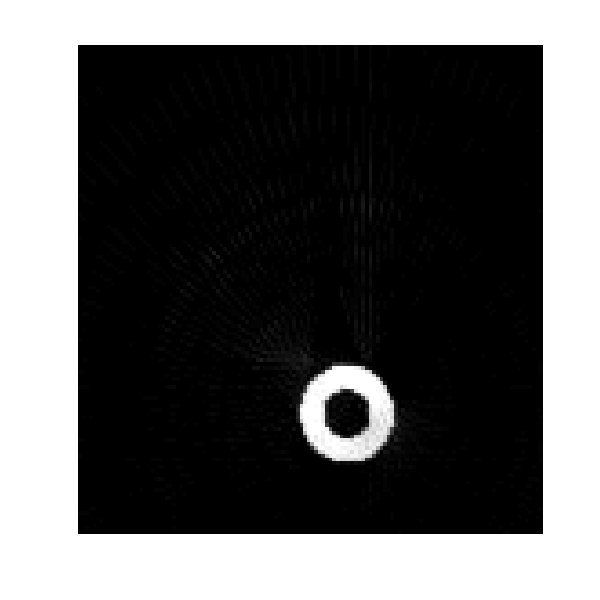}
			\caption{}
			
		\end{subfigure}
		\begin{subfigure}{0.24\textwidth}
			\includegraphics[width=45mm, height=40mm]
			{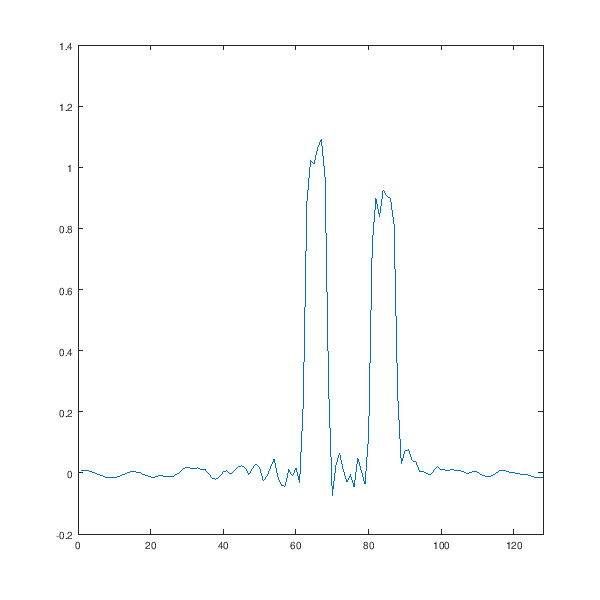}
			\caption{}
			
		\end{subfigure}
		\caption{ strong attenuation, no noise; reconstructions of $f_1, f_2$ using Chang-type methods in 3D (A), (E) and in 2D (C), (G); (B), (D), (F), (H) -- sections along $X$-axis}
		\label{fig:chang.strong.att.nonoise}
	\end{figure}

	
	\subsubsection*{Case with noise}
	\begin{figure}[H]
		\begin{subfigure}{0.24\textwidth}
			\includegraphics[width=45mm, height=45mm]{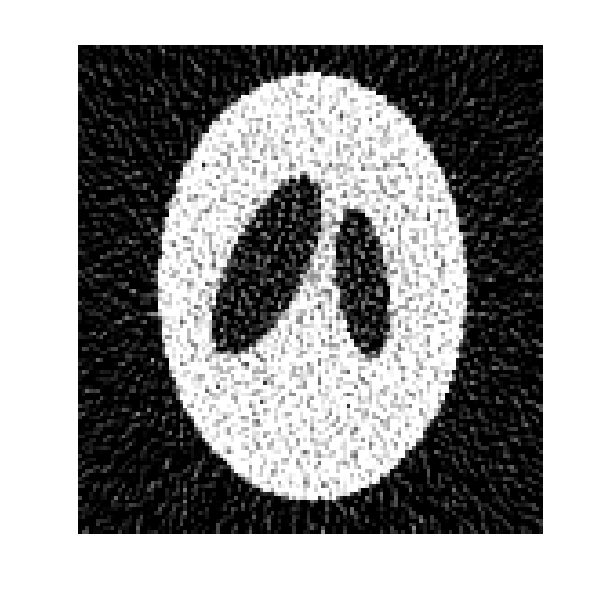}
			\caption{}
			
		\end{subfigure}
		\begin{subfigure}{0.24\textwidth}
			\includegraphics[width=45mm, height=45mm]
			{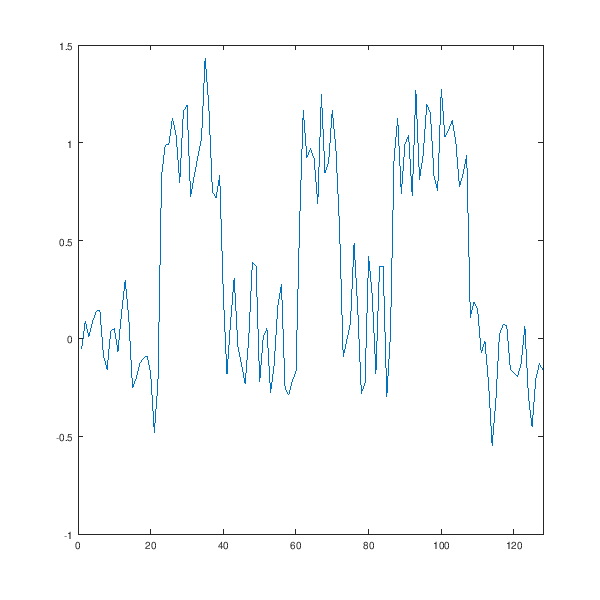}
			\caption{}
			
		\end{subfigure}
		\begin{subfigure}{0.24\textwidth}
			\includegraphics[width=45mm, height=45mm]{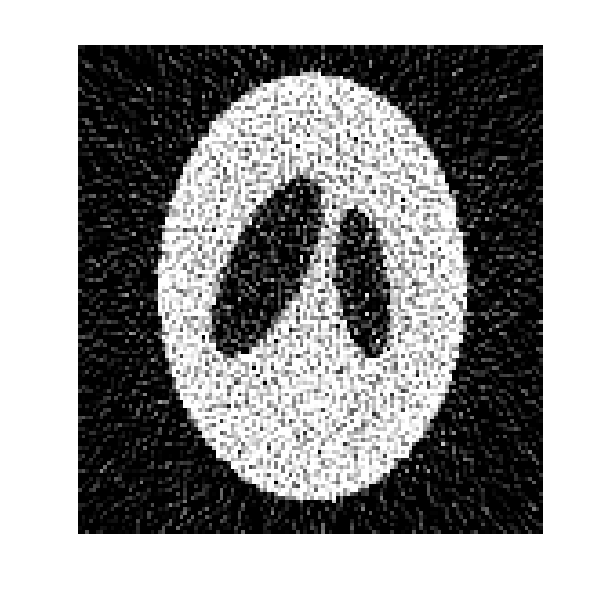}
			\caption{}
			
		\end{subfigure}
		\begin{subfigure}{0.24\textwidth}
			\includegraphics[width=45mm, height=45mm]
			{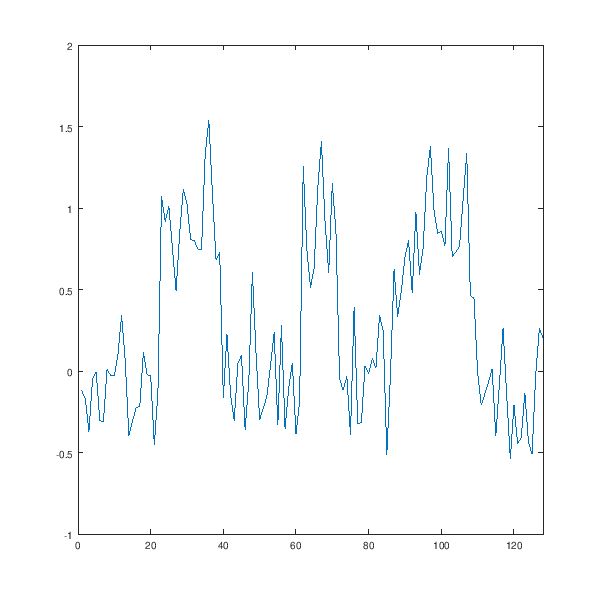}
			\caption{}
			
		\end{subfigure}
		\newline
		\begin{subfigure}{0.24\textwidth}
			\includegraphics[width=45mm, height=45mm]
			{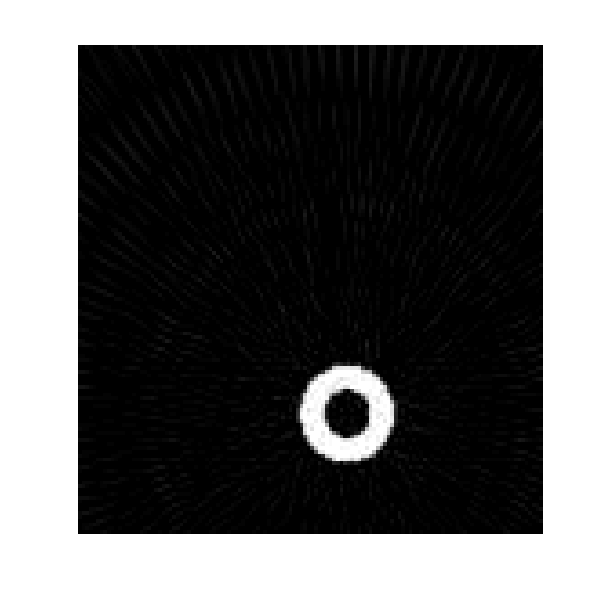}
			\caption{}
			
		\end{subfigure}
		\begin{subfigure}{0.24\textwidth}
			\includegraphics[width=45mm, height=45mm]
			{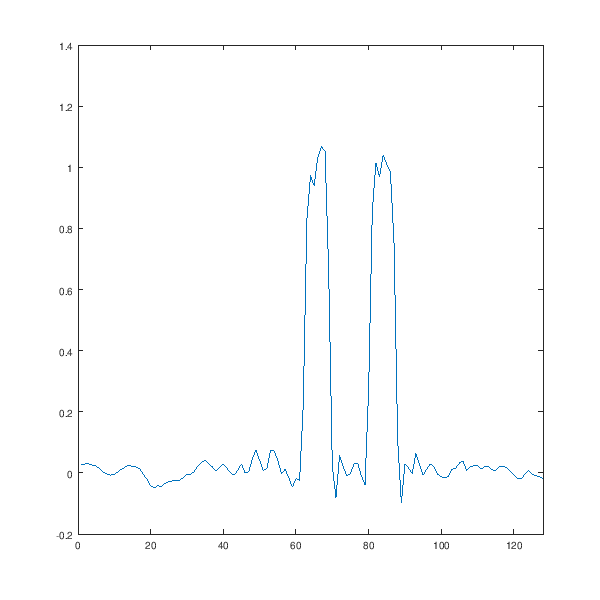}
			\caption{}
			
		\end{subfigure}
		\begin{subfigure}{0.24\textwidth}
			\includegraphics[width=45mm, height=45mm]{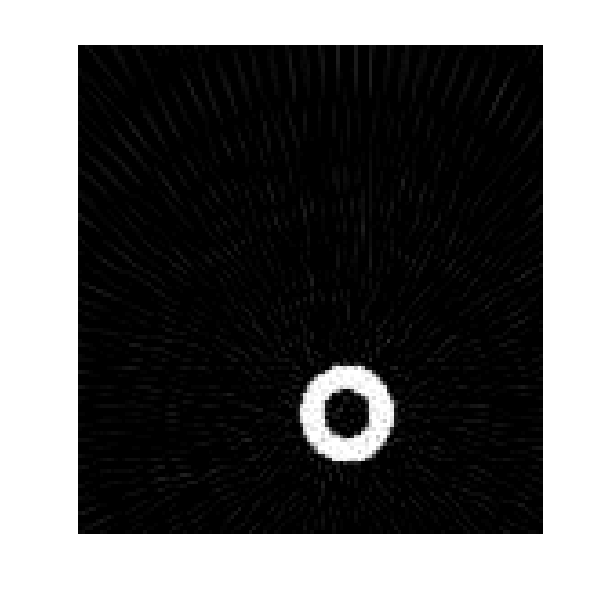}
			\caption{}
			
		\end{subfigure}
		\begin{subfigure}{0.24\textwidth}
			\includegraphics[width=45mm, height=45mm]
			{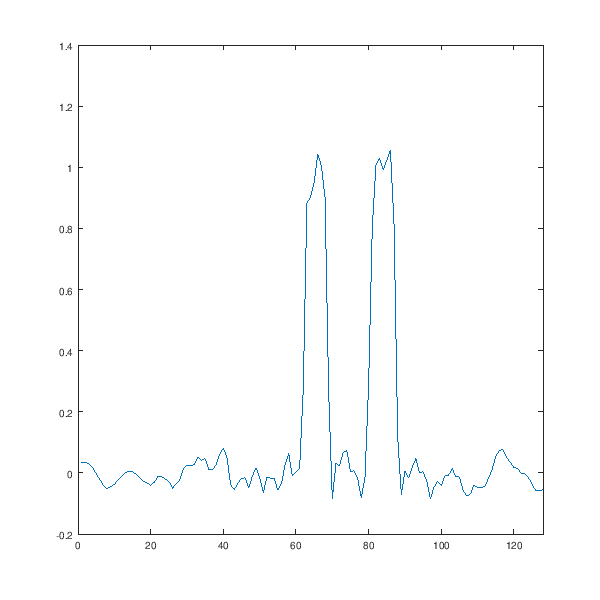}
			\caption{}
			
		\end{subfigure}
		
		\caption{ weak attenuation, weak noise ($n_{weak}=500$); reconstructions of $f_1, f_2$ using Chang-type methods in 3D (A), (E) and in 2D (C), (G); (B), (D), (F), (H) -- sections along $X$-axis}
		\label{fig:chang.weak.att.weak.noise}
	\end{figure}
	
	\begin{figure}[H]
		\begin{subfigure}{0.24\textwidth}
			\includegraphics[width=45mm, height=45mm]{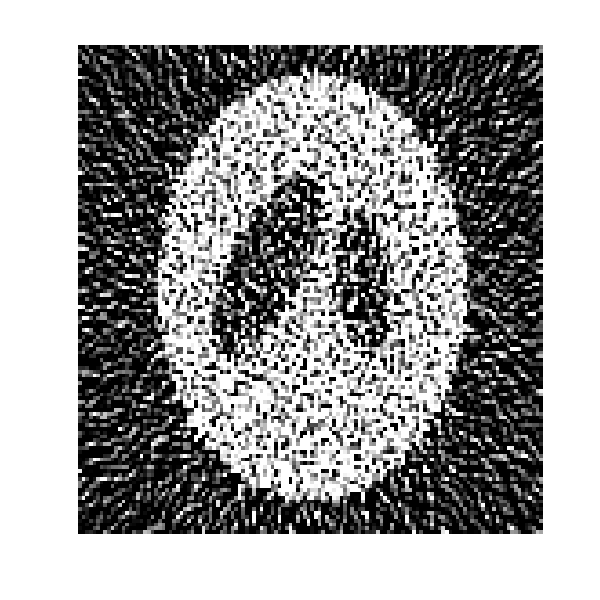}
			\caption{}
			
		\end{subfigure}
		\begin{subfigure}{0.24\textwidth}
			\includegraphics[width=45mm, height=45mm]
			{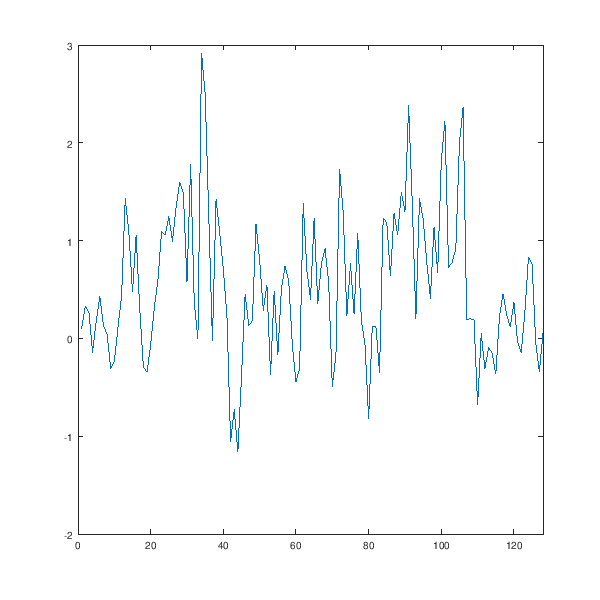}
			\caption{}
			
		\end{subfigure}
		\begin{subfigure}{0.24\textwidth}
			\includegraphics[width=45mm, height=45mm]{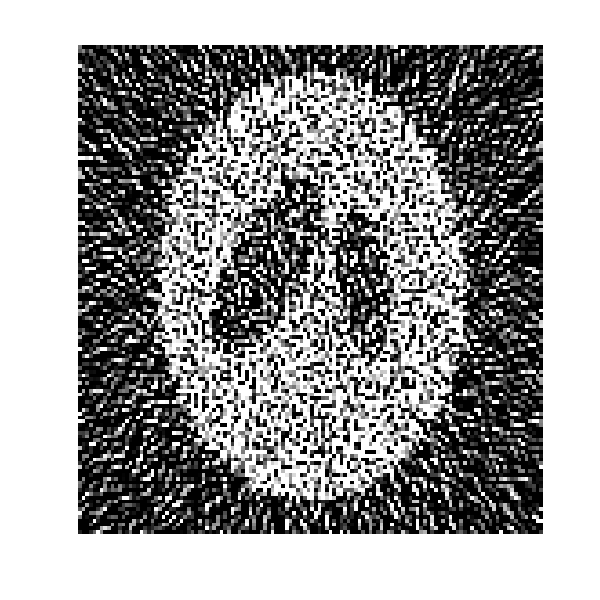}
			\caption{}
			
		\end{subfigure}
		\begin{subfigure}{0.24\textwidth}
			\includegraphics[width=45mm, height=45mm]
			{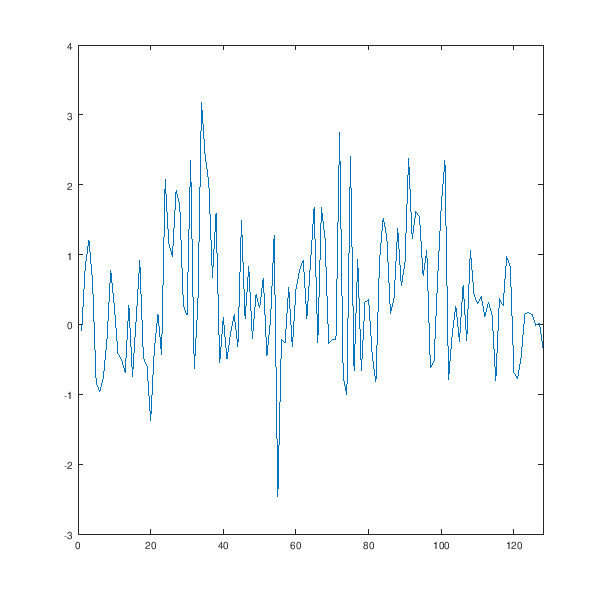}
			\caption{}
			
		\end{subfigure}
		\newline
		\begin{subfigure}{0.24\textwidth}
			\includegraphics[width=45mm, height=45mm]
			{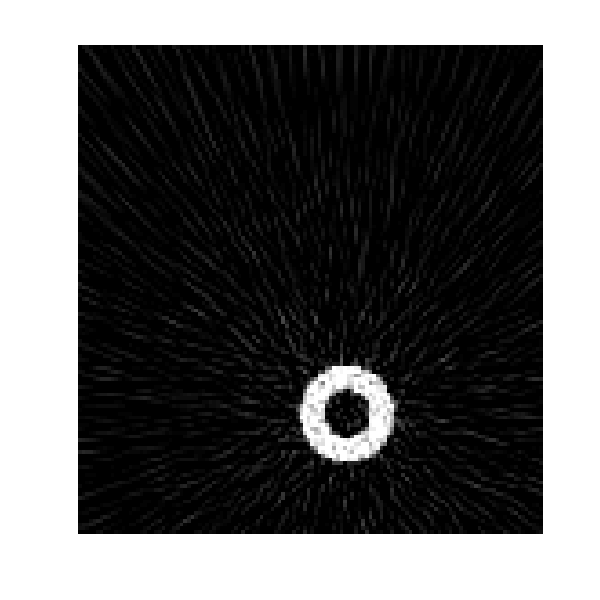}
			\caption{}
			
		\end{subfigure}
		\begin{subfigure}{0.24\textwidth}
			\includegraphics[width=45mm, height=45mm]
			{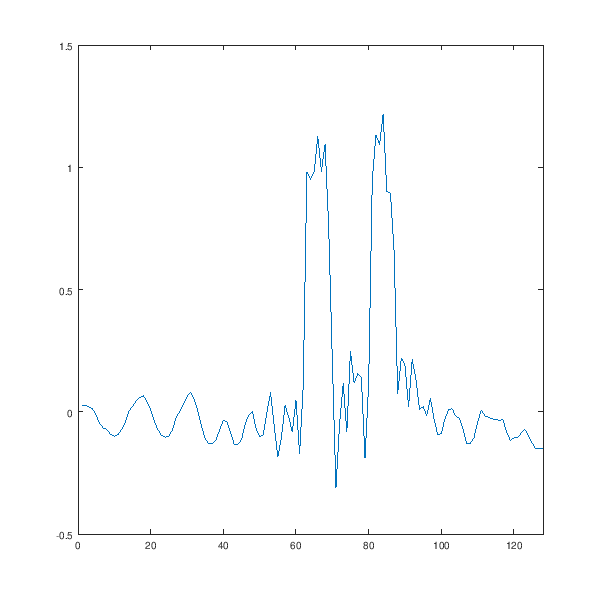}
			\caption{}
			
		\end{subfigure}
		\begin{subfigure}{0.24\textwidth}
			\includegraphics[width=45mm, height=45mm]{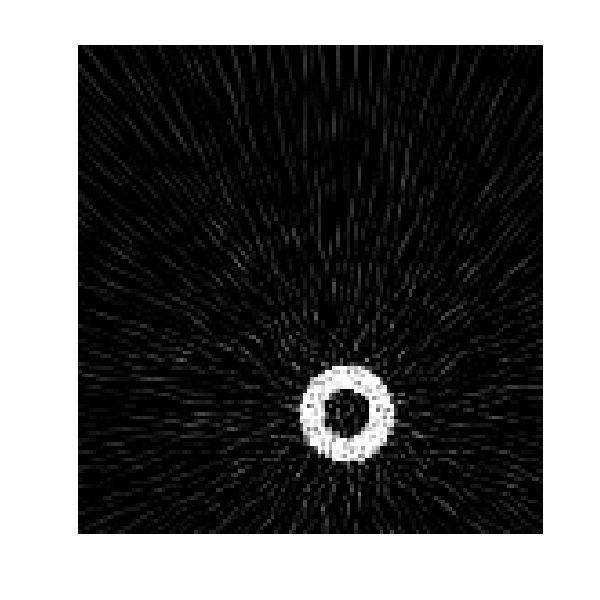}
			\caption{}
			
		\end{subfigure}
		\begin{subfigure}{0.24\textwidth}
			\includegraphics[width=45mm, height=45mm]
			{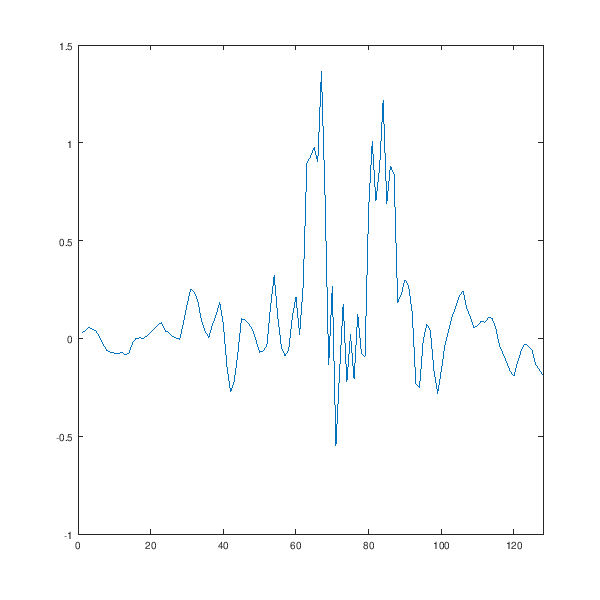}
			\caption{}
			
		\end{subfigure}
		\caption{ weak attenuation, strong noise ($n_{strong}=50$); reconstructions of $f_1, f_2$ using Chang-type methods in 3D (A), (E) and in 2D (C), (G); (B), (D), (F), (H)-sections along $X$-axis}
	\end{figure}

	\begin{figure}[H]
		\begin{subfigure}{0.24\textwidth}
			\includegraphics[width=45mm, height=45mm]{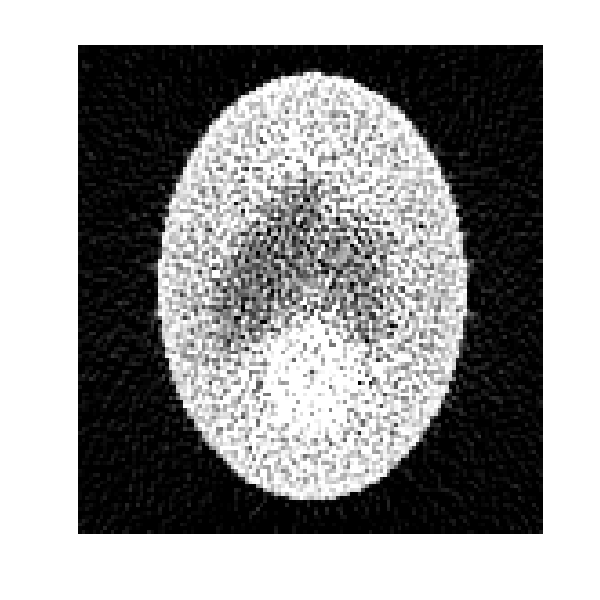}
			\caption{}
			
		\end{subfigure}
		\begin{subfigure}{0.24\textwidth}
			\includegraphics[width=45mm, height=45mm]
			{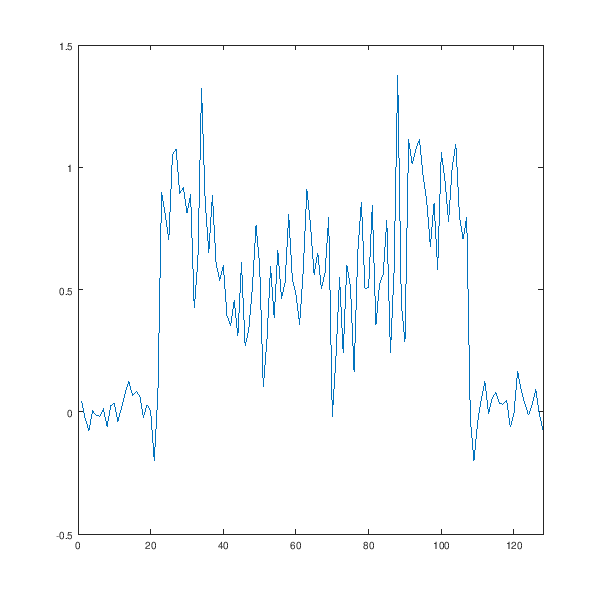}
			\caption{}
			
		\end{subfigure}
		\begin{subfigure}{0.24\textwidth}
			\includegraphics[width=45mm, height=45mm]{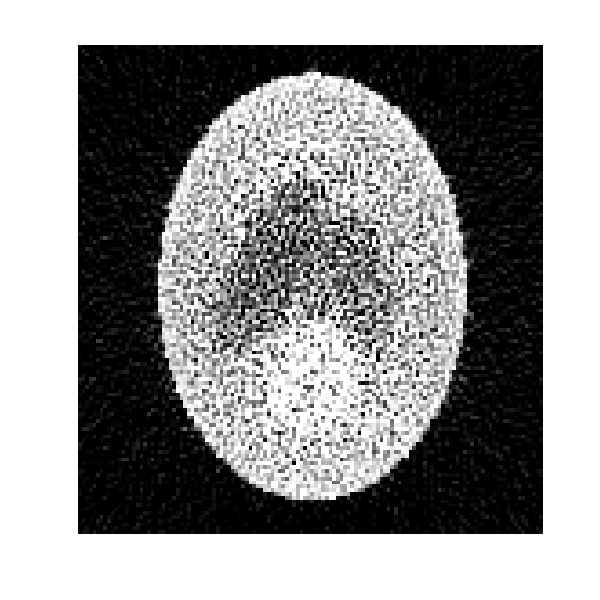}
			\caption{}
			
		\end{subfigure}
		\begin{subfigure}{0.24\textwidth}
			\includegraphics[width=45mm, height=45mm]
			{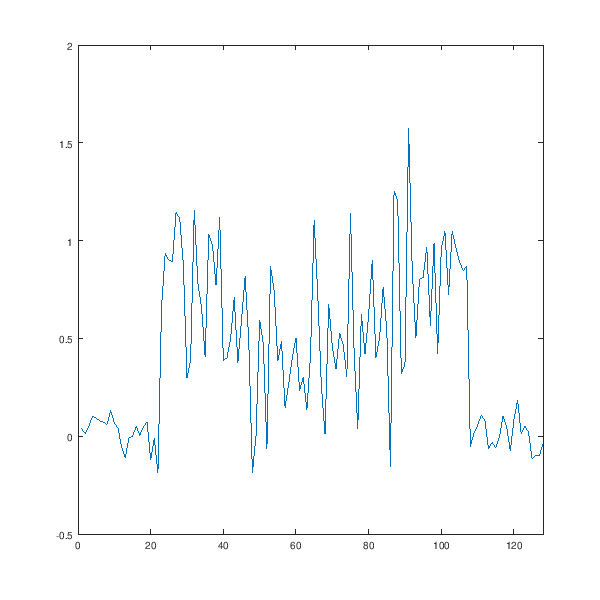}
			\caption{}
			
		\end{subfigure}
		
		\begin{subfigure}{0.24\textwidth}
			\includegraphics[width=45mm, height=45mm]{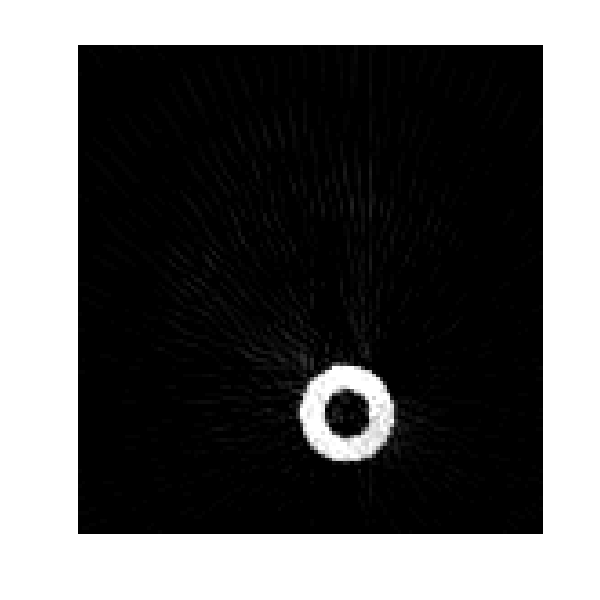}
			\caption{}
			
		\end{subfigure}
		\begin{subfigure}{0.24\textwidth}
			\includegraphics[width=45mm, height=45mm]
			{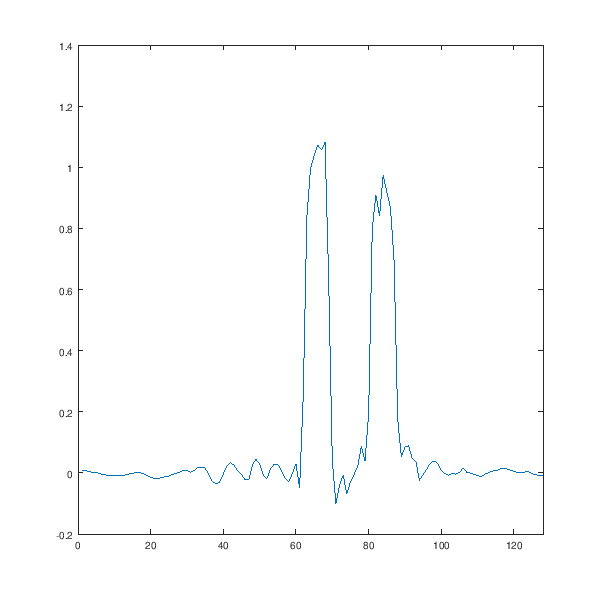}
			\caption{}
			
		\end{subfigure}
		\begin{subfigure}{0.24\textwidth}
			\includegraphics[width=45mm, height=45mm]{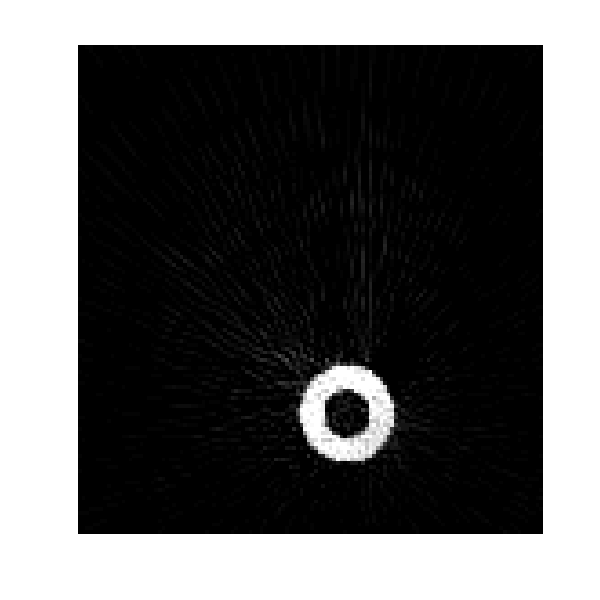}
			\caption{}
			
		\end{subfigure}
		\begin{subfigure}{0.24\textwidth}
			\includegraphics[width=45mm, height=45mm]
			{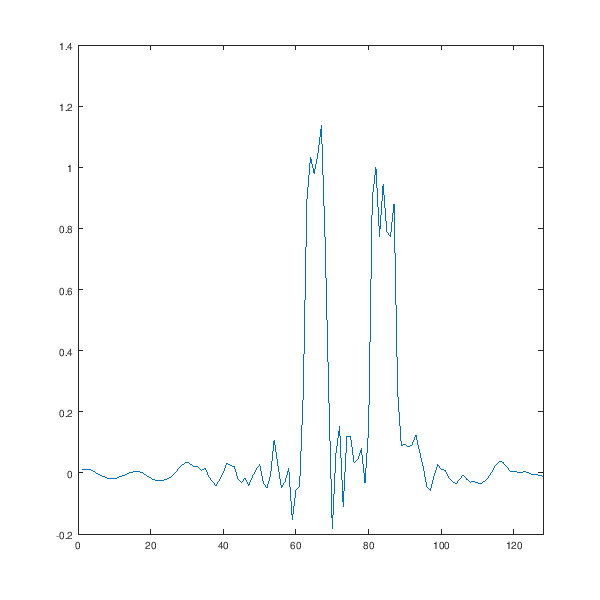}
			\caption{}
			
		\end{subfigure}
		\caption{ strong attenuation, weak noise ($n_{weak}=500$); reconstructions of $f_1, f_2$ using Chang-type formulas in 3D (A), (E) and in 2D (C), (G); (B), (D), (F), (H) -- sections along $X$-axis}
	\end{figure}

	\begin{figure}[H]
		\begin{subfigure}{0.24\textwidth}
			\includegraphics[width=45mm, height=45mm]{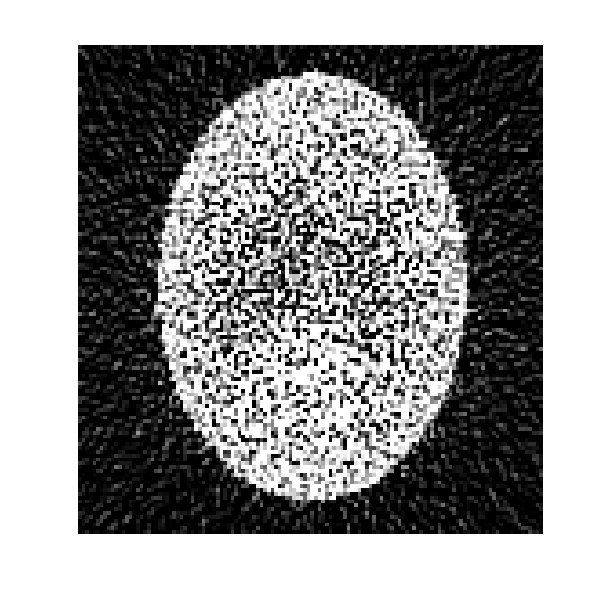}
			\caption{}
			
		\end{subfigure}
		\begin{subfigure}{0.24\textwidth}
			\includegraphics[width=45mm, height=45mm]
			{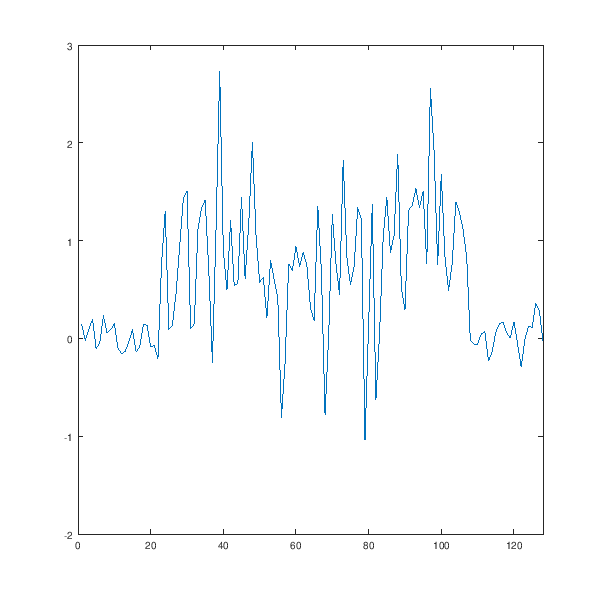}
			\caption{}
			
		\end{subfigure}
		\begin{subfigure}{0.24\textwidth}
			\includegraphics[width=45mm, height=45mm]{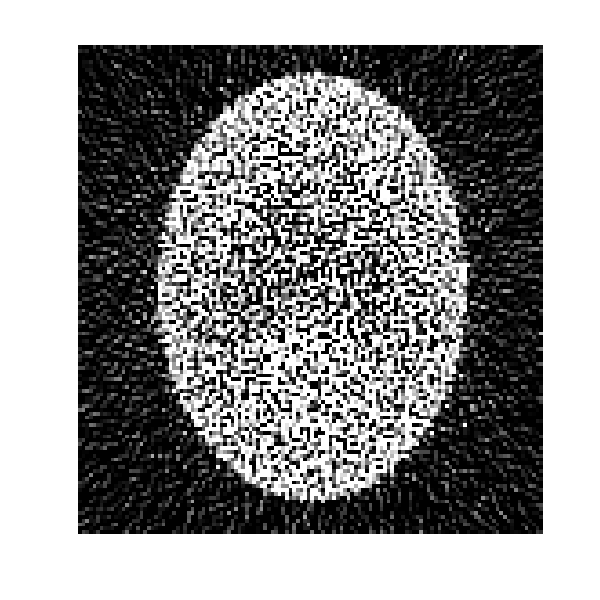}
			\caption{}
			
		\end{subfigure}
		\begin{subfigure}{0.24\textwidth}
			\includegraphics[width=45mm, height=45mm]
			{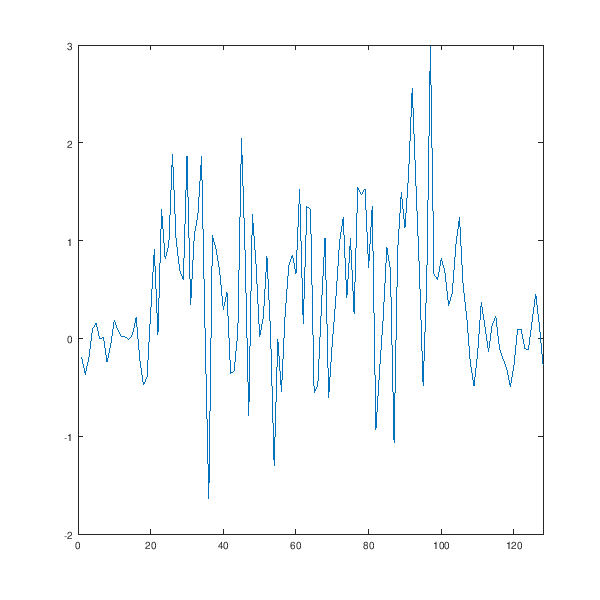}
			\caption{}
			
		\end{subfigure}
		
		\begin{subfigure}{0.24\textwidth}
			\includegraphics[width=45mm, height=45mm]{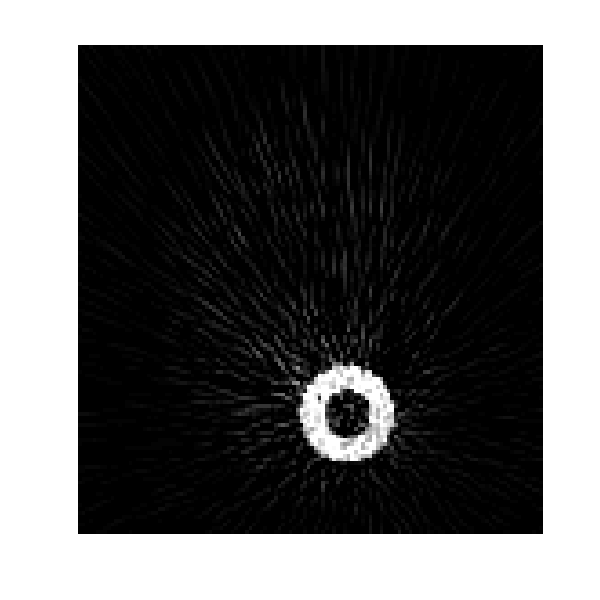}
			\caption{}
			
		\end{subfigure}
		\begin{subfigure}{0.24\textwidth}
			\includegraphics[width=45mm, height=45mm]
			{ph1_att_strong_nweak_chang3d_slice}
			\caption{}
			
		\end{subfigure}
		\begin{subfigure}{0.24\textwidth}
			\includegraphics[width=45mm, height=45mm]{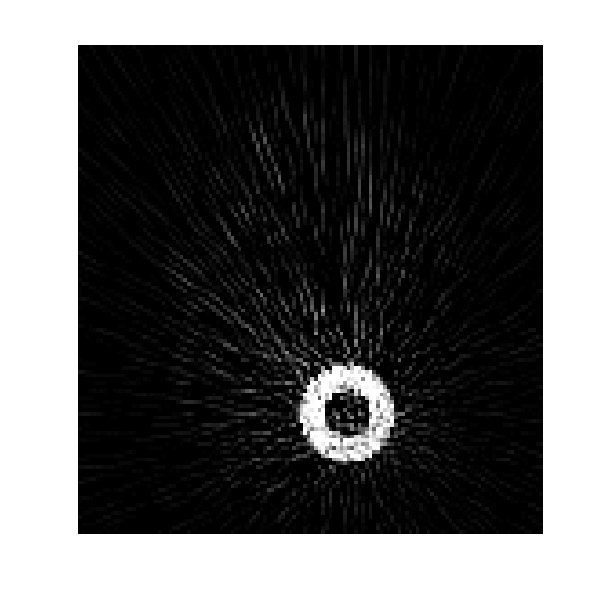}
			\caption{}
			
		\end{subfigure}
		\begin{subfigure}{0.24\textwidth}
			\includegraphics[width=45mm, height=45mm]
			{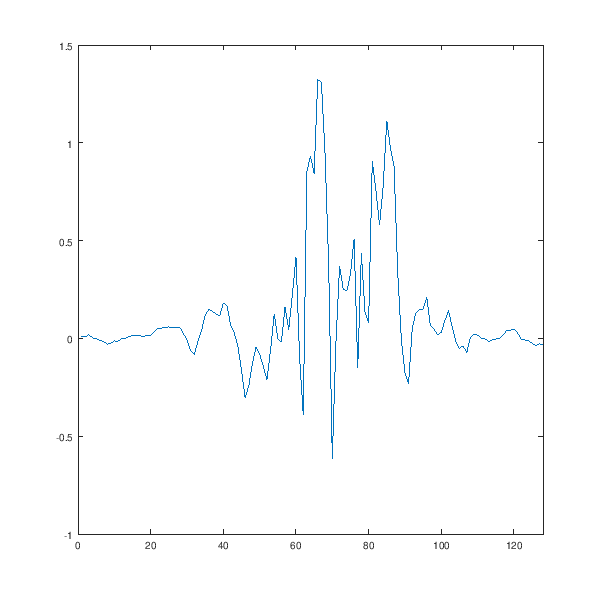}
			\caption{}
			
		\end{subfigure}
		\caption{ strong attenuation, strong noise ($n_{strong}=50$); reconstructions of $f_1, f_2$ using Chang-type methods in 3D (A), (E) and in 2D (C), (G); (B), (D), (F), (H) -- sections along $X$-axis}
		\label{fig:chang.strong.att.strong.noise}
	\end{figure}
	
	Note that in Figures~\ref{fig:chang.strong.att.nonoise} (A), (C) for te noiseless case there is almost no difference between reconstructions using two-dimensional or three-dimensional Chang-type methods. 
	
	At the same time, in  Figures~\ref{fig:chang.weak.att.weak.noise}-\ref{fig:chang.strong.att.strong.noise} one can see that reconstructions obtained using Chang-type methods in 3D look already less noisy than their analogs in 2D. To measure the impact of noise, we computed the relative squared distance between reconstructions from noisy and noise-free data.
	
	Let 
	\begin{align}
	\begin{split}
	&F^{a_j}_{i}, \, i=1,2, \, j = 1, 2,  
	\text{ be the reconstructions of } f_1, f_2 \text{ for}\\
	&\text{strong and weak attenuation levels } a_1, a_2 \text{ without noise and}\\
	&\text{reduced to plane } z=0 \text{ (see Figures~\ref{fig:chang.weak.att.nonoise},~\ref{fig:chang.strong.att.nonoise})}.
	\end{split}
	\end{align}
	Reconstructions from the data with noise will be denoted as follows:
	\begin{align}
	\begin{split}
	&f^{a_j, \, n_k}_{i}, \, i=1,2, \, j = 1, 2,  
	\text{ be the reconstructions of } f_1, f_2 \text{ for strong}\\
	&\text{and weak attenuation levels } a_1, a_2, \text{ noise levels }
	n_1=n_{strong}, n_2=n_{weak}, \text{ respectivley,}\\
	&\text{and reduced to plane } z=0 \text{ (see Figures  \ref{fig:chang.weak.att.weak.noise}-\ref{fig:chang.strong.att.strong.noise})}.
	\end{split}
	\end{align}
	Reconstruction errors $\varepsilon_{f_i, a_j, n_k}$ are defined by the formula
	\begin{eqnarray}\label{eq:chang.error.formula}
	\varepsilon_{f_i, a_j, n_k} = \dfrac{\|f^{a_j, n_k}_{i} - F^{a_j}_{i}\|_2}
	{\|\tilde{f}^{a_j}_{i}\|_2}, \, i = 1,2, \, j = 1,2, \, k = 1,2, \,
	\end{eqnarray}
	where $\|\cdot\|_2$ denotes the Frobenius norm of two-dimensional images seen as matrices of size $N \times N$, where $N$ is the number of pixels per  dimension. The reason to define errors $\varepsilon_{f_i, a_j, n_k}$ as in  \eqref{eq:chang.error.formula} is that Chang-type formulas provide only approximate reconstructions. To measure the effect of noise one must compare approximate reconstructions from noisy data only with approximate reconstructions from the noise-free data.\\
	
	\begin{table}[H]
		\begin{center}
			{\renewcommand{\arraystretch}{1.4}
				\begin{tabular}{| c | c | c | c | c | c | c | c | c | } 
					\hline
					Method / Error & $\varepsilon_{f_1,a_1,n_1}$ & $\varepsilon_{f_1, a_2,n_1}$ & $\varepsilon_{f_1, a_1,n_2}$ & $\varepsilon_{f_1, a_2,n_2}$ 
					& $\varepsilon_{f_2, a_1,n_1}$ & $\varepsilon_{f_2,a_2,n_1}$ & $\varepsilon_{f_2, a_1,n_2}$ & $\varepsilon_{f_2,a_2,n_2}$ \\
					\hline
					2D-method & 1.193 & 1.340 & 0.377 & 0.434 & 0.644 & 0.625 & 0.211 & 0.202 \\
					\hline
					3D-method & 0.779 & 0.942 & 0.251 & 0.299 & 0.438 & 0.432 & 0.137 & 0.135 \\
					\hline
			\end{tabular}}
			\caption{relative errors for Chang-type method in 2D and 3D}
			\label{tab:errors.chang}
		\end{center}
	\end{table}
	From Table~\ref{tab:errors.chang} one can see that three-dimensional Chang-type method with preprocessing outperforms its two-dimensional analog for all cases of activity phantoms, attenuations and noise levels. Moreover, the gain of stability in reconstructions is already visible from Figures~\ref{fig:chang.weak.att.weak.noise}-\ref{fig:chang.strong.att.strong.noise}.

	\subsection{Reconstructions via Kunyansky-type methods}
	\label{subsect.numerical.inversion.iterative}
	For Kunyansky-type methods in two and three dimensions (see also Subsection~\ref{subsect:reconstruct.methods})
	we use parameter value $m=1$. Our choice for $m=1$ is motivated by condition  \eqref{eq:reconstruction-methods.kunyansky.conv-condition}. More precisely, for strong attenuation $a_1$, condition \eqref{eq:reconstruction-methods.kunyansky.conv-condition} is barely satisfied ($\sigma_{W,\, D, \, 1} = 0.89$ for $d=3$ and $0.52$ for $d=2$) therefore, to be able to compare two-dimensional and three-dimensional reconstructions, we keep $m=1$. Finally, for weak attenuation $a_2$, condition   \eqref{eq:reconstruction-methods.kunyansky.conv-condition} is efficiently satisfied ($\sigma_{W, \, D, \, 1} = 0.17$ for $d=3$ and $0.11$ for $d=2$), however, in this case Chang-type methods ($m=0$) give already  almost-perfect reconstructions of $f_1, \, f_2$ (modulo the noise) (see Figures~\ref{fig:chang.weak.att.nonoise} (A), (C), (E), (G)).
	
	\subsubsection*{Noiseless case}
	
	\begin{figure}[H]
		\begin{subfigure}{0.24\textwidth}
			\includegraphics[width=45mm, height=40mm]{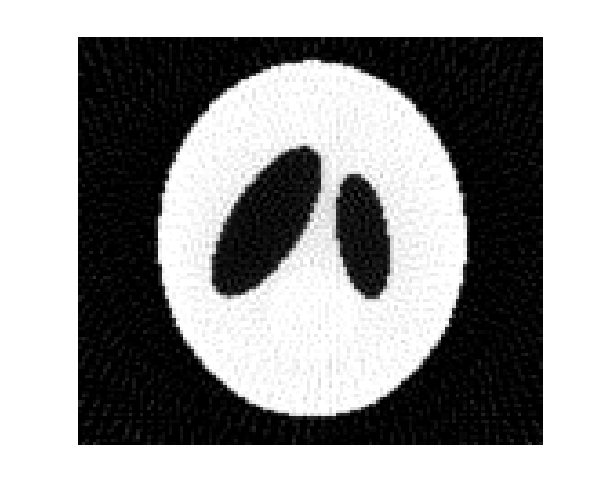}
			\caption{}
			
		\end{subfigure}
		\begin{subfigure}{0.24\textwidth}
			\includegraphics[width=45mm, height=40mm]
			{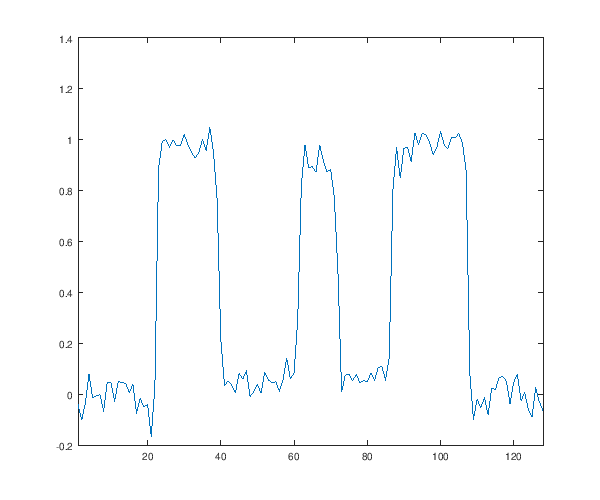}
			\caption{}
			
		\end{subfigure}
		\begin{subfigure}{0.24\textwidth}
			\includegraphics[width=45mm, height=40mm]{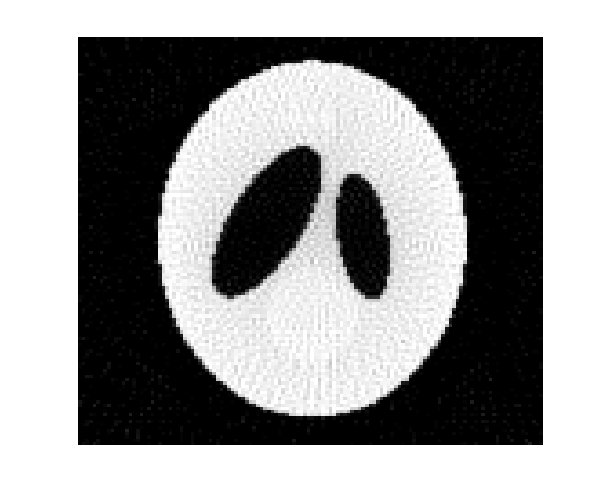}
			\caption{}
			
		\end{subfigure}
		\begin{subfigure}{0.24\textwidth}
			\includegraphics[width=45mm, height=40mm]
			{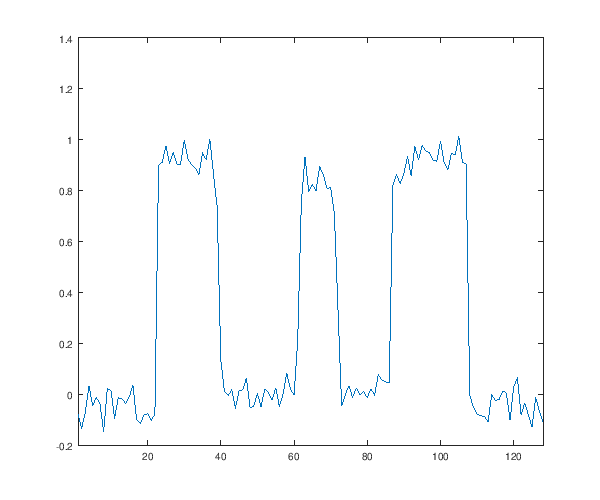}
			\caption{}
			
		\end{subfigure}
		
		\begin{subfigure}{0.24\textwidth}
			\includegraphics[width=45mm, height=40mm]
			{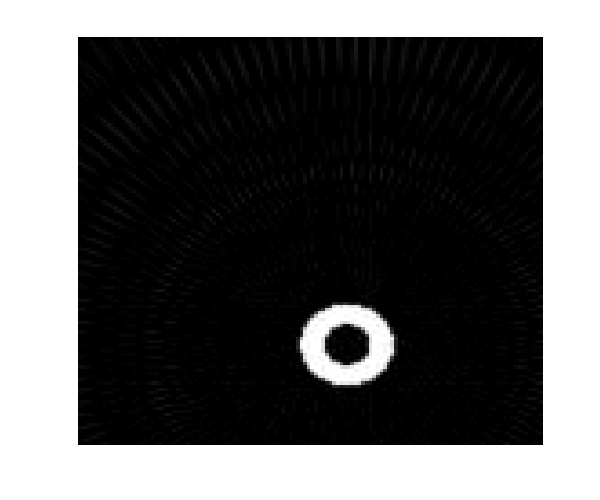}
			\caption{}
			
		\end{subfigure}
		\begin{subfigure}{0.24\textwidth}
			\includegraphics[width=45mm, height=40mm]
			{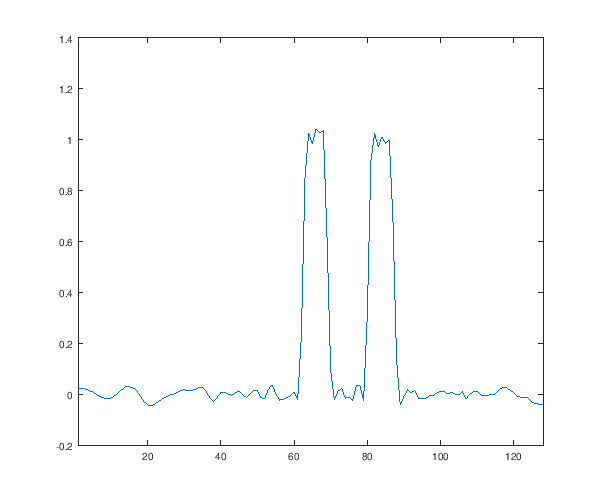}
			\caption{}
			
		\end{subfigure}
		\begin{subfigure}{0.24\textwidth}
			\includegraphics[width=45mm, height=40mm]{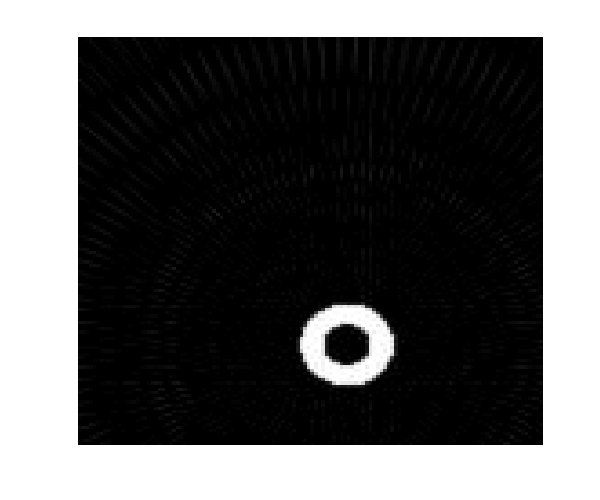}
			\caption{}
			
		\end{subfigure}
		\begin{subfigure}{0.24\textwidth}
			\includegraphics[width=45mm, height=40mm]
			{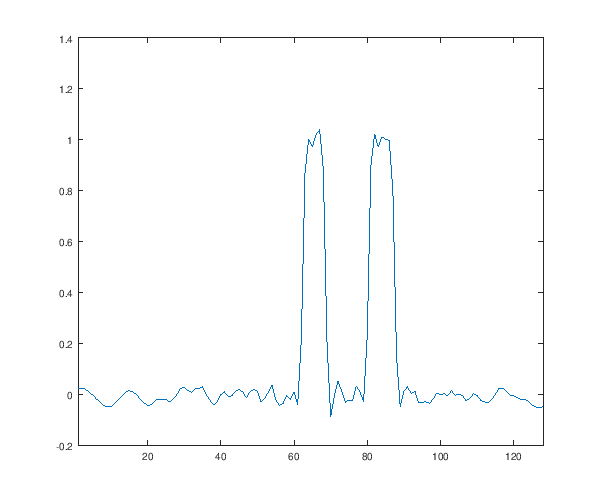}
			\caption{}
			
		\end{subfigure}
		\caption{ weak attenuation $a_2$, no noise; reconstructions of $f_1, f_2$ using iterative Kunyansky-type methods in 3D (A), (E) and in 2D (C), (G); (B), (D), (F), (H) -- sections along $X$-axis}
		\label{fig:iterative.weak.att.nonoise}
	\end{figure}

	\begin{figure}[H]
		\begin{subfigure}{0.24\textwidth}
			\includegraphics[width=45mm, height=40mm]{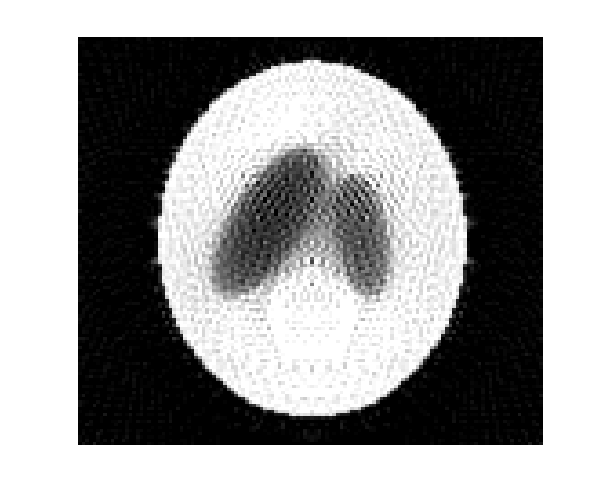}
			\caption{}
			
		\end{subfigure}
		\begin{subfigure}{0.24\textwidth}
			\includegraphics[width=45mm, height=40mm]
			{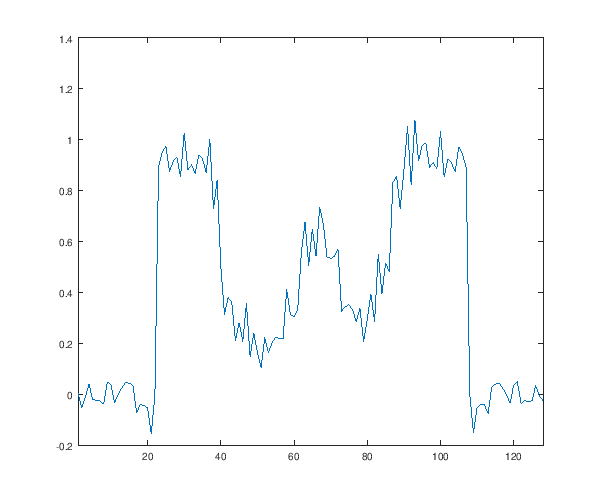}
			\caption{}
			
		\end{subfigure}
		\begin{subfigure}{0.24\textwidth}
			\includegraphics[width=45mm, height=40mm]{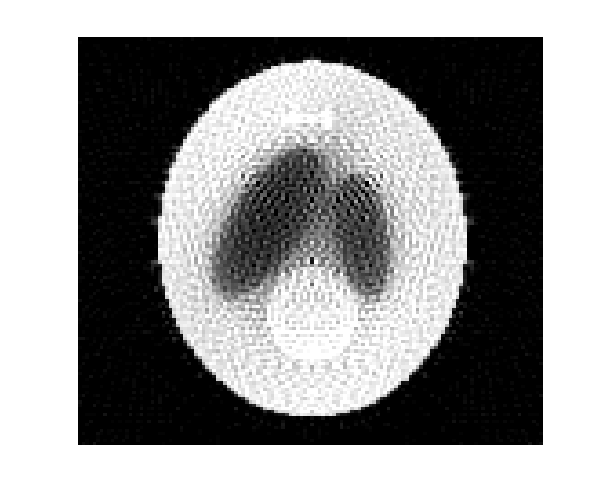}
			\caption{}
			
		\end{subfigure}
		\begin{subfigure}{0.24\textwidth}
			\includegraphics[width=45mm, height=40mm]
			{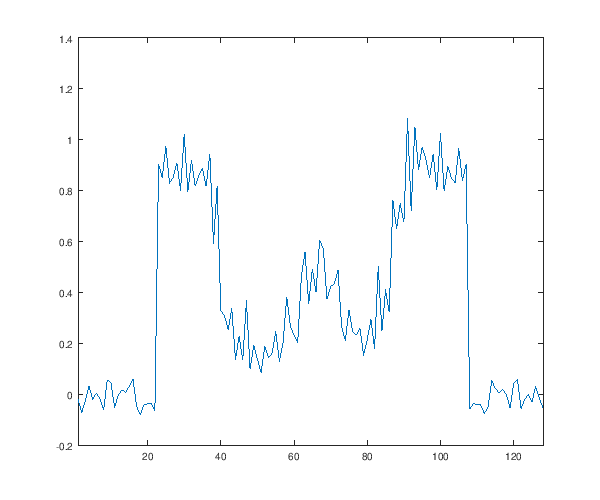}
			\caption{}
			
		\end{subfigure}
		
		\begin{subfigure}{0.24\textwidth}
			\includegraphics[width=45mm, height=40mm]
			{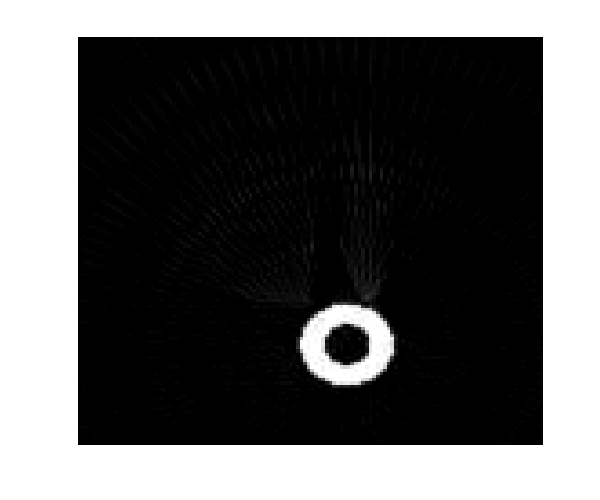}
			\caption{}
			
		\end{subfigure}
		\begin{subfigure}{0.24\textwidth}
			\includegraphics[width=45mm, height=40mm]
			{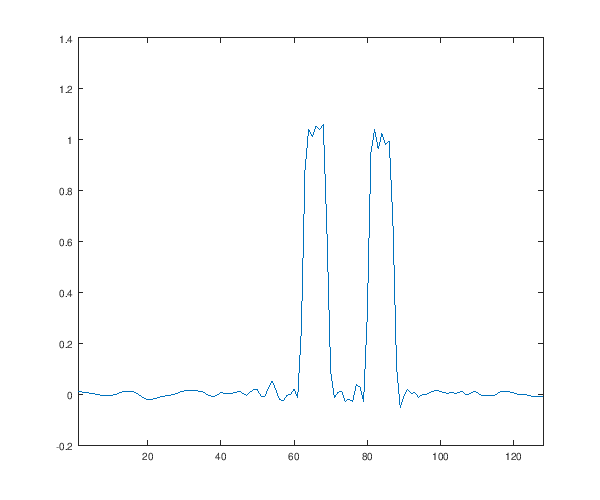}
			\caption{}
			
		\end{subfigure}
		\begin{subfigure}{0.24\textwidth}
			\includegraphics[width=45mm, height=40mm]{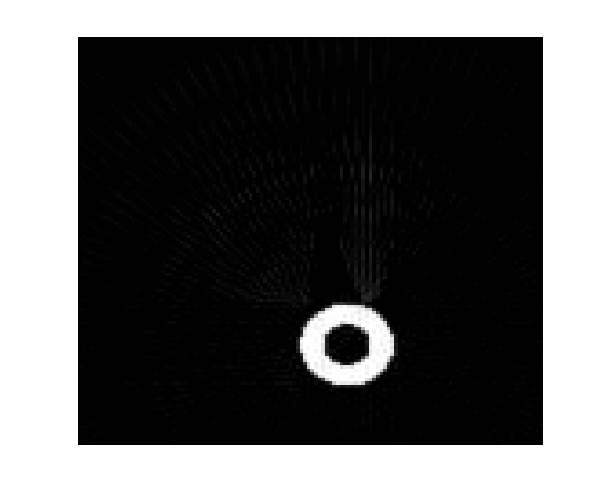}
			\caption{}
			
		\end{subfigure}
		\begin{subfigure}{0.24\textwidth}
			\includegraphics[width=45mm, height=40mm]
			{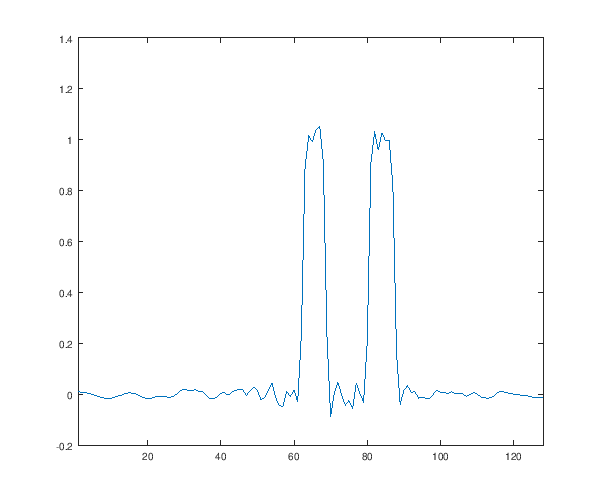}
			\caption{}
			
		\end{subfigure}
		\caption{ strong attenuation $a_1$, no noise; reconstructions of $f_1, f_2$ using iterative Kunyansky-type methods in 3D (A), (E) and in 2D (C), (G); (B), (D), (F), (H) -- sections along $X$-axis}
		\label{fig:iterative.strong.att.nonoise}
	\end{figure}
	
	\subsubsection*{Case with noise}
	
	\begin{figure}[H]
		\begin{subfigure}{0.24\textwidth}
			\includegraphics[width=45mm, height=45mm]{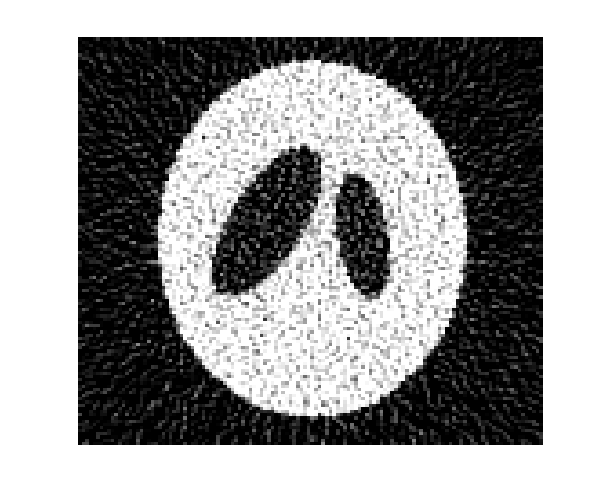}
			\caption{}
			
		\end{subfigure}
		\begin{subfigure}{0.24\textwidth}
			\includegraphics[width=45mm, height=45mm]
			{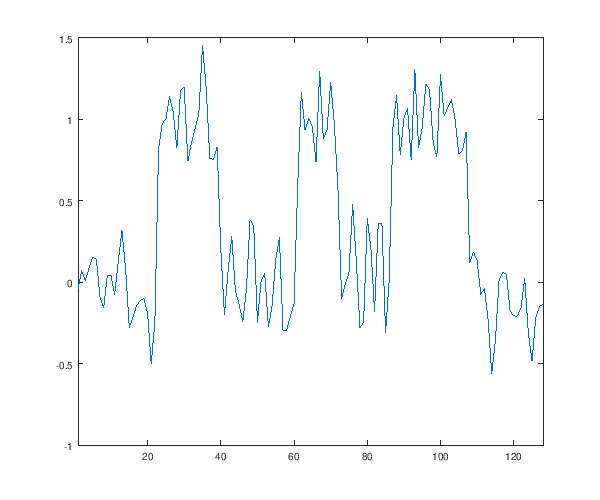}
			\caption{}
			
		\end{subfigure}
		\begin{subfigure}{0.24\textwidth}
			\includegraphics[width=45mm, height=45mm]{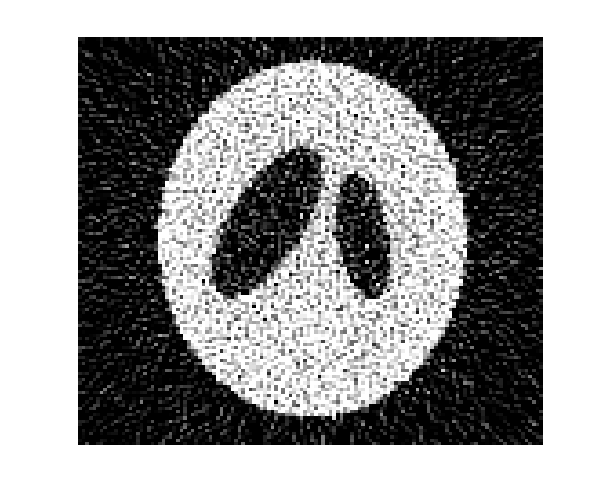}
			\caption{}
			
		\end{subfigure}
		\begin{subfigure}{0.24\textwidth}
			\includegraphics[width=45mm, height=45mm]
			{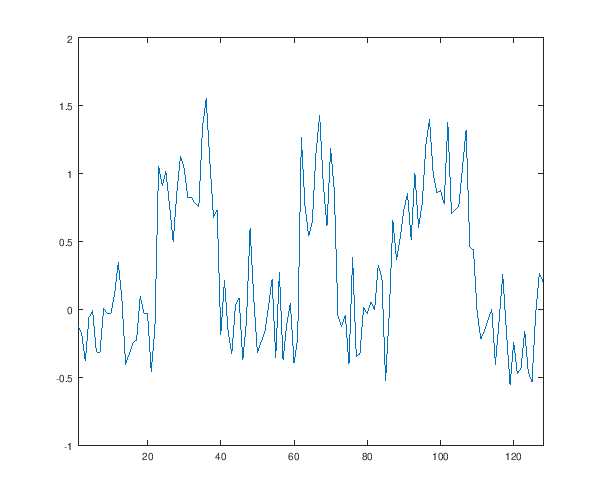}
			\caption{}
			
		\end{subfigure}
		\newline
		\begin{subfigure}{0.24\textwidth}
			\includegraphics[width=45mm, height=45mm]
			{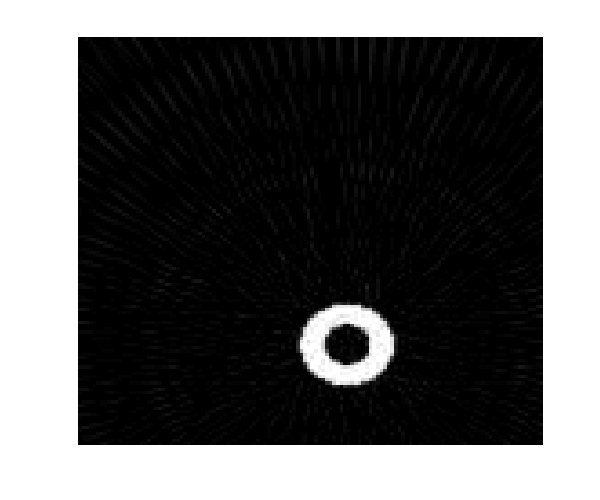}
			\caption{}
			
		\end{subfigure}
		\begin{subfigure}{0.24\textwidth}
			\includegraphics[width=45mm, height=45mm]
			{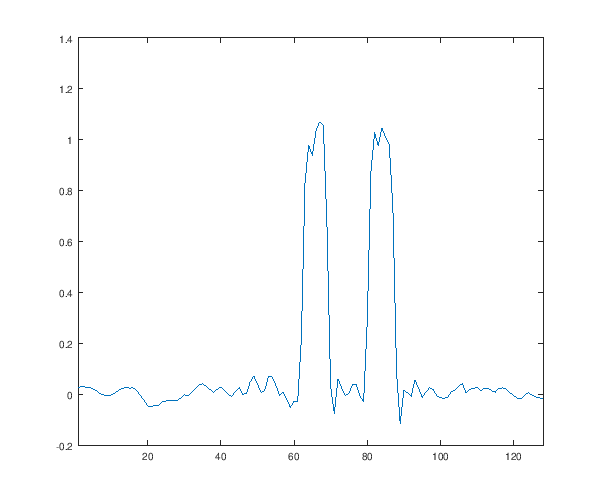}
			\caption{}
			
		\end{subfigure}
		\begin{subfigure}{0.24\textwidth}
			\includegraphics[width=45mm, height=45mm]{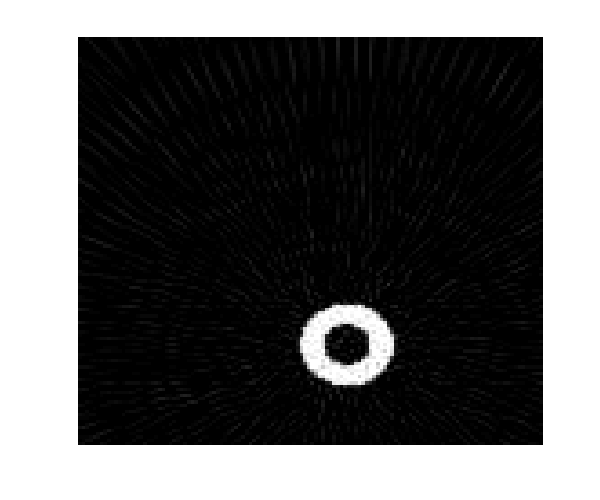}
			\caption{}
			
		\end{subfigure}
		\begin{subfigure}{0.24\textwidth}
			\includegraphics[width=45mm, height=45mm]
			{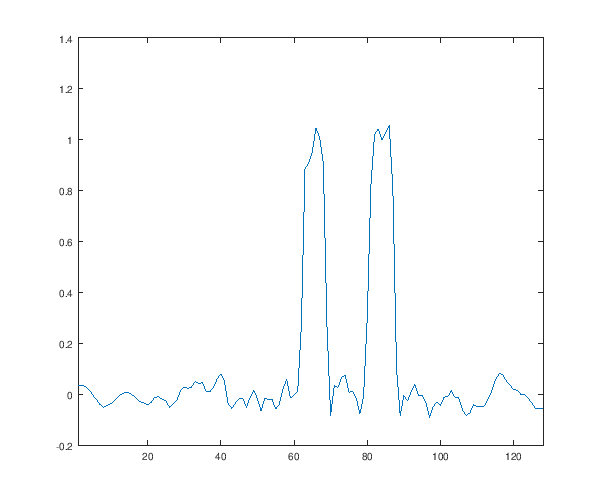}
			\caption{}
			
		\end{subfigure}
		
		\caption{ weak attenuation $a_2$, weak noise ($n_{weak}=500$); reconstructions using Kunyansky-type methods in 3D (A), (E) and in 2D (C), (G); (B), (D), (F), (H) -- sections along X-axis}
		\label{fig:iterative.weak.att.weak.noise}
	\end{figure}
	
	
	\begin{figure}[H]
		\begin{subfigure}{0.24\textwidth}
			\includegraphics[width=45mm, height=45mm]{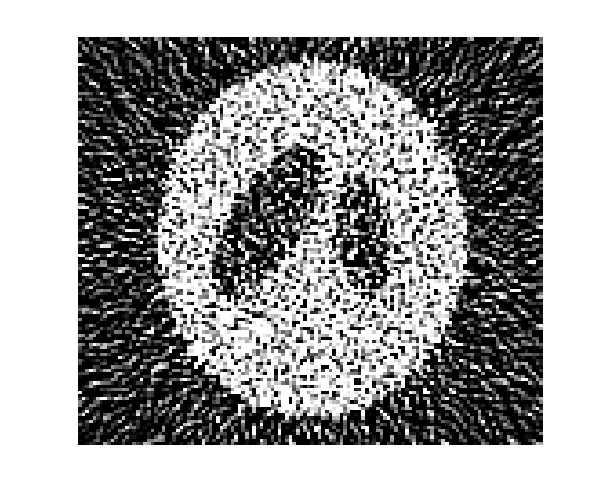}
			\caption{}
			
		\end{subfigure}
		\begin{subfigure}{0.24\textwidth}
			\includegraphics[width=45mm, height=45mm]
			{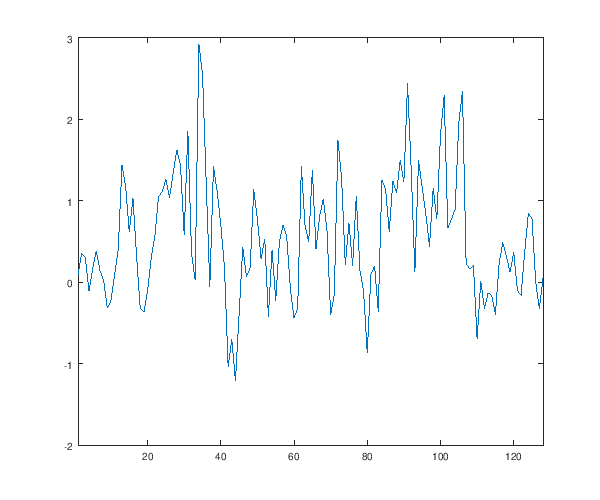}
			\caption{}
			
		\end{subfigure}
		\begin{subfigure}{0.24\textwidth}
			\includegraphics[width=45mm, height=45mm]{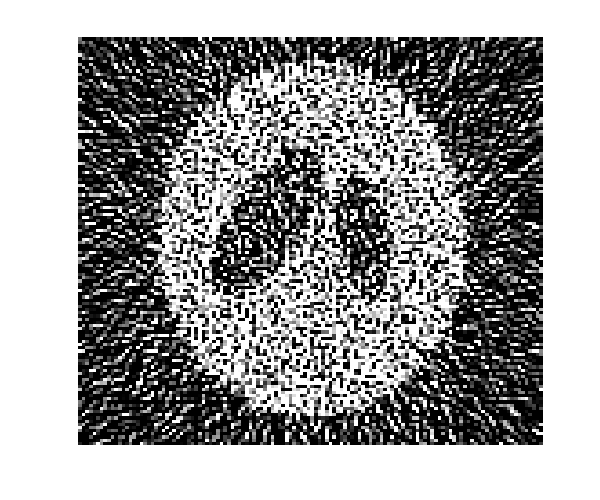}
			\caption{}
			
		\end{subfigure}
		\begin{subfigure}{0.24\textwidth}
			\includegraphics[width=45mm, height=45mm]
			{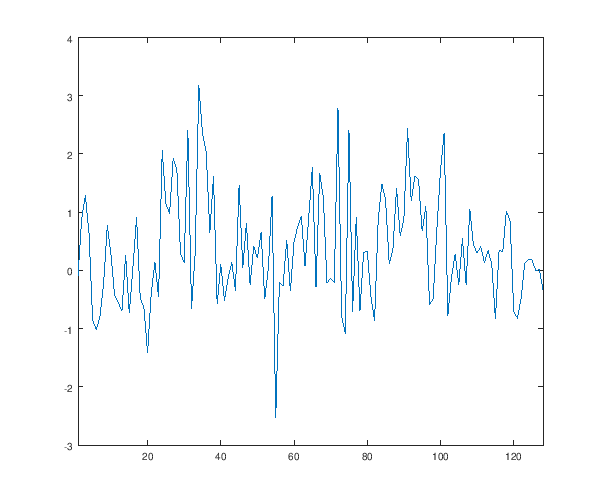}
			\caption{}
			
		\end{subfigure}
		\newline
		\begin{subfigure}{0.24\textwidth}
			\includegraphics[width=45mm, height=45mm]
			{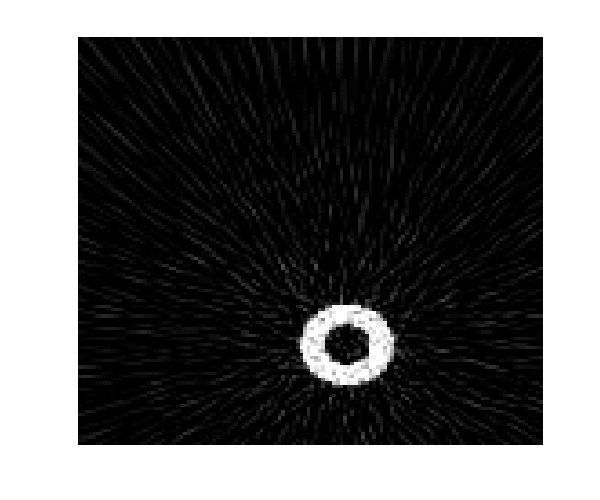}
			\caption{}
			
		\end{subfigure}
		\begin{subfigure}{0.24\textwidth}
			\includegraphics[width=45mm, height=45mm]
			{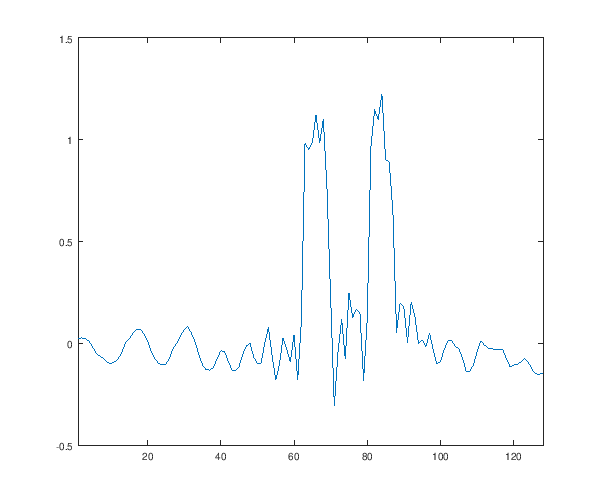}
			\caption{}
			
		\end{subfigure}
		\begin{subfigure}{0.24\textwidth}
			\includegraphics[width=45mm, height=45mm]{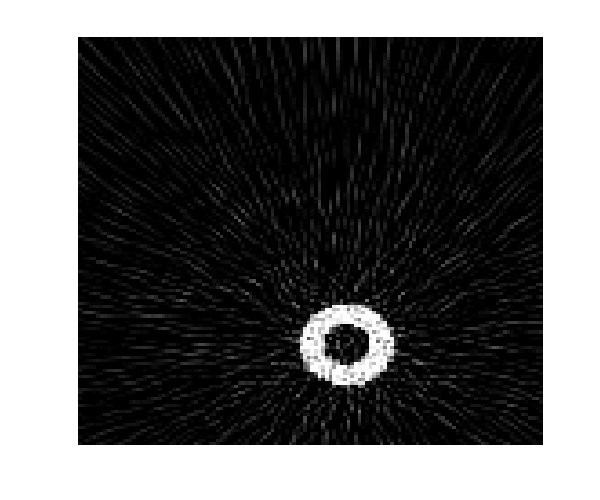}
			\caption{}
			
		\end{subfigure}
		\begin{subfigure}{0.24\textwidth}
			\includegraphics[width=45mm, height=45mm]
			{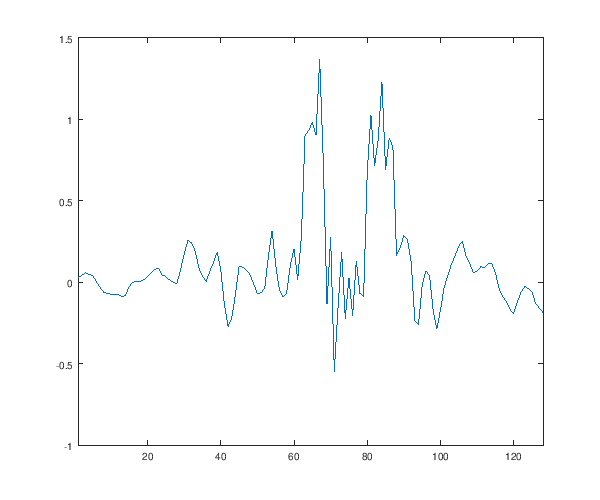}
			\caption{}
			
		\end{subfigure}
		\caption{ weak attenuation $a_2$, strong noise ($n_{strong}=50$); reconstructions using Kunyansky-type methods in 3D (A), (E) and in 2D (C), (G); (B), (D), (F), (H) -- sections along $X$-axis}
	\end{figure}

	\begin{figure}[H]
		\begin{subfigure}{0.24\textwidth}
			\includegraphics[width=45mm, height=45mm]{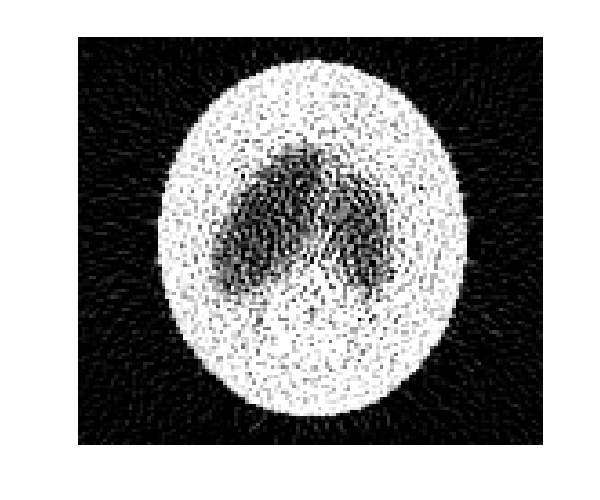}
			\caption{}
			
		\end{subfigure}
		\begin{subfigure}{0.24\textwidth}
			\includegraphics[width=45mm, height=45mm]
			{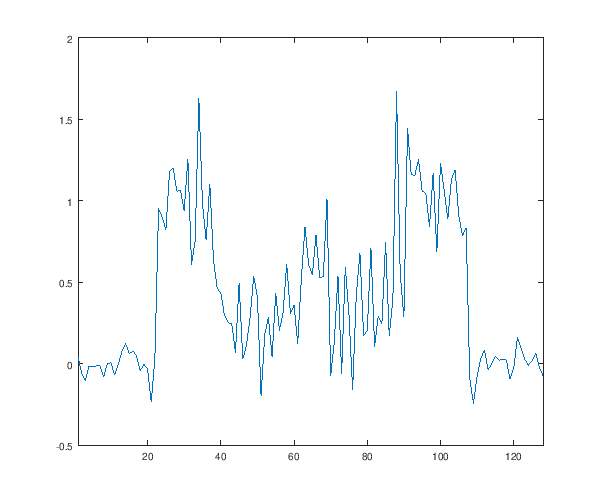}
			\caption{}
			
		\end{subfigure}
		\begin{subfigure}{0.24\textwidth}
			\includegraphics[width=45mm, height=45mm]{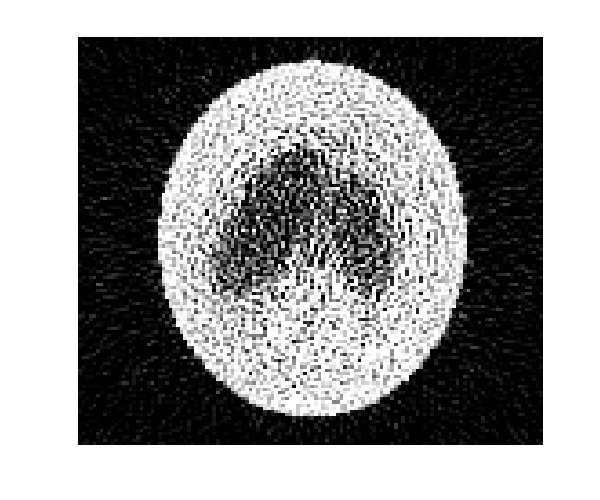}
			\caption{}
			
		\end{subfigure}
		\begin{subfigure}{0.24\textwidth}
			\includegraphics[width=45mm, height=45mm]
			{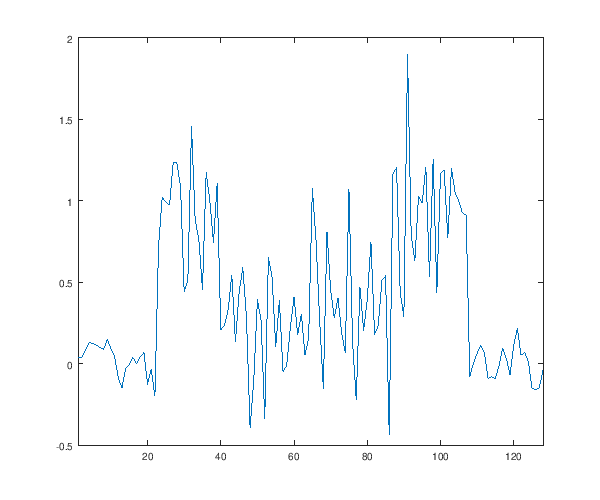}
			\caption{}
			
		\end{subfigure}
		
		\begin{subfigure}{0.24\textwidth}
			\includegraphics[width=45mm, height=45mm]{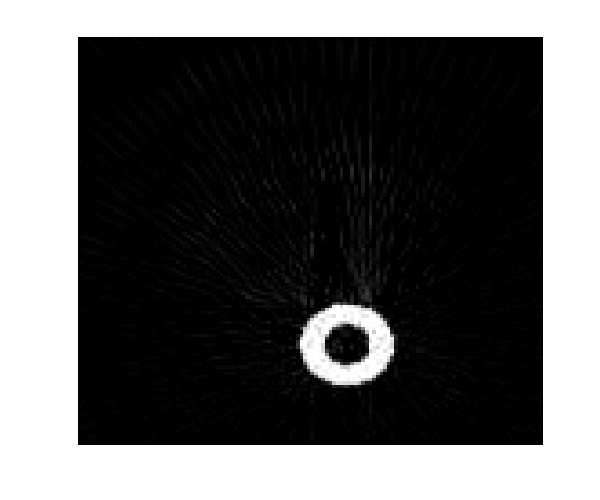}
			\caption{}
			
		\end{subfigure}
		\begin{subfigure}{0.24\textwidth}
			\includegraphics[width=45mm, height=45mm]
			{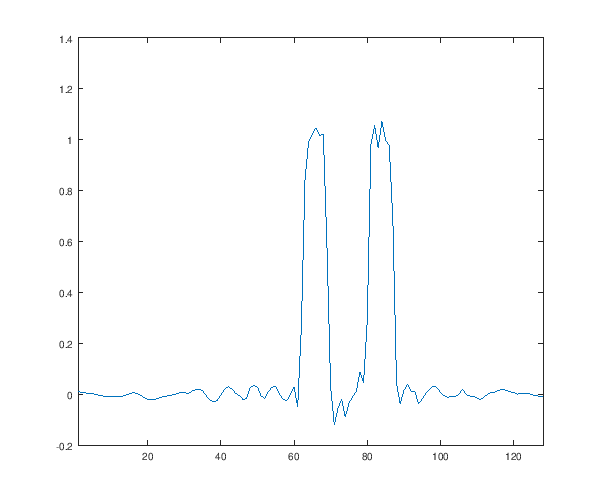}
			\caption{}
			
		\end{subfigure}
		\begin{subfigure}{0.24\textwidth}
			\includegraphics[width=45mm, height=45mm]{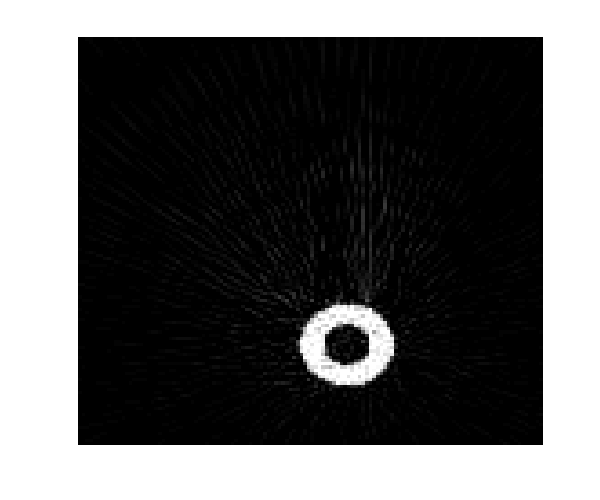}
			\caption{}
			
		\end{subfigure}
		\begin{subfigure}{0.24\textwidth}
			\includegraphics[width=45mm, height=45mm]
			{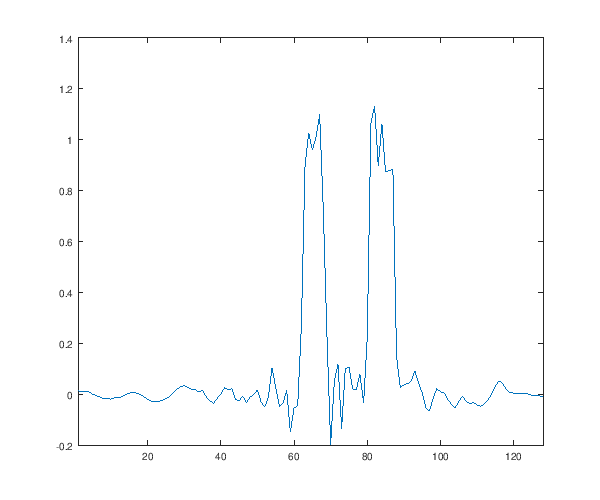}
			\caption{}
			
		\end{subfigure}
		\caption{ strong attenuation $a_1$, weak noise ($n_{weak}=500$); reconstructions using Kunyansky-type methods in 3D (A), (E) and in 2D (C), (G); (B), (D), (F), (H) -- sections along $X$-axis}
	\end{figure}

	\begin{figure}[H]
		\begin{subfigure}{0.24\textwidth}
			\includegraphics[width=45mm, height=45mm]{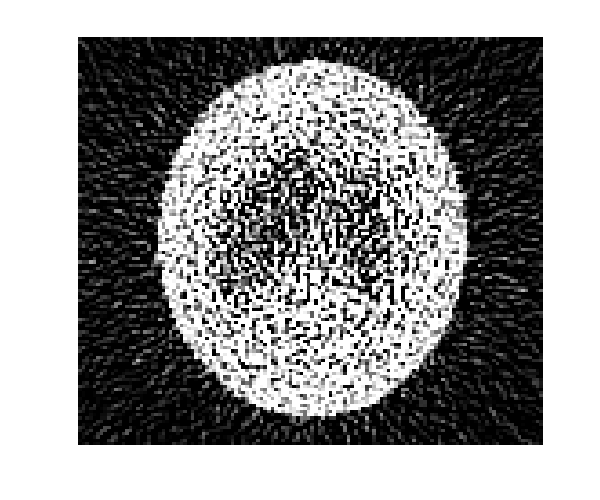}
			\caption{}
			
		\end{subfigure}
		\begin{subfigure}{0.24\textwidth}
			\includegraphics[width=45mm, height=45mm]
			{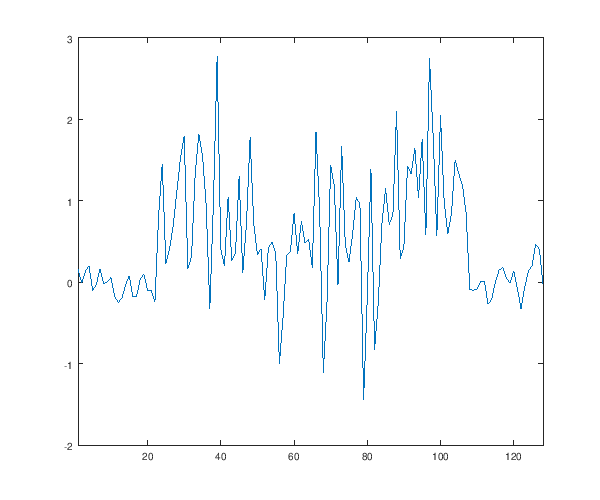}
			\caption{}
			
		\end{subfigure}
		\begin{subfigure}{0.24\textwidth}
			\includegraphics[width=45mm, height=45mm]{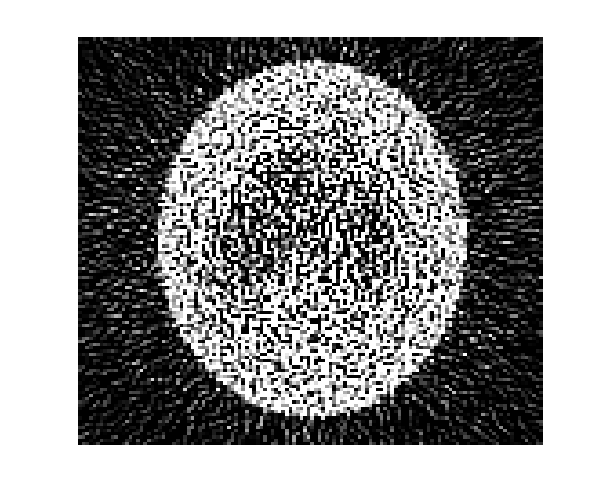}
			\caption{}
			
		\end{subfigure}
		\begin{subfigure}{0.24\textwidth}
			\includegraphics[width=45mm, height=45mm]
			{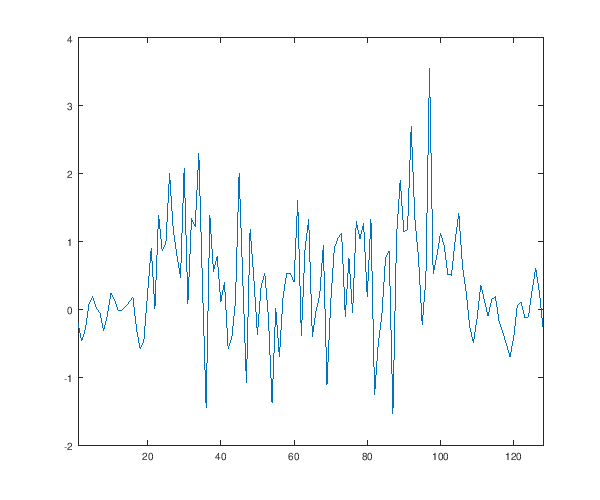}
			\caption{}
			
		\end{subfigure}
		
		\begin{subfigure}{0.24\textwidth}
			\includegraphics[width=45mm, height=45mm]{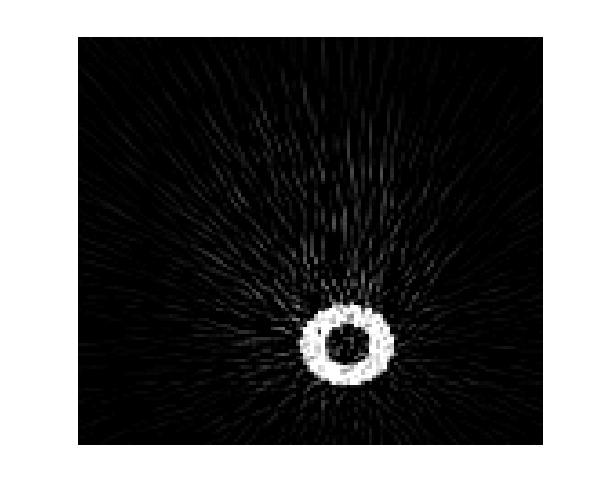}
			\caption{}
			
		\end{subfigure}
		\begin{subfigure}{0.24\textwidth}
			\includegraphics[width=45mm, height=45mm]
			{ph1_att_strong_nweak_iterative3d_slice}
			\caption{}
			
		\end{subfigure}
		\begin{subfigure}{0.24\textwidth}
			\includegraphics[width=45mm, height=45mm]{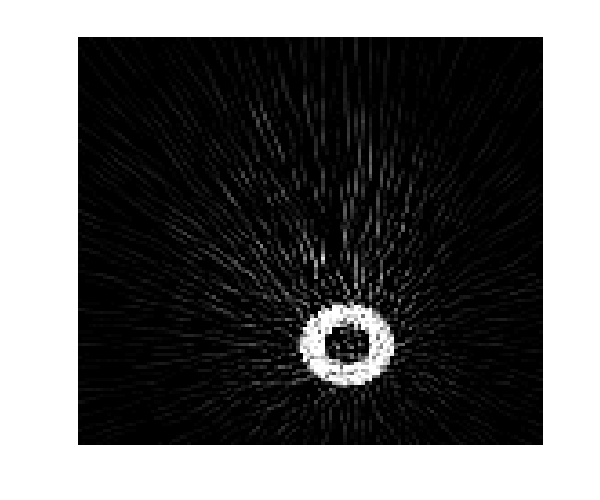}
			\caption{}
			
		\end{subfigure}
		\begin{subfigure}{0.24\textwidth}
			\includegraphics[width=45mm, height=45mm]
			{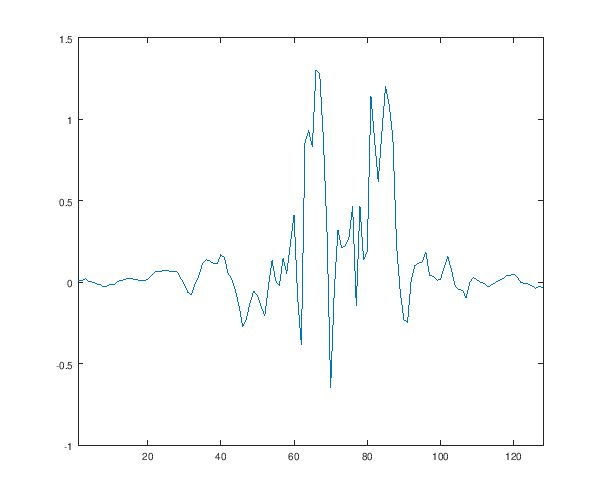}
			\caption{}
			
		\end{subfigure}
		\caption{ strong attenuation $a_1$, strong noise ($n_{strong}=50$); reconstructions using Kunyansky-type methods in 3D (A), (E) and in 2D (C), (G); (B), (D), (F), (H) -- sections along $X$-axis}
		\label{fig:iterative.strong.att.strong.noise}
	\end{figure}
	
	Analogously to the experiment with Chang-type methods, we also compute the relative errors for two-dimensional and three-dimensional reconstructions. 
	\par Let 
	\begin{align}
	\begin{split}
	&F^{a_j}_{i}, \, i=1,2, \, j = 1, 2,  
	\text{ be the reconstructions of } f_1, f_2 \text{ for}\\
	&\text{strong and weak attenuation levels } a_1, a_2 \text{ without noise and}\\
	&\text{reduced to plane } z=0 \text{ (see Figures~\ref{fig:iterative.weak.att.nonoise},~\ref{fig:iterative.strong.att.nonoise})}.
	\end{split}
	\end{align}
	Reconstructions in from noisy data are defined as follows:
	\begin{align}
	\begin{split}
	&f^{a_j, \, n_k}_{i}, \, i=1,2, \, j = 1, 2, \, k = 1,2,  
	\text{ denote the reconstructions of } f_1, f_2\\
	&\text{for strong and weak attenuation levels } a_1, a_2, \\
	&\text{for noise levels } n_1=n_{strong}, n_2=n_{weak}, \text{ respectivley}\\
	&\text{and reduced to plane } z=0 \text{ (see Figures  \ref{fig:iterative.weak.att.weak.noise}-\ref{fig:iterative.strong.att.strong.noise})}.
	\end{split}
	\end{align}
	The reconstruction errors $\varepsilon_{f_i, a_j, n_k}$ are defined by the formula
	\begin{eqnarray}
	\varepsilon_{f_i, a_j, n_k} = \dfrac{\|f^{a_j, n_k}_{i} - F^{a_j}_{i}\|_2}
	{\|F^{a_j}_{i}\|_2}, \, i = 1,2, \, j = 1,2, \, k = 1,2, \,
	\end{eqnarray}
	where $\|\cdot\|_2$ denotes the Frobenius norm of two-dimensional images seen as matrices of size $N \times N$, where $N$ is the number of pixels per  dimension.\\
	
	\begin{table}[H]
		\begin{center}
			{\renewcommand{\arraystretch}{1.4}
				\begin{tabular}{| c | c | c | c | c | c | c | c | c | } 
					\hline
					Method / Error & $\varepsilon_{f_1,a_1,n_1}$ & $\varepsilon_{f_1, a_2,n_1}$ & $\varepsilon_{f_1, a_1,n_2}$ & $\varepsilon_{f_1, a_2,n_2}$ 
					& $\varepsilon_{f_2, a_1,n_1}$ & $\varepsilon_{f_2,a_2,n_1}$ & $\varepsilon_{f_2, a_1,n_2}$ & $\varepsilon_{f_2,a_2,n_2}$ \\
					\hline
					2D-method & 1.279 & 1.415 & 0.437 & 0.438 & 0.714 & 0.634 & 0.254 & 0.205 \\
					\hline
					3D-method & 0.847 & 0.952 & 0.316 & 0.303 & 0.494 & 0.439 & 0.187 & 0.138 \\
					\hline
			\end{tabular}}
			\caption{relative reconstruction errors for Kunyansky-type iterative methods}
			\label{tab:errors.iterative}
		\end{center}
	\end{table}
	From Table~\ref{tab:errors.iterative} one could see that the three-dimensional Kunyansky-type method with preprocessing outperforms its two-dimensional analogue for all cases of activity phantoms, attenuations and noise levels. Moreover, the gain of stability is already visible in Figures~\ref{fig:iterative.weak.att.weak.noise}-\ref{fig:iterative.strong.att.strong.noise}.

	\section{Experiment on real data}\label{sect:numericalexp-real}
	In this experiment we show that, in principle, our approach of reduction of Problem~\ref{pr:problem-pw} to Problem~\ref{pr:problem-rw} can be used in a realistic SPECT setting.
	\par A conventional SPECT procedure was performed on a monkey. The data for this experiment was provided by \textit{Service Hospitalier Fr\'{e}d\'{e}ric Joliot, CEA (Orsay)}. The provided data consisted of two files: first one contained the three-dimensional attenuation map of monkey's head and the second contained emission data in terms of photon counts $N(\gamma)$ along rays $\gamma, \, \gamma\in \Gamma$, where $\Gamma$ was given by  \eqref{eq:numerical.simulation-noiseless.gamma1} for $n_z=n_s=n_{\varphi}=128$. The provided attenuation map $a=a(x), \, x\in \R^3$, was given in a form of a volumetric image of $128\times128\times128$ pixels. 
	\begin{figure}[H]
		\centering
		\begin{subfigure}{0.24\textwidth}
			\includegraphics[width=45mm, height=45mm]{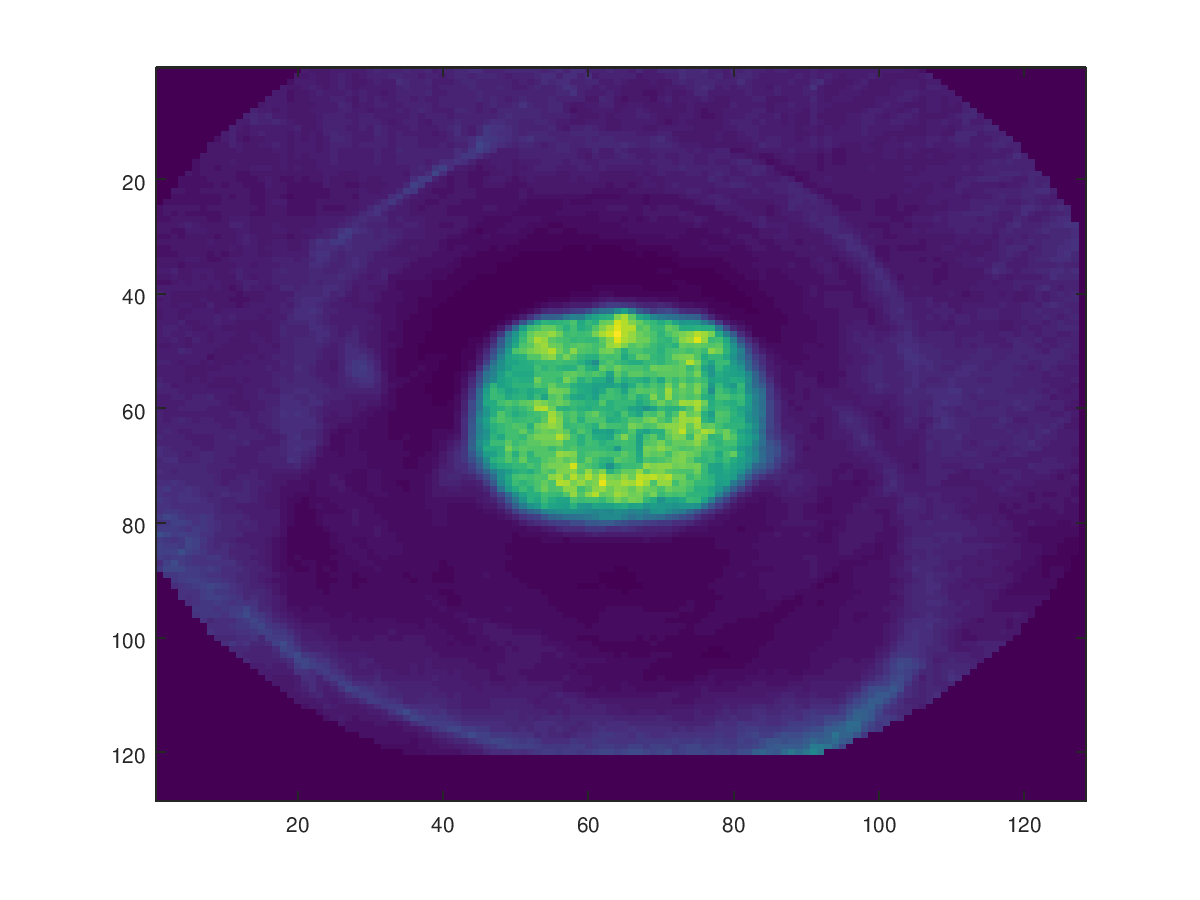}
			\caption{$i_z=50$}
			
		\end{subfigure}
		\begin{subfigure}{0.24\textwidth}
			\includegraphics[width=45mm, height=45mm]
			{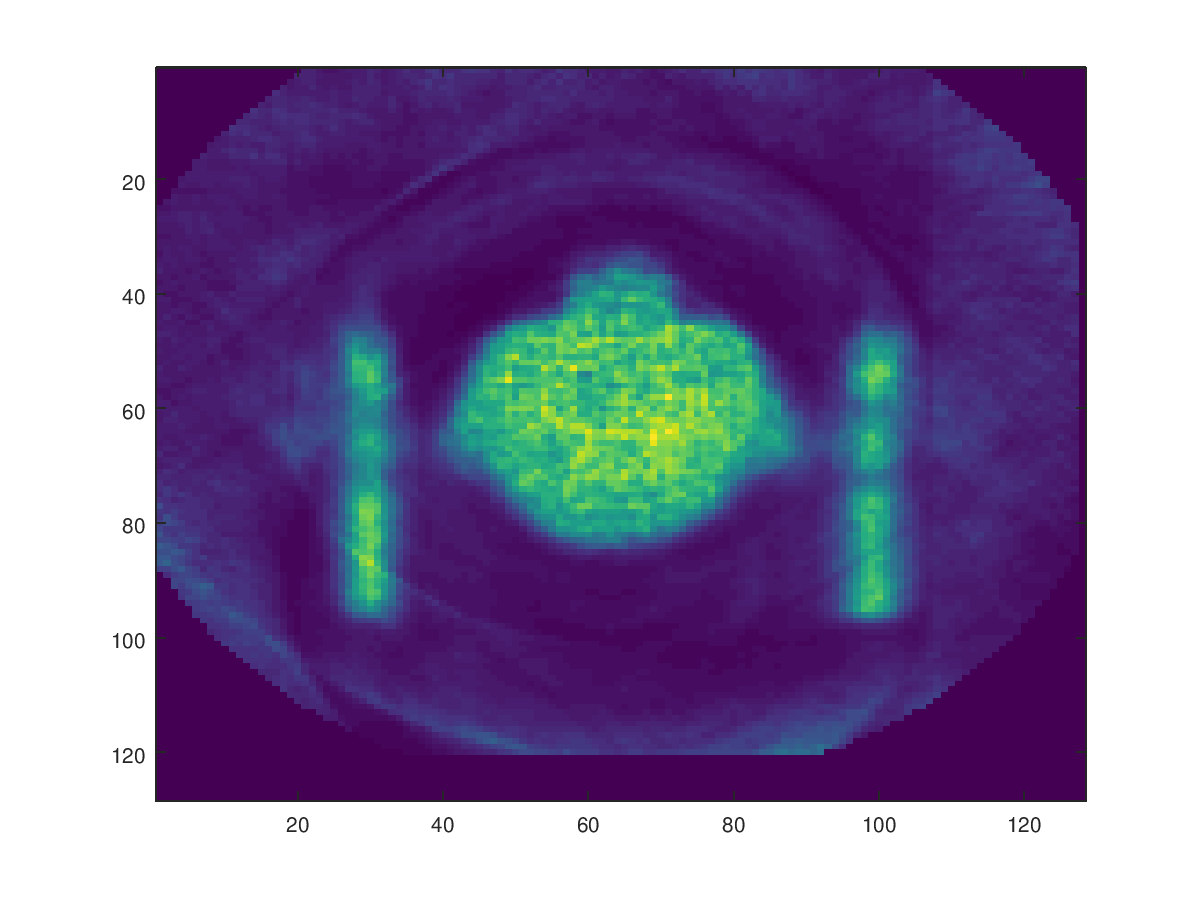}
			\caption{$i_z=60$}
			
		\end{subfigure}
		\begin{subfigure}{0.24\textwidth}
			\includegraphics[width=45mm, height=45mm]{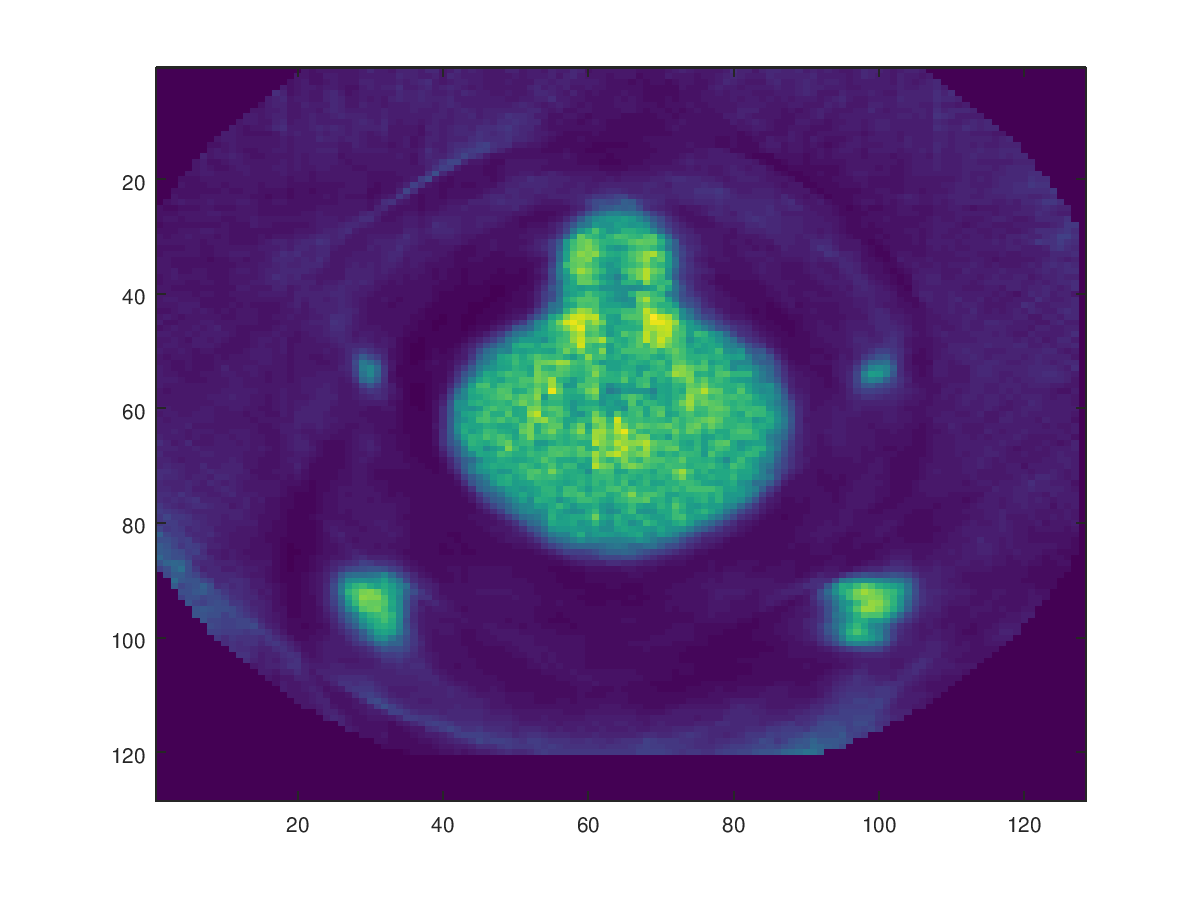}
			\caption{$i_z=70$}
			
		\end{subfigure}
		\begin{subfigure}{0.24\textwidth}
			\includegraphics[width=45mm, height=45mm]
			{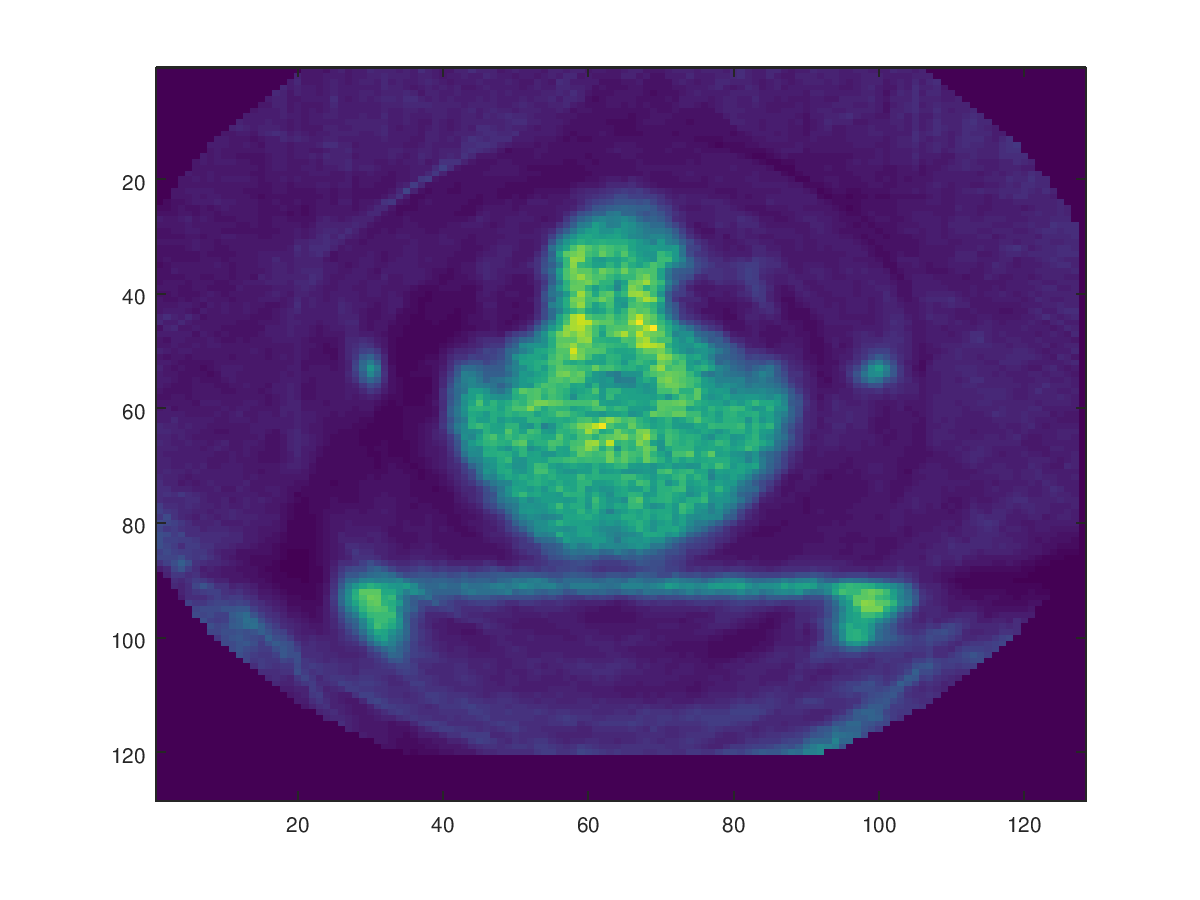}
			\caption{$i_z=80$}
			
		\end{subfigure}
		\caption{ attenuation map of monkey's head; images (A), (B), (C), (D) represent the attenuation map reduced to sequence of planes $z=const$ (smaller values for $i_z$ correspond to the upper part of the monkey's head).
			Green region in the center of image (A) corresponds to brain material, vertical and horizontal lines on images (B), (D) correspond to plates which were used to fix the head of the monkey.}
		\label{fig:monkey.data}
	\end{figure}
	Unfortunately, in the given data the units for the attenuation map were not provided, which is crucial due to non-linear dependence of transform $P_{W_a}$ on the attenuation. To overcome this lack of information we normalized the attenuation map so that the brain material in Figure~\ref{fig:monkey.data} (A) corresponded to attenuation of water $0.15\text{cm}^{-1}$. Also, constant $C$ from \eqref{eq:spect.statistics.params} was not provided, which made it possible to reconstruct the isotope distribution only up to a multiplicative constant. The acquisition per projection was $14$ seconds and the radius of rotation of detectors was 261mm. The total number of registered photons was approximately  $4.6\cdot 10^{6}$ and the maximal number of registered photons per one detector was $60$.
	\begin{figure}[H]
		\centering
		\begin{subfigure}{0.24\textwidth}
			\includegraphics[width=45mm, height=45mm]{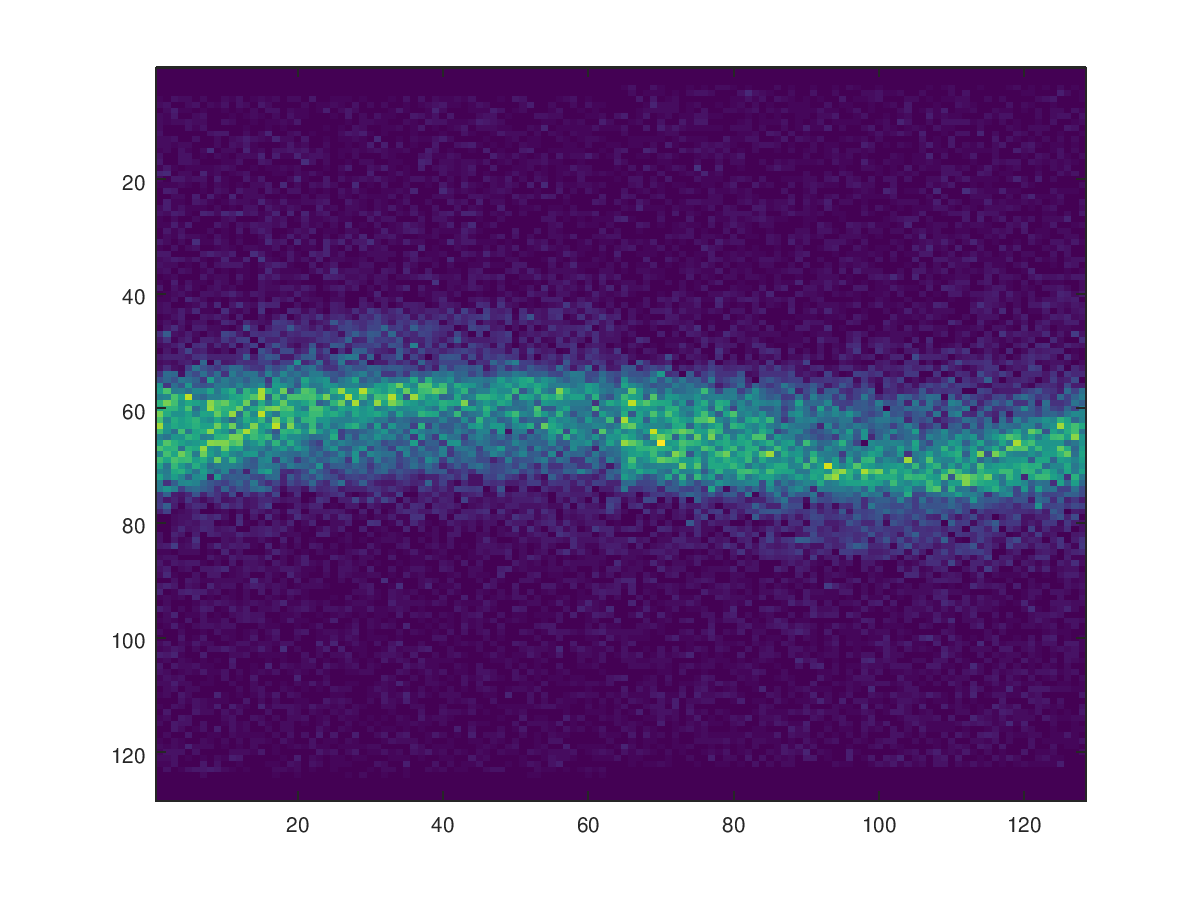}
			\caption{$i_z = 50$}
			
		\end{subfigure}
		\begin{subfigure}{0.24\textwidth}
			\includegraphics[width=45mm, height=45mm]
			{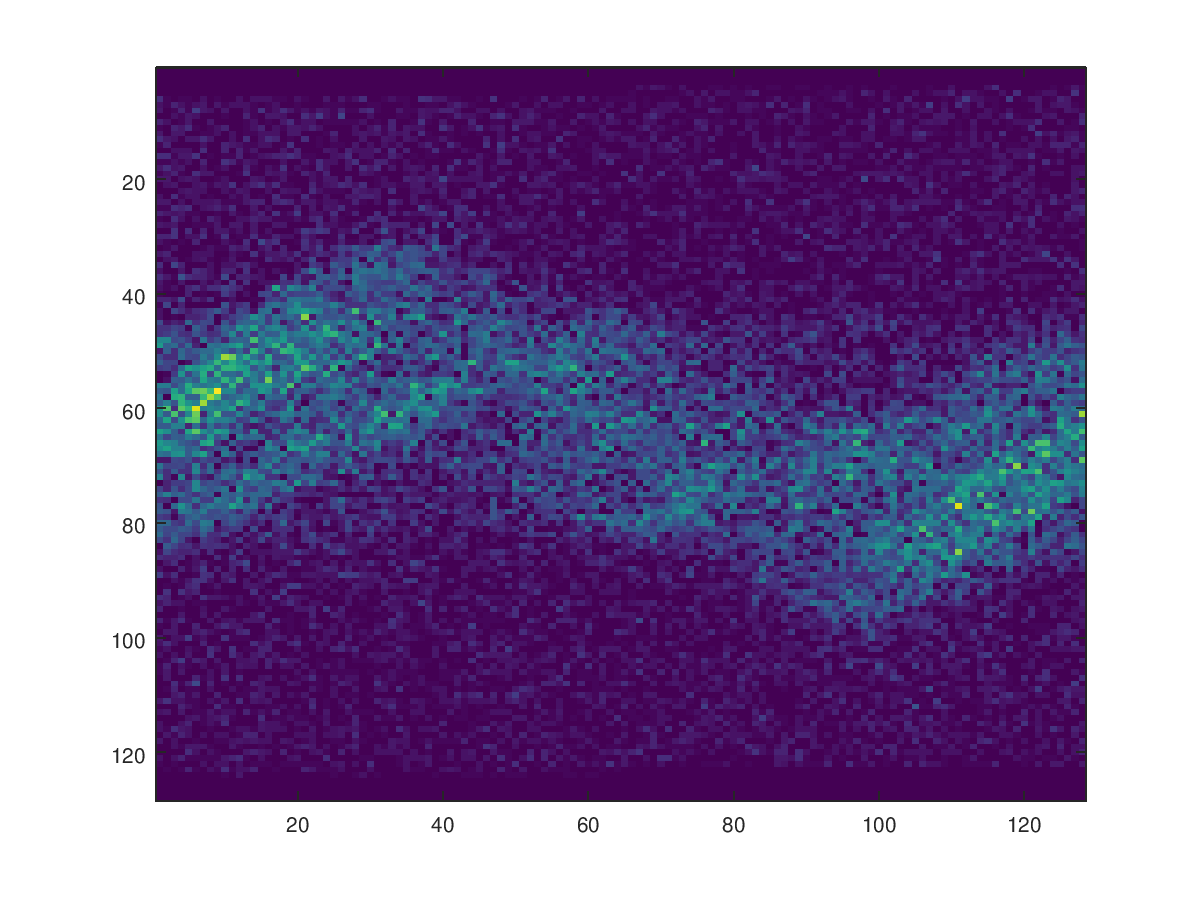}
			\caption{$i_z=60$}
			
		\end{subfigure}
		\begin{subfigure}{0.24\textwidth}
			\includegraphics[width=45mm, height=45mm]{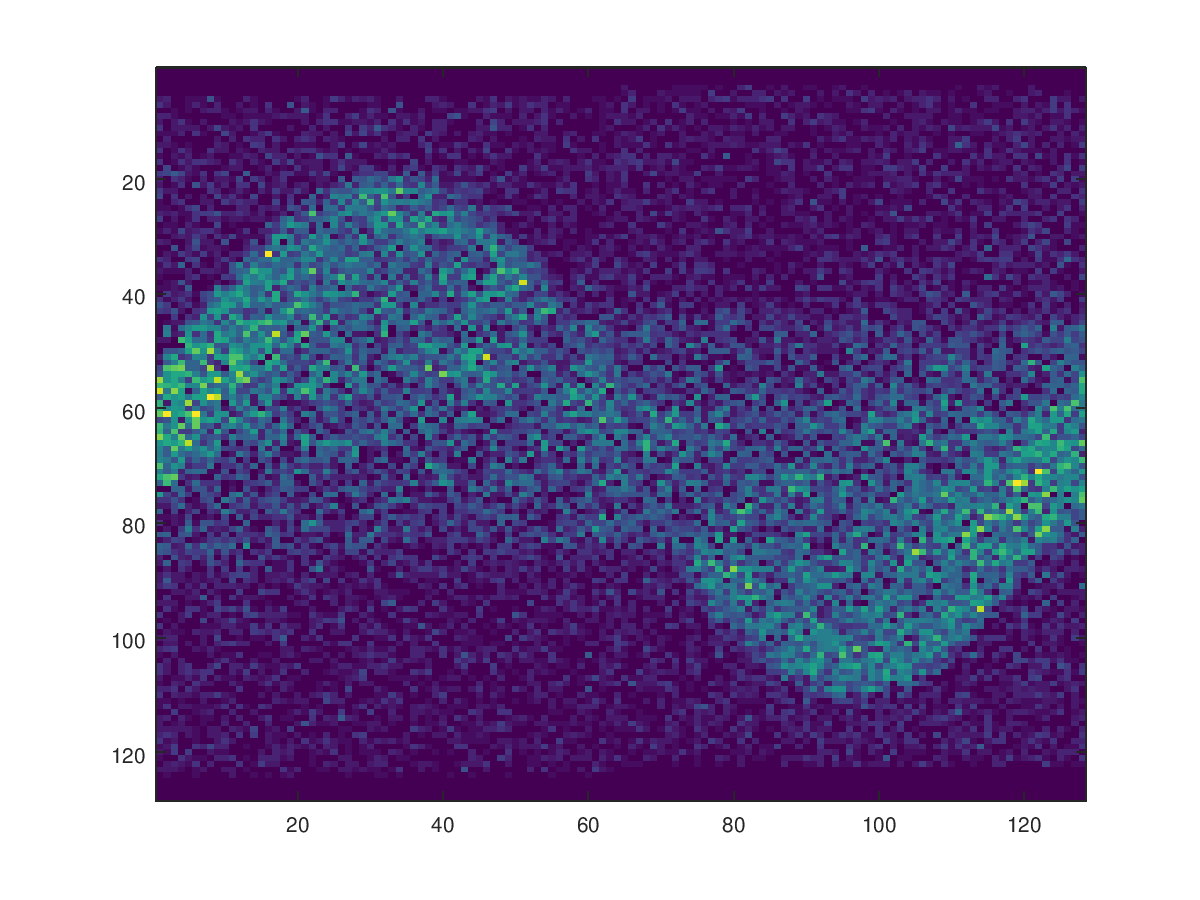}
			\caption{$i_z=70$}
			
		\end{subfigure}
		\begin{subfigure}{0.24\textwidth}
			\includegraphics[width=45mm, height=45mm]
			{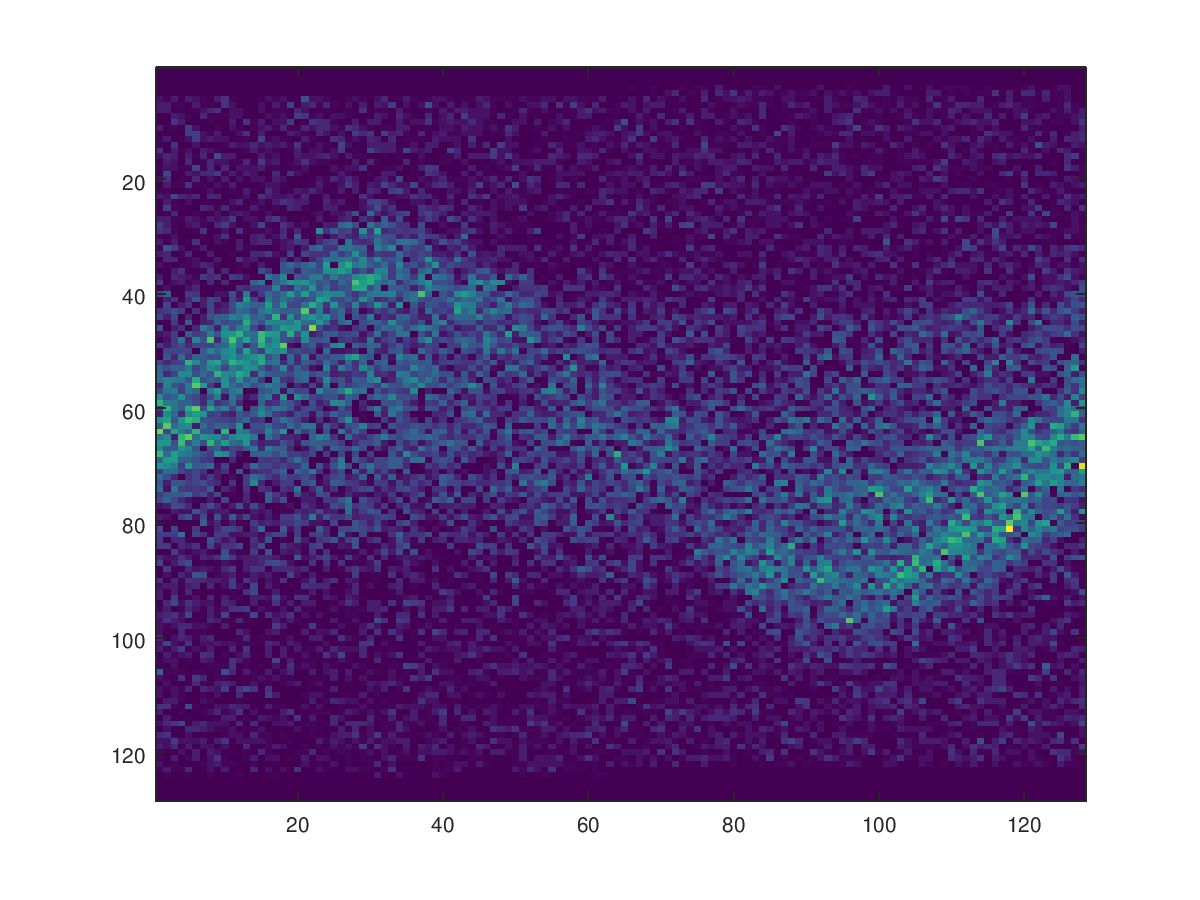}
			\caption{$i_z=80$}
			
		\end{subfigure}
		\caption{ emission data in the experiment on a monkey; images (A)-(D) 
			represent the photon counts $N(\gamma)$, where rays $\gamma$ belong to different slicing planes $z=const$. In each plane $z=const$ rays are parametrized by $(s,\varphi), \, s\in [-1,1], \, \varphi \in [0, 2\pi]$ (see also formula~\eqref{eq:numerical.simulation-noiseless.raygamma}); horizontal axis on (a)-(d) corresponds to variable $\varphi$, vertical axis corresponds to variable $s$.}
		\label{fig:monkey.emission.data}
	\end{figure}
	Applying the preprocessing procedure and Chang-type method for dimension $d=3$ (see Subsection~\ref{subsect:reconstruct.methods}), we obtained approximate reconstructions of the isotope distribution in the brain of the monkey.
	\begin{figure}[H]
		\centering
		\begin{subfigure}{0.24\textwidth}
			\includegraphics[width=45mm, height=45mm]{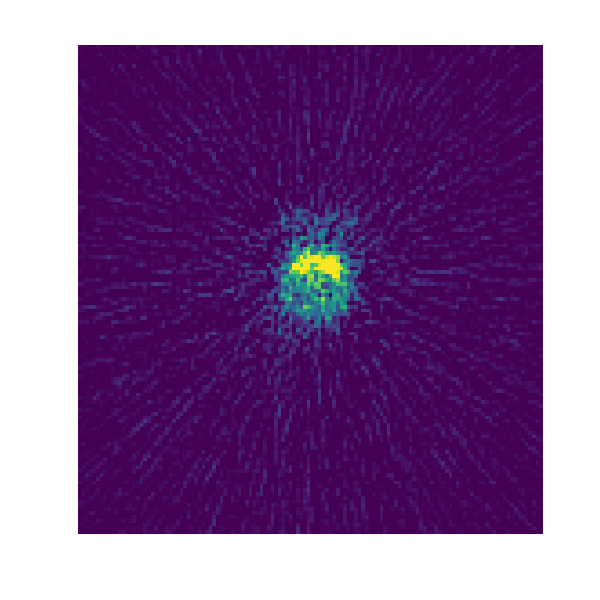}
			\caption{$i_z = 50$}
			
		\end{subfigure}
		\begin{subfigure}{0.24\textwidth}
			\includegraphics[width=45mm, height=45mm]
			{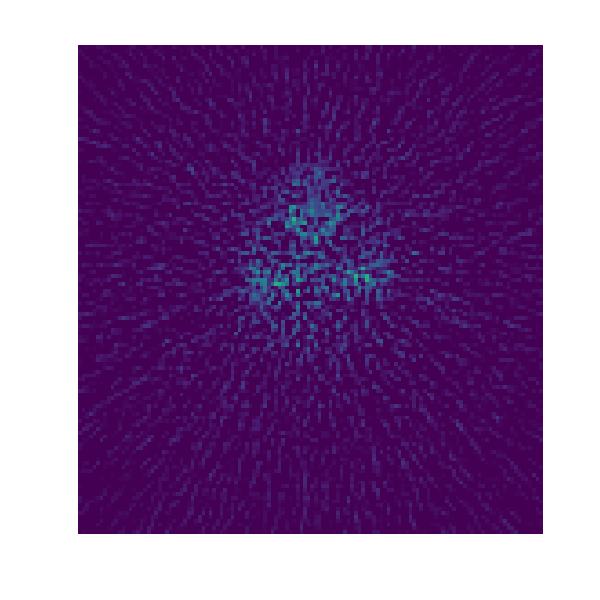}
			\caption{$i_z=60$}
			
		\end{subfigure}
		\begin{subfigure}{0.24\textwidth}
			\includegraphics[width=45mm, height=45mm]{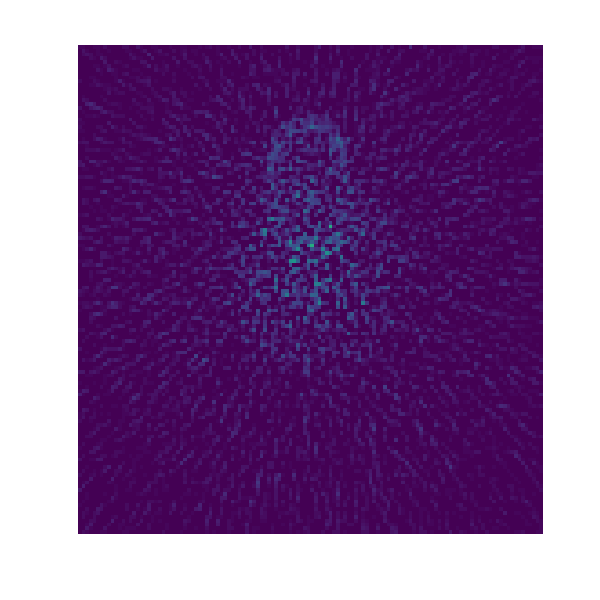}
			\caption{$i_z=70$}
			
		\end{subfigure}
		\begin{subfigure}{0.24\textwidth}
			\includegraphics[width=45mm, height=45mm]
			{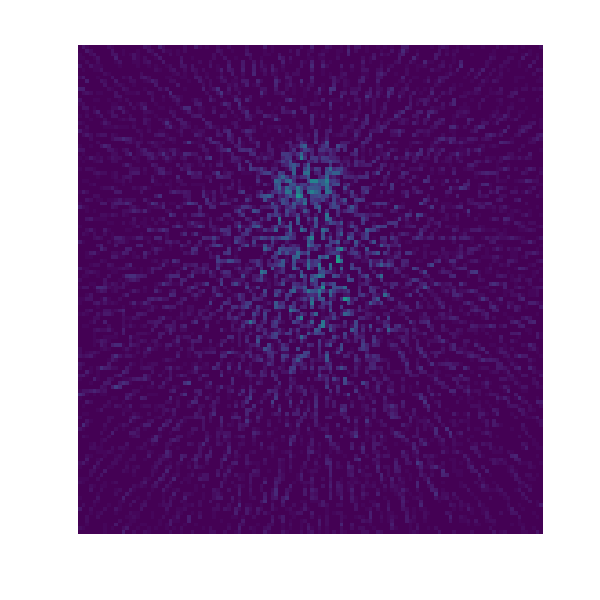}
			\caption{$i_z=80$}
			
		\end{subfigure}
		\caption{ reconstruction for the isotope distribution in different slicing planes $z=const$; note that slicing plane $i_z=50$ in Figure~\ref{fig:monkey.data} (A) corresponds to the area of the frontal lobes.}
		\label{fig:monkey.emission.reconstruction}
	\end{figure}
	In the center of Figure~\ref{fig:monkey.emission.reconstruction}(A) one can notice high concentration of the isotope. This area geometrically corresponds to the area area with frontal lobes (see also Figure~\ref{fig:monkey.data} (a)), where the concentration of the isotope is expected. This supports the idea that such preprocessing can be applied in a real SPECT procedure.

	\section{Discussion}\label{sect:discussion}
	Presented results, especially those in Tables~\ref{tab:errors.chang},~\ref{tab:errors.iterative}, demonstrate that 
	proposed preprocessing procedure can significantly decrease the effect of noise on reconstructions in SPECT. Specifically, we demonstrate this by using reconstruction methods of similar types in two and three dimensions and what is indeed important is that they produce the same reconstructions from raw and preprocessed data (i. e., from $P_{W_a}f$ and $R_wf$, respectively) in the noiseless case. Here, the use of approximate reconstruction methods (i.e. Chang-type, Kunyansky-type methods) can be seen as a first step in analyzing   the proposed procedure. Further numerical tests using reconstruction procedures  better adapted for correcting strong noise in $P_Wf$ and in $R_wf$ (e.g.,  maximum likelihood or bayesian methods) are of definite interest and will be given elsewhere.
	
	One could notice that we do not give any theoretical analysis why the proposed procedure brings some noise regularization. This is also a topic for future research, however, we believe that there are at least two reasons for this effect. 
	
	The first one is due to higher signal-to-noise ratio in data given by $R_wf$ than in one given by $P_{W_a}f$. Indeed, from reduction formulas of  \eqref{eq:reduction.main-formula} and formulas \eqref{eq:spect.statistics}, \eqref{eq:spect.statistics.params} one can see that data modeled $R_{w}f$ has also Poisson distribution but it has higher intensities than in data modeled by $P_{W_a}f$. We recall that for Poisson distribution signal-to-noise ratio is proportional to $\sqrt{\lambda}$, where $\lambda$ is the intensity. So one can argue that higher signal-to-noise ratio reflects the improvement in reconstructions in presence of noise. At the same time the reduction of  Problem~\ref{pr:problem-pw} to Problem~\ref{pr:problem-rw} brings more numerical instability in reconstructions. Indeed, on each slice Problem~\ref{pr:problem-pw} is less ill-posed (ill-posedeness degree is at least $1/2$) than Problem~\ref{pr:problem-rw} in the whole space (ill-posedeness degree is at least $1$). So there are two opposite effects of the reduction which require further investigation. 
	
	The second possible reason is due to filtering of the reconstructions in $z$-direction. Indeed, in slice-by-slice reconstructions only the data given by $P_Wf$ for rays reduced to planes $z=const$ are used. At the same time, in formulas of   \eqref{eq:reduction.main-formula} one can see that for each plane $(s,\theta)\in \Pi$ value for $R_wf(s,\theta)$ is composed from values of $P_{W}f$ for rays varying in $z$-direction. This can be seen as an application of some non-linear filter on reconstructions in $z$-direction. However, in this case our preprocessing procedure has an advantage over standard filtering schemes because it is based on a geometric identity and has no free parameters to optimize.
	
	Finally, from a statistical point of view the reduction by formulas of  \eqref{eq:reduction.main-formula} introduces correlations between variables $R_wf(s,\theta)$ for pairs of intersecting planes $(s,\theta) \in \Pi$. Therefore, prior to apply some advanced reconstruction algorithms (for example, EM or bayesian algorithms) to data given by $R_wf$ one has to decide how to consider these new correlations. This also raises new questions for future research which will be given elsewhere.
	
	\section{Appendix: a few remarks on implementations}
	\label{sect.sect.numerical.implementation}
	There were only two crucial steps for our numerical implementations\footnote{Some details of our implementations can be found at GitHub repository:
		\noindent\href{https://github.com/fedor-goncharov/wrt-project}{github.com/fedor-goncharov/wrt-project}}:
	
	\begin{enumerate}
		\item \underline{Numerical implementation of formulas \eqref{eq:numerical.reduction.formulas-discr}, \eqref{eq:numerical.weight.formulas-discr} for reduction of $P_{W_a}f$ to $R_{w} f$}. This issue was already commented in Subsection~\ref{subsect:preprocessing.discretizations}.
		We used quadratic and spline interpolations in variables $z,s$, respectively, to sample $P_{W_a}f$ for missing rays. By doing so we were  obtaining the data which did not belong to the image of operator $P_{W_a}$. To our knowledge, efficient interpolation of data even for $Pf$, where $P$ is the classical ray transform, is still an open question.
		
		\item \underline{Inversion of classical Radon transform $R$ for $d=3$}. The reconstruction methods from Subsection~\ref{subsect:reconstruct.methods} are based on inversion of the classical Radon transform in dimension $d=3$. There exist many open-access libraries for efficient computations of $R^{-1}$ for $d=2$ (for example, in MATLAB/Octave, C/C++ or Python), but for $d=3$ we did not find any. Because of that we have implemented our own numerical version of $R^{-1}$ for $d=3$, based on the well-known Projection Theorem (see \cite{natterer1986book}, Chapter 2, Theorem 1.1) and using a very nice NUFFT library (i.e., Non-uniform Fast Fourier Transform) for MATLAB/Octave developed in TU Chemnitz~\cite{keiner2009nfft}.
	\end{enumerate}

	\section*{Acknowledgements}
	The materials for the present work were obtained in the framework of research conducted under the direction of professor R. G. Novikov. I am also very grateful for the comments and remarks of professor L.~Kunyansky during the conference AIP-2019 and also during the author's thesis defense procedure.
	Also I would like to express my gratitude to \mbox{J.-P.~Guillement} from Universit\'{e} de Nantes for providing codes for SPECT reconstructions in 2D and for related discussions on numerical methods in SPECT.


\begin{thebibliography}{99}
		
		\bibitem{bal2011spect}
		G.~Bal, A. Jollivet,  
		\newblock{Combined source and attenuation reconstructions in SPECT}.
		\newblock{Tomography and Inverse Transport Theory, Contemp. Math},  559, 13-27, 2011.
		
		\bibitem{bocoum2019structured}
		M.~Bocoum, J.-L.~Gennisson, J.-B.~Laudereau, A.~Louchet-Chauvet, J.-M.~ Tualle, F.~Ramaz. 
		\newblock{Structured ultrasound-modulated optical tomography},
		\newblock{\em Applied optics} 58(8): 1933-1940, 2019.
		
		\bibitem{beylkin1984inversion}
		G.~Beylkin,
		\newblock{The inversion problem and applications of the generalized Radon transform}.
		\newblock{\em Communications on pure and applied mathematics}, 37(5): 579-599, 1984.
		
		\bibitem{boman2011localnoninj}
		J.~Boman, 
		\newblock{Local non-injectivity for weighted Radon transforms}.
		\newblock{\em Contemp. Math} 559: 39-47, 2011.
		
		
		\bibitem{chang1978correciton}
		L.-T.~Chang, 
		\newblock{A method for attenuation correction in radionuclide computed tomography}. 
		\newblock{\em IEEE Transactions on Nuclear Science}, 25(1): 638-643, 1978.
		
		\bibitem{deans2007book}
		S.~R.~Deans,
		\newblock{The Radon transform and some of its applications}, 
		\newblock{Courier Corporation}, 2007.
		
		\bibitem{filipovic2019pet}
		M.~Filipović, E.~Barat, T.~Dautremer, C.~Comtat, S.~Stute. 
		\newblock{PET reconstruction of the posterior image probability, including multimodal images.} 
		\newblock{\em IEEE transactions on medical imaging}, 2018.
		
		
		\bibitem{frigyik2008xray}
		B.~Frigyik, P.~Stefanov, G.~Uhlmann,
		\newblock{The X-ray transform for a generic family of curves and weights},
		\newblock{\em Journal of Geometric Analysis}, 18(1), 2008.
		
		
		\bibitem{gach2008shepplogan}
		H.~M.~Gach, C.~Tanase, F.~Boada. 
		\newblock{2D \& 3D Shepp-Logan phantom standards for MRI}. 
		\newblock{\em 19th International Conference on Systems Engineering}, pp. 521-526, 2008.
		
		\bibitem{ggg1959geomhom} 
		I. M~Gelfand, M.~I.~Graev. 
		Geometry of homogeneous spaces, representations of groups in
		homogeneous spaces and related questions of integral geometry.
		\newblock{\em Trudy Moskov. Mat. Obshch.}, 8:321–390, 1959.
		
		
		\bibitem{ggg2003selected}
		I.~M.~Gelfand, S.~Gindikin, M.~Graev,
		\newblock{Selected topics in integral geometry}.
		\newblock{\em American Mathematical Soc.}, Vol. 220, 2003
		
		\bibitem{ggv1962generalized}
		I.~M.~Gelfand, M.~I.~Graev, N.~J.~Vilenkin.
		\newblock{Obobshchennye funktsii, Vyp. 5. Integralnaya
			geometriya i svyazannye s nei voprosy teorii predstavlenii},
		\newblock{\em Gosudarstv. Izdat. Fiz.-Mat. Lit.}, 1962.
		
		\bibitem{goncharov2016analog}
		F.~O.~Goncharov, R.~G.~Novikov, 
		\newblock{An analog of Chang inversion formula for weighted Radon transforms in multidimensions}.
		\newblock{\em Eurasian Journal of Mathematical and Computer Applications}, 4(2): 23-32, 2016.
		
		\bibitem{goncharov2017iterative}
		F.~O.~Goncharov,
		\newblock{An iterative inversion of weighted Radon transforms along hyperplanes}.
		\newblock{\em Inverse Problems}, 33(12): 124005, 2017.
		
		\bibitem{guillement2002spect}
		J.-P.~Guillement, F.~Jauberteau, L.~Kunyansky, R.~Novikov, R.~Trebossen. 
		\newblock{On single-photon emission computed tomography imaging based on an exact formula for the nonuniform attenuation correction}, 
		\newblock{\em Inverse Problems 18, no. 6}, 2002.
		
		\bibitem{guillement2008wiener}
		J.-P.~Guillement, R.~G. Novikov, 
		\newblock{On Wiener type filters in SPECT.}
		\newblock{\em Inverse Problems}, 24(2), 2008
		
		\bibitem{guillement2014finite}
		J.-P.~Guillement, R.~G.~Novikov, 
		\newblock{Inversion of weighted Radon transforms via finite Fourier series weight approximations}. 
		\newblock{\em Inverse Problems in Science and Engineering},
		22(5): 787-802, 2014.
		
		
		\bibitem{hudson1994accelerated}
		H. M.~Hudson, R.~S.Larkin,
		\newblock{Accelerated image reconstruction using ordered subsets of projection data.} 
		\newblock{\em IEEE transactions on medical imaging} 13(4):601-609, 1994.
		
		\bibitem{harishchandra1958}
		Harish-Chandra,
		\newblock {Spherical functions on a semi-simple Lie group I}.
		\newblock{\em American Journal of Mathematics}, 80:241–310, 1958.
		
		\bibitem{harishchandra1958a} 
		Harish-Chandra,
		\newblock{Spherical functions on a semi-simple Lie group II}.
		\newblock{\em American Journal of Mathematics}, 80: 553–613, 1958.
		
		\bibitem{helgason1965}
		S.~Helgason.
		\newblock{Radon-Fourier transforms on symmetric spaces and related group representations}.
		\newblock{\em Bull. Amer. Math. Soc.}, 71: 757–763, 1965.
		
		\bibitem{helgason1965a}
		S.~Helgason,
		\newblock{The Radon transform on Euclidean spaces, compact two-point homogeneous
			spaces and Grassman manifolds}.
		\newblock{\em Acta Math.}, 113:153–180, 1965.
		
		\bibitem{helgason1966}
		S.~Helgason.
		\newblock{A duality in integral geometry on symmetric spaces}.
		\newblock{\em In Proc. U.S.-Japan Seminar in Differential Geometry, Kyoto}, 1965.
		
		\bibitem{helgason70}
		S.~Helgason.
		\newblock{A duality for symmetric spaces, with applications to group representations.}
		\newblock{\em Advances in Math.}, 5:1–154, 1970.
		
		\bibitem{helgason1999}
		S.~Helgason,
		\newblock {The Radon transform}.
		\newblock{Vol.2. Boston: Birkh\"{a}user}, 1999.
		
		
		\bibitem{ilmavirta2016}
		J.~Ilmavirta, 
		\newblock{On Radon transforms on compact Lie groups}.
		\newblock{\em Proceedings of the American Mathematical Society}, 144(2): 681-691, 2016.
		
		\bibitem{ilmavirta2019}
		J.~Ilmavirta, J.~Railo,
		\newblock{Geodesic ray transform with matrix weights for piecewise constant functions}
		\newblock{\em arXiv preprint}, arXiv:1901.030525, 2019. 
		
		
		\bibitem{john1955book}
		F.~John,
		\newblock {Plane Waves and Spherical Means Applied to Partial Differential Equations}.
		\newblock {Dover Publications}, 1955.
		
		
		\bibitem{keiner2009nfft}
		J.~Keiner, S.~Kunis, D.~Potts, 
		\newblock{Using NFFT 3 - a software library for various nonequispaced fast Fourier transforms}.
		\newblock{\em ACM Trans. Math. Software}, 36, Article 19, 1-30, 2009.
		
		\bibitem{kunyansky1992generalized}
		Kunyansky, L.,
		\newblock{Generalized and attenuated {R}adon transforms: restorative approach to the numerical inversion}.
		\newblock {\em Inverse Problems}, 8(5): 809-819, 1992.
		
		\bibitem{kunyansky2001formula}
		L.~Kunyansky, 
		\newblock{A new SPECT reconstruction algorithm based on Novikov explicit inversion formula}.
		\newblock{\em Inverse Problems}, 17(2), 2001.
		
		\bibitem{kuchment2014book}
		P.~Kuchement, 
		\newblock{The Radon transform and medical imaging},
		\newblock{\em SIAM}, Vol. 85., 2014.
		
		\bibitem{miqueles2011}
		E.~Miqueles, A.~R.~De Pierro, 
		\newblock{Iterative reconstruction in X-ray fluorescence tomography based on Radon transform}.
		\newblock{\em IEEE Transactions on medical imaging}, 30(2): 438-450, 2011.
		
		
		\bibitem{miller1987newslant}
		D.~Miller, M.~Oristaglio, G.~Beylkin. 
		\newblock{A new slant on seismic imaging: Migration and integral geometry.}, 
		\newblock{\em Geophysics 52, no. 7: 943-964}, 1987.
		
		\bibitem{natterer1986book}
		F.~Natterer, 
		\newblock{The mathematics of computerized tomography}, 
		\newblock{\em SIAM}, Vol. 32, 1986.
		
		\bibitem{natterer2001inversion}
		F.~Natterer, 
		\newblock{Inversion of the attenuated Radon transform}, 
		\newblock{\em Inverse Problems}, 17(1), 2001.
		
		\bibitem{nguyen2009compton}
		M.~K.~Nguyen, T.T.~Truong, D.~Driol, H.~Zaidi, 
		\newblock{On a novel approach to Compton scattered emission imaging}, \newblock{\em IEEE Transactions on Nuclear Science}, 56(3): 1430-1437, 2009.
		
		
		\bibitem{novikov2002formula}
		R.~G.~Novikov, 
		\newblock{An inversion formula for the attenuated X-ray transform}.
		\newblock{\em Arkiv f\"{o}r Matematik}, 40(1): 145-167, 2002.
		
		
		\bibitem{novikov2011chang}
		R.~G.~Novikov,
		\newblock{Weighted Radon transforms for which Chang's approximate inversion formula is exact}. 
		\newblock{\em Russian Mathematical Surveys} 66(2), 2011.
		
		\bibitem{novikov2014weighted}
		R.~G. Novikov.
		\newblock {W}eighted {R}adon transforms and first order differential systems on
		the plane.
		\newblock {\em Moscow mathematical journal}, 14(4):807--823, 2014.
		
		
		\bibitem{paternain2016geodesic}
		G.~P.~Paternain, M.~Salo, G.~Uhlmann, Z.~Hanming
		\newblock{The geodesic X-ray transform with matrix weights},
		\newblock{\em arXiv preprint arXiv:1605.07894}, 2016.
		
		
		\bibitem{sharafutdinov2012}
		V.~A.~Sharafutdinov,
		\newblock{Integral geometry of tensor fields}, 
		\newblock{\em Walter de Gruyter}, Vol.~1, 2012.
		
		
		\bibitem{shepp1974fourier}
		L. Shepp, B. Logan, 
		\newblock{The Fourier reconstruction of a head section}.
		\newblock{\em IEEE Transactions on nuclear science}, 21(3): 21-43, 1974.
		
		\bibitem{vslambrock2015pet}
		K.~V.~Slambrouck, S.~Stute, C.~Comtat, M.~Sibomana, F.H.~van Velden, R.~Boellaard, J.~Nuyts,
		\newblock{Bias reduction for low-statistics PET: maximum likelihood reconstruction with a modified Poisson distribution}.
		\newblock{\em IEEE Transactions on medical imaging},
		34(1): 126-136, 2015.
		
		
		\bibitem{slambrouck2014reconstruction}
		K.~V.~Slambrouck, J.~Nuyts, 
		\newblock{Reconstruction scheme for accelerated maximum likelihood reconstruction: The patchwork structure}, 
		\newblock{\em IEEE Transactions on Nuclear Science}, 61(1):173-181, 2014.
		
		\bibitem{strichartz1982variations}
		R.~S.~Strichartz,
		\newblock{Radon inversion–variations on a theme}.
		\newblock{\em The American Mathematical Monthly}, 89(6): 377-423, 1982.
		
		\bibitem{strichartz1991harmonic}
		R.~S.~Strichartz, 
		\newblock{$L^p$ harmonic analysis and Radon transforms on the Heisenberg group}.
		\newblock{\em Journal of functional Analysis}, 96(2): 350-406, 1991.
		
		\bibitem{quinto1980measures}
		E.~T. Quinto, 
		\newblock{The dependence of the generalized Radon transforms on defining measures}.
		\newblock{\em Transactions on the American Mathematical Society}, 257(2): 331-346, 1980.
		
		\bibitem{quinto1983rotinv}
		E.~T.~Quinto,
		\newblock{The invertibility of rotation invariant Radon transforms}.
		\newblock{\em Journal of Mathematical Analysis and Applications}, 91(2): 510-522, 1983.
		
		\bibitem{quinto2006}
		E.~T. Quinto, 
		\newblock{An introduction to X-ray tomography and Radon transforms}.
		\newblock{\em Proceedings of symposia in Applied Mathematics, Vol. 63}, 2006.
		
		\bibitem{quinto2014}
		E.~T. Quinto, G.~Ambartsoumian, R.~Felea, V.~Krishnan, C.~Nolan, 
		\newblock{Microlocal Analysis and Imaging, a short note in ``The mathematics of the planet Earth''}.
		\newblock{Chapter 7, pages 8-11}, Springer, Berlin, New York, 2014.
		
		\bibitem{quinto2018}
		E.~T.~Quinto, C.~Grathwohl, P.~Kunstmann, A.~Rieder,
		\newblock{Microlocal Analysis of imaging operators for effective common offset seismic reconstruction}.
		\newblock{\em Inverse Problems}, 34(11), 2018.
		
		
	\end{thebibliography}
\end{document}